\documentclass[preprint]{article}

\usepackage{color}

\RequirePackage{amsthm,amsmath,natbib}
\RequirePackage[colorlinks]{hyperref}

\usepackage{graphicx}
\usepackage[utf8]{inputenc}
\usepackage[T1]{fontenc}
\usepackage{aeguill}
\usepackage{amssymb}
\usepackage[utf8]{inputenc}
\usepackage[T1]{fontenc}
\usepackage{aeguill}
\usepackage{amssymb}
\usepackage{graphicx}
\usepackage[a4paper,pdftex=true]{geometry}

\usepackage{natbib}
\usepackage{bbm}
\usepackage{amsmath}
\usepackage{amssymb}
\usepackage{hyperref}
\usepackage[a4paper,pdftex=true]{geometry}
\usepackage{fancyvrb}

\newtheorem{theorem}{Theorem}
\newtheorem{lemma}[theorem]{Lemma}
\newtheorem{corollary}[theorem]{Corollary}

\newtheorem{proposition}[theorem]{Proposition}

\newtheorem{remark}{Remark}

\newtheorem{definition}{Definition}

\newcommand{\RR}[1][2]{\mathbb{R}^{#1}}
\newcommand{\ZZ}[1][2]{\mathbb{Z}^{#1}}
\newcommand{\NN}[1][]{\mathbb{N}^{#1}}

\newcommand{\Vect}[1]{
  \mathbf{#1}
}

\newcommand{\tr}[1]{
 {#1}^{\!  T}
}

\newcommand{\Mat}[1]{
  \underline{\mathbf{#1}}
}

\newcommand{\Esp}{
  \mathbf{E}
}

\newcommand{\Var}{
  Var
}

\newcommand{\mSp}{\mathbbm{M}}

\newcommand{\mm}{\lambda^{\mE}}
\newcommand{\mE}{\mathbbm{m}}

\newcommand{\cSp}[1][]{{\Omega}#1}

\newcommand{\cSptt}{\overline{\Omega}}

\newcommand{\sSp}{\mathbbm{S}}
\newcommand{\sB}{\mathcal{B}}
\newcommand{\sm}{\mu}

\newcommand{\xm}[2][x]{#1^{#2}}

\newcommand{\ism}[2][\Lambda\times\mSp]{\int_{#1} #2 \sm(d\xm{m})}

\newcommand{\cD}[1][]{\left\lceil{#1}\right\rceil}

\newcommand{\VE}{V}

\newcommand{\Par}{\theta}                    
\newcommand{\ParV}{\Vect{\theta}}            

\newcommand{\ParVT}{\Vect{\theta}^\star}     
\newcommand{\SpPar}{\Vect{\Theta}}           

\newcommand{\Dom}[1][0]{\Delta_{#1}}

\newcommand{\VIPar}[3]{\VE \left(#1|#2  ; #3 \right)}

\newcommand{\SEx}[2]{v_{#1}(#2)}
\newcommand{\SExI}[3]{v_{#1}(#2|#3)}

\newcommand{\SExIV}[2]{\Vect{v} \left(#1 | #2  \right)}

\newcommand{\cSEx}[2]{\kappa^{(#1)}_{#2}}

\newcommand{\Dtt}{\overline{D}}
\newcommand{\Ltt}{\overline{\Lambda}}

\newcommand{\rem}{\ensuremath{\mathcal{R}}}

\title{ Residuals and goodness-of-fit tests for stationary marked Gibbs point processes}

\author{Jean-François Coeurjolly$^{1,2}$ and Frédéric Lavancier$^3$ \\
$^1$ Gipsa-lab, CNRS, Grenoble, France \\
$^2$ Laboratory Jean Kuntzmann, Grenoble, France.\\
$^3$ Laboratoire Jean Leray, Nantes, France. }

\date{}
\begin{document}

\maketitle






\begin{center}
\begin{minipage}{10cm}{\small \centerline{\bf Abstract} 
The inspection of residuals is a fundamental step to investigate the quality of adjustment of a parametric model to data.
For spatial point processes, the concept of residuals has been recently proposed by \citeauthor{A-BadTurMolHaz05} \cite{A-BadTurMolHaz05} as an empirical counterpart of the {\it Campbell equilibrium} equation for marked Gibbs point processes. The present paper focuses on stationary marked Gibbs point processes and deals with asymptotic properties of residuals for such processes. In particular, the consistency and the asymptotic normality are obtained for a wide class of residuals including the classical ones (raw residuals, inverse residuals, Pearson residuals). Based on these asymptotic results, we define goodness-of-fit tests with Type-I error theoretically controlled. One of these tests constitutes an extension of the quadrat counting test widely used to test the null hypothesis of a homogeneous Poisson point process.
}\end{minipage}\end{center}




\begin{center}\begin{minipage}{14cm}{ \small

}\end{minipage}\end{center}

\section{Introduction}

Recent works on statistical methods for spatial point pattern makes parametric inference feasible for a wide range of models, see \cite{R-Mol08} for an overview of this topic and more generally the books of \cite{B-DalVer88}, \cite{B-StoKenMecRus87} \cite{B-MolWaa03} or \cite{B-IllPenSto08} for a survey on spatial point processes. The question is then to know whether the model is well-fitted to data or not. For classical parametric models, this is usually done via the inspection of residuals. They play a central role in parametric inference, see \cite{B-Atk85} for instance. This notion  is quite complex for spatial point processes and has been recently proposed by \citeauthor{A-BadTurMolHaz05} \cite{A-BadTurMolHaz05}, following ideas from a previous work of \citeauthor{A-StoGra91} \cite{A-StoGra91}.

 The definition of residuals for spatial point processes is a natural generalization of the well-known residuals for point processes in one-dimensional time, used in survival analysis (see \cite{B-FleHar91} or \cite{B-AndBorGilKei93} for an overview). For example, a simple measure of the adequacy of a one-dimensional point process model consists in computing the difference between the number of events in an interval $[0,t]$ and the conditional intensity (or hazard rate of the lifetime distribution) parametrically estimated and integrated from 0 to $t$. The extension in higher dimension requires further developments due to the lack of natural ordering. It may be done for point processes admitting a conditional density with respect to the Poisson process. These point processes correspond to the Gibbs measures. 
The equilibrium in one dimension  between the number of events and the integrated hazard rate  may be replaced in higher dimension by the {\it Campbell equilibrium} equation or {\it Georgii-Nguyen-Zessin} formula (see \cite{A-Geo76}, \cite{A-NguZes79b} and Section~\ref{sec-desc}), which is the basis for defining the class of $h-$residuals where $h$ represents a test function. In particular,  \citeauthor{A-BadTurMolHaz05} \cite{A-BadTurMolHaz05} consider different choices of $h$ leading to the so-called raw residuals, inverse residuals and Pearson residuals, and  show that they share similarities with the residuals obtained for generalized linear models.

Thanks to various diagnostic plots developped in the seminal paper \cite{A-BadTurMolHaz05} and implementation within the \texttt{R} package \texttt{spatstat} \cite{A-BadTur05}, residuals appear to be a very convenient tool in practice. Some properties of the residuals process are exhibited in \cite{A-BadTurMolHaz05} and \cite{A-BadMolPak08}, including a conditional independence property and variance formulae in particular cases. In these two papers, the authors conjecture that a strong law of large numbers and a central limit theorem should hold for the residuals as the sampling window expands.

Our paper addresses this question for $d-$dimensional stationary marked Gibbs point processes.  We obtain the strong consistency and the asymptotic normality in several contexts for a large class of test functions $h$. The $h-$residuals crucially depend on an estimate of the parameter vector. We consider the natural framework where the estimate is computed with the same data over which the $h-$residuals are assessed. The assumptions are very general and we show that they are fulfilled for several classical models, including  the area interaction point process, the multi-Strauss marked point process, the Strauss type disc process, the Geyer's triplet point process, etc. The assumptions on the estimator are quite natural and we show that they are fulfilled in particular by the maximum pseudolikelihood estimator (in short MPLE) (see \cite{A-BadTur00} for instance), for which asymptotic properties are now well-known (see \cite{A-JenMol91}, \cite{A-JenKun94}, \cite{A-BilCoeDro08}, \cite{A-DerLav09} and \cite{A-CoeDro09}).

Moreover, based on these asymptotic results, we propose statistical goodness-of-fit tests for which the Type-I error is asymptotically controlled. To the best of our knowledge, this is the first attempt in this direction. Such tests exist for rejecting the assumption of a homogeneous or inhomogeneous Poisson point process, but for general marked Gibbs point processes, the existing validation methods  are either graphical (for example by using the QQ-plot proposed by \cite{A-BadTurMolHaz05}) or rely on Monte-Carlo based simulations. We present two tests based on the computation of the residuals on different subdomains of the observation window. They extend in a very natural way the quadrat counting test for homogeneous Poisson distributions (see \cite{B-Dig03} for instance). Besides, we present a test which combines  several different $h-$residuals (associated to different functions $h$), computed on the entire observation window. The next step will be to implement these testing procedures to assess their power, compare them and reveal their limits. A thorough study will require extensive simulations and should deserve a separate paper.\\


The rest of the paper is organized as follows. Section~\ref{sec-prel} gathers the main notation used in this paper and briefly displays the general background. The definition of marked Gibbs point processes is given. They depend exclusively on the choice of an energy function or equivalently, of  a local energy function. All the assumptions are based on this function. The {\it Georgii-Nguyen-Zessin} formula is recalled, leading to the definition of the $h-$innovations and $h-$residuals. Some examples are presented including the classical residuals considered by \cite{A-BadTurMolHaz05} and  new ones connected to the well-known empty space function (or spherical contact distribution), denoted in the literature by $F$, see  \cite{B-MolWaa03} for instance.

Section~\ref{sec-asymp} deals with asymptotic properties and presents our main results. A parametric  stationary $d$-dimensional marked Gibbs point process is observed in a domain, denoted by $\Lambda_n$, assumed to increase up to $\RR[d]$. Sufficient conditions expressed in terms of the test function and the local energy function are given in order to derive the strong consistency result. We also propose an asymptotic control in probability of the departure of the $h$-residuals process from the $h$-innovations  through the departure of the estimate from the true parameter vector (see Proposition~\ref{prop-equiv}). This allows us to deduce asymptotic normality results. Two different frameworks are considered: for the first one, the initial domain is splitted into a fixed finite number of subdomains (with volumes aimed at converging to $+\infty$) and we consider the vector composed of the $h-$residuals computed on each subdomain. 
For the second framework, we consider the vector composed of the $h_j-$residuals (for $j=1,\ldots,s$) computed on the same domain $\Lambda_n$, where  $h_1,\ldots,h_s$ are different test functions.

The asymptotic normality results depend on unknown asymptotic covariance matrices. The important question of estimating these matrices is addressed in Section~\ref{sec-estMat}. We give a general condition under which these matrices are definite-positive and propose a consistent estimate.

Section~\ref{sec-gof} exploits the asymptotic results obtained before. Some goodness-of-fit tests are proposed, based on normalized residuals computed in the two previous frameworks. They are shown to converge to some $\chi^2$ distribution. Framework 1 leads to a generalization of the quadrat counting test for homogeneous Poisson distributions. Framework 2 yields a test which combines the information coming from several residuals, as for instance residuals coming from the estimation of the empty space function at several points.

The different assumptions made in the previous sections to obtain asymptotic results are discussed in Section~\ref{sec-disc}.  When considering classical test functions, exponential family models and the MPLE, the regularity and integrability type assumptions are shown to be satisfied for a wide class of examples. The testing procedures require moreover an identifiability condition to provide a proper normalization. 
Proposition~\ref{prop-lInn} shows that this condition is easy to check for the first proposed test and appears to be not restrictive in this case. For the other tests, checking this condition depends more specifically on the model and the test function.  We show,  in Proposition~\ref{prop-exemples}, how this condition can be verified on two examples of models with several choices of test functions.



In Section \ref{non-hereditary}, the very special situation where the energy function is not hereditary is considered. The GNZ formula is not valid any more in this setting but, provided a slight modification, it has been recently extended  in \cite{A-DerLav09}.  This leads to a natural generalization of the residuals to the non-hereditary setting.

Proofs of our main results are postponed to Section~\ref{sec-proofs}. The main material is composed of an ergodic theorem obtained by~\cite{A-NguZes79} and a new multivariate central limit theorem for spatial processes. Our setting actually involves some non stationary conditional centered random fields. A general central limit theorem adapted to this context has been obtained in \cite{A-ComJan98} for self-normalized sums (see also \cite{A-JenKun94} in the stationary case and without self-normalization). But, contrary to these papers where the observation domain is assumed to be of the form $[-n,n]^d$, we consider domains that may increase continuously up to $\RR[d]$. This particularity, which seems more relevant, requires an extension of the results in \cite{A-ComJan98} and \cite{A-JenKun94} to triangular arrays. This new central limit theorem is presented in Appendix~\ref{annexe-tcl}.

\section{Background on marked Gibbs point processes and definition of residuals}\label{sec-prel}

\subsection{General notation, configuration space}
We denote by  $\sB(\RR[d])$ the space of bounded Borel sets in $\RR[d]$. For any $\Lambda\in\sB(\RR[d])$, $\Lambda^c$ denotes the complementary set of $\Lambda$ inside $\RR[d]$.
The norm $|.|$ will be used without ambiguity for different kind of objects. For a vector $\Vect x$, $|\Vect x|$ represents the uniform norm of $\Vect x$; For a countable set $\mathcal J$, $|\mathcal J|$ represents the number of elements belonging to $\mathcal J$; For a set $\Delta\in\sB(\RR[d])$, $|\Delta|$ is  the volume of $\Delta$.

Let $\Mat M$ be a matrix, we denote by $\|\Mat M\|$ the Frobenius norm of $\Mat M$ defined by
$\|\Mat M\|^2=Tr(\tr{\Mat M}\Mat M)$, where $Tr$ is the trace operator. For a vector $\Vect x$, $\|\Vect x\|$ is simply its euclidean norm.

For all $\Vect x\in\RR[d]$ and $\rho>0$, $\mathcal{B}(\Vect x,\rho):=\{\Vect y,\ |\Vect y-\Vect x|<\rho\}$. Let us also consider the short notation, for $i\in\ZZ[d]$,  $\mathbbm{B}_i(\rho) = \mathcal{B}(i,\rho) \cap \ZZ[d]$.

The space $\RR[d]$ is endowed with the Borel $\sigma$-algebra and the Lebesgue measure $\lambda$. Let $\mSp$ be a measurable space, which aims at being the mark space, endowed with the $\sigma$-algebra $\mathcal M$ and the probability measure $\lambda^{\mE}$. The state space of the point processes will be $\sSp:=\RR[d]\times\mSp$ measured by $\mu:=\lambda\otimes\lambda^{\mE}$.
We shall denote for short $x^m=(x,m)$ an element of $\sSp$. 

The space of point configurations will be denoted by $\Omega=\Omega(\sSp)$. This is the set of simple integer-valued measures on $\sSp$.  It is endowed with the $\sigma$-algebra $\mathcal F$ generated by the sets $\{\varphi\in\Omega,\  \varphi(\Lambda\times A)=n\}$ for all $n\in\NN$, for all  $A\in\mathcal M$ and for all $\Lambda\in\sB(\RR[d])$. For any $x^m\in\sSp$ and $\varphi\in\Omega$, we denote $x^m\in\varphi$ if $\varphi(x^m)>0$. For any $\varphi\in\Omega$ and any $\Lambda\in\sB(\RR[d])$, we denote $\varphi_{\Lambda}:=\varphi_{\Lambda\times\mSp}$ the projection of $\varphi$ onto $\Lambda\times\mSp$, which is just the mesure $\sum_{x^m\in\varphi\cap(\Lambda\times\mSp)}\delta_{x^m}$, where $\delta_x$ is the Dirac measure at $x$. We will use without ambiguity some set notation for elements in $\Omega$, e.g. $\varphi\cup\{x^m\}=\varphi\cup x^m:=\varphi+\delta_{x^m}$ and 
for $x^m\in\varphi$, $\varphi\smallsetminus \{x^m\}=\varphi\smallsetminus x^m:=\varphi-\delta_{x^m}$. For any $\Lambda\in\sB(\RR[d])$, the number of elements of $\varphi_{\Lambda}$ is denoted by $|\varphi_{\Lambda}|:=\varphi(\Lambda\times\mSp)$.

\subsection{Marked Gibbs point processes}

The framework of this paper is restricted to stationary marked Gibbs point processes. Since we are interested in asymptotic properties, we consider these point processes on the infinite volume $\RR[d]$. Let us briefly recall their definition.

A marked point process $\Phi$ is a $\Omega$-valued random variable, with probability distribution $P$ on $(\Omega,\mathcal F)$. 
The most prominent marked point process is the marked Poisson process $\pi^{\nu}$ with intensity measure $\nu$ on $\RR[d]$ (and mark density $\lambda^{\mE}$). The homogeneous marked Poisson process arises when $\nu=z\lambda$, with $z>0$.

 Let $\ParV \in \SpPar$, where $\SpPar$ is some compact set of $\RR[p]$ (for some $p\geq 1$). For any $\Lambda\in\sB(\RR[d])$, let us consider the parametric function $V_{\Lambda}(.;\ParV)$ from $\Omega$ into $\RR[]\cup\{+\infty\}$. For fixed $\ParV$, ($V_{\Lambda}(.;\ParV))_{\Lambda\in\sB(\RR[d])}$ constitutes a compatible family of energies if, for every $\Lambda\subset\Lambda'$ in $\sB(\RR[d])$, there exists a measurable function $\psi_{\Lambda,\Lambda'}$ from $\Omega$ into $\RR[]\cup\{+\infty\}$ such that 
\begin{equation}\label{compatibility}
 \forall\varphi\in\Omega\quad V_{\Lambda'}(\varphi;\ParV)=V_{\Lambda}(\varphi;\ParV)+\psi_{\Lambda,\Lambda'}(\varphi_{\Lambda^c};\ParV).
 \end{equation}
  From a physical point of view, $V_{\Lambda}(\varphi_{\Lambda};\ParV)$ is the energy of $\varphi_{\Lambda}$ in $\Lambda$ given the outside configuration $\varphi_{\Lambda^c}$.
 The following definition is the classical way to define Gibbs measures through their conditional specifications (see \cite{B-Pre76}).
 
 \begin{definition}A probability measure $P_{\ParV}$ on $\Omega$  is a marked Gibbs measure for the compatible family of energies $(V_{\Lambda}(.;\ParV))_{\Lambda\in\sB(\RR[d])}$ and the intensity $\nu$ if for every $\Lambda\in\sB(\RR[d])$, for $P_{\ParV}$-almost every outside configuration $\varphi_{\Lambda^c}$, the law of $P_{\ParV}$ given $\varphi_{\Lambda^c}$ admits the following conditional density with respect to $\pi^{\nu}$:
 $$f_{\Lambda}(\varphi_{\Lambda}|\varphi_{\Lambda^c};\ParV)=\frac{1}{Z_{\Lambda}(\varphi_{\Lambda^c};\ParV)}e^{-V_{\Lambda}(\varphi;\ParV)},$$
where $Z_{\Lambda}(\varphi_{\Lambda^c};\ParV)$ is a normalization called the partition function.
\end{definition}

The existence of a Gibbs measure on $\Omega$ which satisfies these conditional specifications is a difficult issue. We do not want to open this discussion here and we will assume that the Gibbs measures we consider exist. We refer the interested reader to \cite{B-Rue69,B-Pre76,A-BerBilDro99,A-Der05,A-DerDroGeo09},  see also Section \ref{sec-disc} for several examples.

In this article, we focus on stationary marked point processes on $\sSp$, i.e. on point processes admitting a conditional density with respect to the homogeneous marked Poisson process $\pi$. Moreover, without loss of generality, the intensity of the Poisson process, $z$, is fixed to~1. We assume in a first step that the family of energies is hereditary, which means that for any $\Lambda\in\sB(\RR[d])$, for any $\varphi\in\Omega$, and for all $x^m\in\Lambda\times\mSp$, 
\begin{equation}\label{heredite}V_{\Lambda}(\varphi;\ParV))=+\infty \Rightarrow V_{\Lambda}(\varphi \cup\{x^m\};\ParV))=+\infty,\end{equation}
or equivalently, for all $x^m\in\varphi_{\Lambda}$,  $f_{\Lambda}(\varphi_{\Lambda}|\varphi_{\Lambda^c};\ParV)>0 \Rightarrow f_{\Lambda}(\varphi_{\Lambda}\smallsetminus\{x^m\}|\varphi_{\Lambda^c};\ParV)>0$. The non-hereditary case will be considered in Section \ref{non-hereditary}. The main assumption is then the following.

\begin{list}{}{}
\item \textbf{[Mod-E]}: For any $\ParV \in \SpPar$,  the  compatible family of energies $(V_{\Lambda}(.;\ParV))_{\Lambda\in\sB(\RR[d])}$ is hereditary, invariant by translation, and such that an associated Gibbs measure $P_{\ParV}$ exists and is stationary. Our data consist in the realization of a point process with Gibbs measure $P_{\ParVT}$. The vector $\ParVT$ is thus the true parameter to be estimated, assumed to be in $\mathring{\SpPar}$.
\end{list}

The local energy to insert a marked point $x^{m}$ into the configuration $\varphi$ is defined for any $\Lambda$ containing $x^m$ by 
$$V^{}\left( x^{m}|\varphi; \ParV \right):=V_{\Lambda}(\varphi \cup\{x^m\})-V_{\Lambda}(\varphi).$$
From the compatibility of the family of energies, i.e. (\ref{compatibility}), this definition does not depend on $\Lambda$. We restrict our study to finite-range interaction point processes, which is the main limitation of this paper.

\begin{list}{}{}
\item \textbf{[Mod-L]}: There exists $D\geq 0$ such that for all $(m,\varphi) \in \mSp\times \Omega$
$$\VIPar{0^{m}}{\varphi}{\ParV} = \VIPar{0^{m}}{\varphi_{\mathcal B(0,D)}}{\ParV}.$$
\end{list}

\subsection{Definitions of residuals for spatial point processes} \label{sec-desc}

The basic ingredient for the definition of residuals is the so-called GNZ formula stated below.
\begin{theorem}[Georgii-Nguyen-Zessin Formula]
Under \textbf{[Mod-E]}, for any function $h(\cdot,\cdot;\ParV): \sSp\times \Omega\to \RR[]$ (eventually depending on some parameter $\ParV$) such that the following quantities are finite, then
\begin{equation}\label{GNZnonstat}
\Esp\left( \ism[ \mathbb{R}^d \times \mSp]{h\left(x^m,\Phi;\ParV\right) e^{- \VIPar{x^m}{\Phi}{\ParVT}}} \right) = 
\Esp\left( \sum_{x^m \in \Phi } h\left(x^m,\Phi\setminus x^m;\ParV\right)  \right) ,
\end{equation}
where $\Esp$ denotes the expectation with respect to $P_{\ParVT}$.
\end{theorem}
For stationary marked Gibbs point processes, (\ref{GNZnonstat}) reduces to
\begin{equation}\label{GNZ}
\Esp\left( h\left(0^M,\Phi;\ParV\right) e^{- \VIPar{0^M}{\Phi}{\ParVT}} \right) = 
\Esp\left( h\left(0^M,\Phi\setminus 0^M;\ParV\right)  \right) 
\end{equation}
where $M$ denotes a random variable with probability distribution $\mm$. 
The following definition is based on empirical versions of both terms appearing in~(\ref{GNZ}).
\begin{definition} For any bounded domain $\Lambda$, let us define the $h-$innovations (denoted by $I_\Lambda$) and the $h-$residuals (denoted by $R_\Lambda$ and depending on  an estimate $\widehat{\ParV}$ of  $\ParVT$) by
\begin{eqnarray*}
I_{\Lambda}\left( \varphi;h,\ParVT \right) &:=&  \ism[\Lambda\times \mSp]{h \left( x^m , \varphi ; \ParVT \right) e^{-\VIPar{x^m}{\varphi}{\ParVT}}} -  \sum_{x^m\in \varphi_\Lambda} h \left( x^m , \varphi\setminus x^m ; \ParVT \right)\\
R_{\Lambda}\left( \varphi;h,\widehat{\ParV} \right) &:= & \ism[\Lambda\times \mSp]{h \left( x^m , \varphi ; \widehat{\ParV} \right) e^{-\VIPar{x^m}{\varphi}{\widehat{\ParV}}}} - \sum_{x^m\in \varphi_\Lambda} h \left( x^m , \varphi\setminus x^m ; \widehat{\ParV} \right).
\end{eqnarray*}
\end{definition}
From a practical point of view, the last notion is the most interesting since it provides a computable measure. The main examples considered by \citeauthor{A-BadTurMolHaz05} in \cite{A-BadTurMolHaz05} (in the context of stationary point processes) are obtained by setting $h \left( x^m , \varphi ; \ParV \right)=1$ for the raw residuals, $h \left( x^m , \varphi ; \ParV \right)=e^{\VIPar{x^m}{\varphi}{\ParV}}$ for the inverse residuals and $h \left( x^m , \varphi ; \ParV \right)=e^{\VIPar{x^m}{\varphi}{\ParV}/2}$ for the Pearson residuals. In particular, one may note that the raw residuals constitutes a difference of two estimates of the intensity of the point process (up to a normalisation by $|\Lambda|$): the first one is a parametric one and depends on the model while the second one is a nonparametric one (since it is equal to $|\varphi_\Lambda|$). Another more evolved example is to consider the function defined for $r>0$ by
$$
h_r(x^m, \varphi;\ParV): = \mathbf{1}_{[0,r]} (d(x^m,\varphi)) \; e^{\VIPar{x^m}{\varphi}{\ParV}} 
$$
where $d(x^m,\varphi)=\inf_{y^m\in \varphi} \|y-x\|$. Considering this function leads to 
\begin{equation} \label{eq-Resempty}
R_\Lambda(\varphi;h_r,\widehat\ParV) = \ism{\mathbf{1}_{[0,r]} (d(x^m,\varphi))} - \sum_{x^m\in\varphi_\Lambda} h_r(x^m,\varphi\setminus x^m;\widehat{\ParV}).
\end{equation}
Then for a large window $R(\varphi;h_h,\widehat{\ParV})/|\Lambda|$ leads to a difference of two estimates of the well-known empty space function $F$ at distance $r$. Recall that  for a marked stationary point process (see \cite{B-MolWaa03} for instance) 
$$
F(r) := P\left( d(0^M ,\Phi) \leq r \right).
$$
The first term in the right hand side of \eqref{eq-Resempty} corresponds to the natural nonparametric estimator of $F(r)$ while the second one  is a parametric  estimator of $F(r)$.




\section{Asymptotic properties} \label{sec-asymp}
From now on, we assume that the point process satisfies \textbf{[Mod-E]} and \textbf{[Mod-L]}, that is \textbf{[Mod]}. The realization of $\Phi\sim P_{\ParVT}$ is assumed to be observed in a domain $\Lambda_n\oplus D^{+}$, with $D^{+}\geq D$, aimed at growing up to $\RR[d]$ as $n\to +\infty$. According to the locality assumption \textbf{[Mod-L]},  we are thus ensured that the $h-$innovations and $h-$residuals can be computed. 

The aim of this section is to present several asymptotic properties for  $I_{\Lambda_n}$ and $R_{\Lambda_n}$. We prove their consistency and we propose two asymptotic normality results  within different frameworks:
\begin{itemize}
\item Framework~1: for a fixed test function $h$, $\Lambda_n$ is a cube, divided into a fixed finite number of sub-cubes (which will increase with $\Lambda_n$). The purpose is then to obtain the asymptotic normality for the vector composed of the $h-$residuals computed in each sub-cube.
\item Framework~2: we consider $h_1,\ldots,h_s$ $s$ different test functions and the aim is to obtain the asymptotic normality of the vector composed of the $h_j-$residuals computed on~$\Lambda_n$. 
\end{itemize}
In both frameworks, an estimate of $\ParVT$ is involved. We assume that it is computed from the full domain $\Lambda_n$ with the same data used to evaluate the $h-$residuals, which is a natural setting in practice.
Moreover, contrary to the previous works dealing with asymptotic properties on Gibbs point processes ({\it e.g.} \cite{A-JenKun94}, \cite{A-ComJan98} or \cite{A-BilCoeDro08}), where $\Lambda_n$  is assumed to be of the discrete form $[-n,n]^d$, we consider general domains  that may grow continuously up to $\RR[d]$. 

The asymptotic results obtained in this section are the basis to derive  goodness-of-fit tests, as presented in Section~\ref{sec-gof}.





\subsection{Consistency of the residuals process}

We obtain consistency results for $I_{\tilde\Lambda_n}\left( \Phi;h,\ParVT \right)$ and $R_{\tilde\Lambda_n}\left( \Phi;h,\widehat{\ParV}_n(\Phi) \right)$, where for all $n\geq 1$, $\tilde\Lambda_n \subset \Lambda_n$, $(\tilde\Lambda_n)_{n\geq 1}$ and $(\Lambda_n)_{n\geq 1}$ are regular sequences whose size increases to $\infty$. 

The assumption \textbf{[C]} gathers the two following assumptions:
\begin{itemize}
\item[\textbf{[C1]}] 
$$
\Esp \left( \left| h \left( 0^M , \Phi ; \ParVT \right) \right| e^{- V^{}\left( 0^M|\Phi; \ParVT \right)}
\right) <+\infty.
$$
\item[\textbf{[C2]}] For all $(m,\varphi) \in \mSp \times \Omega$, the functions $h \left( 0^m , \varphi ; \ParV \right)$ and  $f \left( 0^m , \varphi ; \ParV \right):= h \left( 0^m , \varphi ; \ParV \right) e^{- V^{}\left( 0^m|\varphi; \ParV \right)}$ are continuously differentiable with respect to $\ParV$ in a neighborhood $\mathcal{V}(\ParVT)$ of $\ParVT$ and
$$
\Esp \left( \left\|  \Vect{f}^{(1)}_{} \left( 0^M , \Phi ;   \ParVT \right)\right\| \right) <+\infty
\quad \mbox{ and }\quad
\Esp \left( \left\|  \Vect{h}^{(1)}_{} \left( 0^M , \Phi ;   \ParVT \right)\right\| e^{- V^{}\left( 0^M|\Phi; \ParVT \right)}\right) <+\infty,
$$
where $\Vect f^{(1)}$ denotes the gradient vector of $f$ with respect to $\ParV$.
\end{itemize}
Concerning the residuals process, we also need to assume
\begin{itemize}
\item[\textbf{[E1]}] The estimator $\widehat{\ParV}_n(\varphi)$ of $\ParVT$, computed from the full observation domain $\Lambda_n$,  converges for $P_{\ParVT}-$a.e. $\varphi$ towards $\ParVT$, as $n\to +\infty$.
\end{itemize}

\begin{proposition} \label{prop-cons}
Assuming~\textbf{[Mod]}, we have as $n \to +\infty$
\begin{itemize}
\item[(a)] Under \textbf{[C1]}: for $P_{\ParVT}-$a.e. $\varphi$, $|\tilde\Lambda_n|^{-1} I_{\tilde\Lambda_n}\left( \varphi;h,\ParVT \right)$ converges towards 0.
\item[(b)] Under \textbf{[C]} and \textbf{[E1]}: for $P_{\ParVT}-$a.e. $\varphi$, $|\tilde\Lambda_n|^{-1} R_{\tilde\Lambda_n}\left( \varphi;h,\widehat{\ParV}_n(\varphi) \right)$ converges towards 0.
\end{itemize}
\end{proposition}

\begin{remark}
Assumption  \textbf{[Mod-L]}, while useful to allow the computation of the residuals in practice,  is actually useless to prove their consistency.
\end{remark}


\subsection{Asymptotic control in probability of the residuals process}

We provide in this section a control for the departure of the residuals from the innovations and $(\widehat{\ParV}_n-\ParVT)$. This is a crucial result to investigate the asymptotic normality of the residuals.  We need the folllowing assumptions. 
\begin{itemize}
\item[\textbf{[N1]}] For all $(m,\varphi) \in \mSp \times \Omega$, the functions $h \left( 0^m , \varphi ; \ParV \right)$ and  $f \left( 0^m , \varphi ; \ParV \right)$ (defined in \textbf{[C1]}) are twice continuously differentiable with respect to $\ParV$ in a neighborhood $\mathcal{V}(\ParVT)$ of $\ParVT$ and 
$$
\Esp \left( \left\|  \Vect{\underline{f}}^{(2)}_{} \left( 0^M , \Phi ;  \ParVT \right)\right\| \right) <+\infty
\quad \mbox{ and }\quad
\Esp \left( \left\|  \Vect{\underline{h}}^{(2)}_{} \left( 0^M , \Phi ;  \ParVT \right)\right\|e^{- V^{}\left( 0^M|\Phi; \ParVT \right)} \right) <+\infty,
$$
where $  \Vect{\underline{g}}^{(2)}_{} \left( 0^m , \varphi ;  \ParVT \right) = \left(\frac{\partial^2}{\partial\theta_j \partial \theta_k} g \left( 0^m , \varphi ; \ParVT \right)\right)_{1\leq j,k\leq p}$ for $g=f,h$.
\item[\textbf{[E2]}] There exists a random vector $\Vect{T}$ such that the following convergence holds as $n\to+\infty$
$$|\Lambda_n|^{1/2} \left( \widehat{\ParV}_n(\Phi)-\ParVT\right) \stackrel{d}{\longrightarrow} \Vect{T}.
$$
\end{itemize}

\begin{proposition} \label{prop-equiv}
Under assumptions \textbf{[C]}, \textbf{[N1]} and \textbf{[E1-2]}, assuming that $|\tilde\Lambda_n|=\mathcal{O}(|\Lambda_n|)$, then as $n\to +\infty$, 
\begin{equation}\label{eq-proba}
R_{\tilde\Lambda_n}\left( \Phi;h,\widehat{\ParV}_n(\Phi) \right) 
= I_{\tilde\Lambda_n}\left( \Phi;h,\ParVT \right) 
- |\tilde\Lambda_n| \tr{\left( \widehat{\ParV}_n(\Phi)-\ParVT \right)} \Vect{\mathcal{E}}\left(h;\ParVT\right) + o_{P}(|\tilde\Lambda_n|^{1/2}) ,
\end{equation}
where $\Vect{\mathcal{E}}\left(h;\ParVT\right)$ is the vector defined by
\begin{equation} \label{eq-defE}
\Vect{\mathcal{E}}\left(h;\ParVT\right) := \Esp\left( 
h \left( 0^M , \Phi ; \ParVT \right)  \Vect{V}^{(1)}\left( 0^M|\Phi; \ParVT \right)
e^{- V^{}\left( 0^M|\Phi; \ParVT \right)}
\right).
\end{equation}
\end{proposition}

The notation $X_n(\Phi)=o_P(w_n)$ means that $w_n^{-1}X_n(\Phi)$ converges in probability towards 0 as $n$ tends to infinity.

\begin{remark}
Note that for exponential family models, $\Vect{V}^{(1)}(x^m|\varphi;\ParVT)$ corresponds to the vector of sufficient statistics (see Section \ref{sec-disc} for more details).
\end{remark}

\subsection{Assumptions required for the asymptotic normality results}
Apart from the assumptions \textbf{[Mod]}, \textbf{[C]} and \textbf{[N1]} on the model, we will need to assume \textbf{[N2-4]} below. All these assumptions are fulfilled by many models as proved in Section \ref{sec-disc}.

\begin{itemize}
\item[\textbf{[N2]}] For any bounded domain $\Lambda$, for any $\ParV \in \mathcal{V}(\ParVT)$,
$$\Esp\left( \left|{I}_{\Lambda}\left( \Phi;h,\ParVT \right) \right|^3\right)<+\infty.$$

\item[\textbf{[N3]}] For any sequence of bounded domains  $\Gamma_n$ such that $\Gamma_n\to 0$ when $n\to\infty$, for any $\ParV \in \mathcal{V}(\ParVT)$, $$\Esp\left( {I}_{\Gamma_n}\left( \Phi;h,\ParV \right)^2\right)\longrightarrow 0.$$
 
\item[\textbf{[N4]}] For any $\varphi \in \Omega$ and any bounded domain $\Lambda$, $I_{\Lambda}\left( \varphi ; \ParV \right)$ depends only on $\varphi_{\Lambda \oplus D}$.
\end{itemize}

Concerning the properties required for the estimator $\widehat{\ParV}_n$, we need its consistency through \textbf{[E1]} and to refine \textbf{[E2]} into   \textbf{[E2(bis)]} below. Note that the maximum pseudolikelihood estimator satisfies these assumptions for many models (see section \ref{sec-hypEst}). 

\begin{itemize}
\item[\textbf{[E2(bis)]}] The estimate admits the following expansion
$$
\widehat{\ParV}_n(\Phi) - \ParVT = \frac{1}{|\Lambda_n|} \Vect{U}_{\Lambda_n}\left( \Phi ; \ParVT \right) + o_{P}( |\Lambda_n|^{-1/2}),
$$
where, for any $\ParV \in \mathcal{V}(\ParVT)$,
\begin{itemize}
\item[(i)] for any $\varphi\in \Omega$ and for two disjoint bounded domains $\Lambda_1,\Lambda_2$, 
$$\Vect{U}_{\Lambda_1\cup \Lambda_2}\left( \varphi ; \ParV \right)= \Vect{U}_{\Lambda_1}\left( \varphi ; \ParV \right)+ \Vect{U}_{\Lambda_2}\left( \varphi ; \ParV \right),$$
\item[(ii)] for all $j=1,\ldots,p$ and any bounded domain $\Lambda$
$$\Esp\left(\left|\left( \Vect{U}_{\Lambda}\left( \Phi ; \ParV \right)
\right)_j\right|^3 \right) <+\infty,$$
\item[(iii)] for all $j=1,\ldots,p$ and for any bounded domain $\Lambda$
$$
\Esp \left( \left. \left( \Vect{U}_{\Lambda}\left( \Phi ; \ParV \right)
\right)_j \right| \Phi_{\Lambda^c}\right) =0,
$$

\item[(iv)] for all $j=1,\ldots,p$ and for any sequence of bounded domains  $\Gamma_n$, 
$$\Esp\left(\left( \Vect{U}_{\Gamma_n}\left( \Phi ; \ParV \right)
\right)_j^2 \right)\longrightarrow 0\quad\text{as}\quad\Gamma_n\to 0,$$

\item[(v)] for any $\varphi \in \Omega$ and any bounded domain $\Lambda$, $\Vect{U}_{\Lambda}\left( \varphi ; \ParV \right)$ depends only on $\varphi_{\Lambda \oplus D}$.
\end{itemize}
\end{itemize}

\begin{remark} \label{rem-E2bis}
Assumption \textbf{[E2(bis)]} implies \textbf{[E2]}. Indeed, under this assumption one may apply Theorem~2.1 of \cite{A-JenKun94} and assert: there exists a matrix $\Mat{\Sigma}$ such that $|\Lambda_n|^{-1/2} \Vect{U}_{\Lambda_n}\left( \Phi ; \ParVT \right) \stackrel{d}{\to} \mathcal{N}(0,\Mat{\Sigma})$, as $n\to +\infty$. 
\end{remark}

\subsection{Asymptotic normality of the $h-$residuals computed on subdomains of $\Lambda_n$}\label{sec-fwk1}

In this framework, we give ourself a test function $h$ and we compute the $h-$residuals on disjoint subdomains of $\Lambda_n$. In this context, we assume that the domain $\Lambda_n$ is a cube and is divided into a fixed number of subdomains as follows
$$\Lambda_n := \bigcup_{j \in \mathcal J} \Lambda_{j,n}$$ where $\mathcal J$ is a finite set and all the $\Lambda_{j,n}$ are disjoint cubes with the same volume $ \left| \Lambda_{0,n}\right|$ increasing up to $+\infty$. Let us denote by ${\Vect{R}}_{\mathcal{J},n} \left(\varphi;h,\widehat{\ParV}_n\right)$ the vector of the residuals computed on each subdomain, i.e. ${\Vect{R}}_{\mathcal{J},n} \left(\varphi;h,\widehat{\ParV}_n\right)= \left( R_{\Lambda_{j,n}}\left( \varphi;h,\widehat{\ParV}_n\right) \right)_{j\in \mathcal J}$.

According to Proposition~\ref{prop-equiv} and in view of \textbf{[E2(bis)]}, we introduce the following notation 
\begin{equation} \label{def-Rinfty}
R_{\infty,\Lambda}(\varphi;h,\ParV):=I_{\Lambda}(\varphi;h,\ParV)- \tr{{\Vect{U}}_{\Lambda}\left( \varphi ; \ParV \right)} \Vect{\mathcal{E}}\left( h;\ParV\right)
\end{equation}
for any $\varphi\in \Omega$, for any bounded domain $\Lambda$ and for any $\ParV \in \SpPar$. 

\begin{proposition} \label{prop-fwk1}
Assume that
\begin{itemize}
\item The parametric model satisfies \textbf{[Mod]}.
\item The energy function and the test function $h$ satisfy \textbf{[C]} and \textbf{[N1-4]}.
\item The energy function and the estimate $\widehat{\ParV}_n$ satisfy \textbf{[E1]} and \textbf{[E2(bis)]}.
\end{itemize}
Then, the following convergence in distribution holds, as $n\to +\infty$
\begin{equation} \label{eq-fwk1}
{} \left| \Lambda_{0,n}\right|^{-1/2}{\Vect{R}}_{\mathcal{J},n} \left(\Phi;h,\widehat{\ParV}_n\right) \stackrel{d}{\longrightarrow} \mathcal{N}\left( 0,\Mat{\Sigma}_1(\ParVT) \right),
\end{equation}
where $\Mat{\Sigma}_1(\ParVT)=\lambda_{Inn}\;  \Mat{I}_{|\mathcal J|} +|\mathcal{J}|^{-1}(\lambda_{Res} -\lambda_{Inn} )\;\Mat{J}$ with  $\Mat{J}=\Vect{e}\tr{\Vect{e}}$ and $\Vect{e}=\tr{(1,\ldots,1)}$. The constants $\lambda_{Inn}$ and $\lambda_{Res}$ are respectively defined by
\begin{eqnarray}
\lambda_{Inn}&=& D^{-d} \sum_{|k|\leq 1} \Esp\left(
{I}_{\Delta_0(D)}\left( \Phi;h,\ParVT \right)
{I}_{\Delta_k(D)}\left( \Phi;h,\ParVT \right)
 \right) ,\label{eq-defLambdaInn}\\
\lambda_{Res} &=&  D^{-d} \sum_{|k|\leq 1} 
\Esp\left( 
R_{\infty,\Delta_0(D)}( \Phi;h,\ParVT) R_{\infty,\Delta_k(D)}( \Phi;h,\ParVT)
\right) 
 \label{eq-defLambdaRes},
\end{eqnarray}
 where, for all $k\in \ZZ[d]$, $\Delta_k(D)$ is the cube centered at $kD$ with side-length $D$. 
\end{proposition}

From this asymptotic normality result, we can deduce the convergence for the norm of the centered residuals.  This is the basis for a generalization of the quadrat counting test discussed in Section \ref{sec-gof}. We denote by $\overline{\Vect{R}}_{\mathcal{J},n}(\varphi;h)$ the mean residuals over all subdomains, that is
$\overline{\Vect{R}}_{\mathcal{J},n}(\varphi;h)=|\mathcal J|^{-1}\sum_{j\in\mathcal J} R_{\Lambda_{j,n}}\left( \varphi;h,\widehat{\ParV}_n\right)$.

\begin{corollary}\label{cor-quadrat}
Under the assumptions of Proposition \ref{prop-fwk1},
\begin{equation} \label{eq-cor}
|\Lambda_{0,n}|^{-1} \| \Vect{R}_{\mathcal{J},n}(\Phi;h)- \overline{\Vect{R}}_{\mathcal{J},n}(\Phi;h) \|^2 \stackrel{d}{\longrightarrow} \lambda_{Inn}\ \chi^2_{|\mathcal{J}|-1}.
\end{equation}

\end{corollary}

\begin{proof}
An easy computation shows that $\lambda_{Inn}$ and $\lambda_{Res}$ are the two eigenvalues of $\Mat{\Sigma}_1(\ParVT)$ with respective order $|\mathcal J|-1$ and $1$. Let $\Mat{P}_{Inn}$ be the matrix of orthonormalized eigenvectors associated to $\lambda_{Inn}$. This matrix of size $(|\mathcal{J}|,|\mathcal{J}|-1)$ satisfies by definition $\tr{\Mat{P}}_{Inn}\Mat{P}_{Inn} = \Mat{I}_{|\mathcal{J}|-1}$ and, from (\ref{eq-fwk1}), $|\Lambda_{0,n}|^{-1}  \; \| \tr{\Mat{P}}_{Inn} \Vect{R}_{\mathcal{J},n}(\varphi;h)\|^2 \stackrel{d}{\to} \lambda_{Inn}\ \chi^2_{|\mathcal{J}|-1}$. Moreover, it is easy to check that ${\Mat{P}_{Inn}}\tr{\Mat{P}}_{Inn} = \Mat{I}_{|\mathcal{J}|}-|\mathcal{J}|^{-1} \Mat{J}_{|\mathcal{J}|}$ which leads to $\| \tr{\Mat{P}}_{Inn} \Vect{R}_{\mathcal{J},n}(\varphi;h)\|^2= \| \Vect{R}_{\mathcal{J},n}(\varphi;h)- \overline{\Vect{R}}_{\mathcal{J},n}(\varphi;h) \|^2$.
\end{proof}

\begin{remark}
The asymptotic covariance matrix $\Mat{\Sigma}_1(\ParVT)$ and  $\lambda_{Inn}$ involve  only the covariance structure of the innovations (or the residuals) in a finite box around $0$.  This comes from the locality assumption \textbf{[Mod-L]}, also involved in \textbf{[N4]} and \textbf{[E2(bis)]}. A challenging task in practice is to estimate  $\lambda_{Inn}$ and $\lambda_{Res}$ (and so $\Mat{\Sigma}_1(\ParVT)$), this issue is investigated in Section \ref{sec-estMat}. 
\end{remark}

\subsection{Asymptotic normality of the $(h_j)_{j=1,\ldots,s}-$residuals computed on $\Lambda_n$}\label{sec-fwk2}
In this framework, we consider $s$ different test functions and we compute all $h_j-$residuals on the same domain $\Lambda_n$, which is assumed to be a cube growing up to $\RR[d]$ when $n\to+\infty$. 

We present an asymptotic normality result for the random vector $\left( R_{\Lambda_{n}}\left( \Phi;h_j,\widehat{\ParV}_n\right) \right)_{j=1,\ldots,s}$.

\begin{proposition} \label{prop-fwk2}
Assume that
\begin{itemize}
\item The parametric model satisfies \textbf{[Mod]}.
\item The energy function and the test functions $h_j$ (for $j=1,\ldots,s$) satisfy \textbf{[C]} and \textbf{[N1-4]}.
\item The energy function and the estimate $\widehat{\ParV}_n$ satisfy \textbf{[E1]} and \textbf{[E2(bis)]}.
\end{itemize}
Then, the following convergence in distribution holds, as $n\to +\infty$
\begin{equation} \label{eq-fwk2}
{} \left| \Lambda_{n}\right|^{-1/2} \left( R_{\Lambda_{n}}\left( \Phi;h_j,\widehat{\ParV}_n\right) \right)_{j=1,\ldots,s} \stackrel{d}{\longrightarrow} \mathcal{N}\left( 0,\Mat{\Sigma}_2(\ParVT) \right),
\end{equation}
where $\Mat{\Sigma}_2(\ParVT)$ is the $(s, s)$ matrix given by
\begin{equation} \label{eq-defSigma2}
\Mat{\Sigma}_2(\ParVT) = 
D^{-d} \sum_{|k| \leq 1} 
\Esp\left( 
\Vect{R}_{\infty,\Delta_0(D)}(\Phi;\Vect{h},\ParVT) \; \tr{\Vect{R}_{\infty,\Delta_k(D)}(\Phi;\Vect{h},\ParVT)}
\right),
\end{equation}
where $\Vect{R}_{\infty,\Lambda}(\varphi,\Vect{h},\ParVT):=\left( R_{\infty,\Lambda} (\varphi;h_j,\ParVT) \right)_{j=1,\ldots,s}$, see (\ref{def-Rinfty}), and where, for all $k\in \ZZ[d]$, $\Delta_k(D)$ is the cube centered at $kD$ with side-length $D$.

\end{proposition}


\section{Estimation and positivity of the asymptotic covariance matrices} \label{sec-estMat}

\subsection{Statement of the problem}


The aim of this section is to provide a condition under which, on the one hand the matrices $\Mat{\Sigma}_1(\ParVT)$ and $\Mat{\Sigma}_2(\ParVT)$, defined in Propositions~\ref{prop-fwk1} and~\ref{prop-fwk2}, are positive-definite, and  on the other hand  $\lambda_{Inn}$, involved  in Corollary \ref{cor-quadrat}, is positive. Then we define estimators of  $\Mat{\Sigma}_1^{-1/2}(\ParVT)$, $\lambda_{Inn}^{-1}$ and $\Mat{\Sigma}_2^{-1/2}(\ParVT)$. As a consequence, we will be in position to normalize and estimate the quantities arising in (\ref{eq-fwk1}), (\ref{eq-cor}) and (\ref{eq-fwk2}) so that they converge to a free law.

Before this, let us focus on the particular form of the matrix $\Mat{\Sigma}_1(\ParVT)$. This $(|\mathcal J|,|\mathcal J|)$ matrix  has two eigenvalues $\lambda_{Inn}$ and $\lambda_{Res}$ (respectively defined by~\eqref{eq-defLambdaInn} and~\eqref{eq-defLambdaRes}), whose multiplicity is $|\mathcal J|-1$ for $\lambda_{Inn}$ and 1 for $\lambda_{Res}$.  By using the Gram-Schmidt process for orthonormalizing the eigenvectors of $\Mat{\Sigma}_1(\ParVT)$, one obtains the explicit form for the squared inverse of this matrix, provided $\lambda_{Inn}$ and $\lambda_{Res}$  do not vanish:
$$
\Mat{\Sigma}_1^{-1/2}(\ParVT) = \frac{1}{\sqrt{\lambda_{Inn}}} \; \Mat{I}_{|\mathcal{J}|} \; + \; \frac1{|\mathcal{J}|} \left( \frac1{\sqrt{\lambda_{Res}}}-\frac1{\sqrt{\lambda_{Inn}}} \right) \Mat{J},
$$
where $\Mat{J}=\Vect{e}\tr{\Vect{e}}$ and $\Vect{e}=\tr{(1,\ldots,1)}$. Therefore, estimating $\Mat{\Sigma}_1^{-1/2}(\ParVT)$ can be reduced to the estimation of these two eigenvalues $\lambda_{Inn}$ and $\lambda_{Res}$. 

Consequently, the estimation of  $\lambda_{Inn}$ and the covariance matrices $\Mat{\Sigma}_1(\ParVT)$ and $\Mat{\Sigma}_2(\ParVT)$  is achieved by estimating (\ref{eq-defLambdaInn}), (\ref{eq-defLambdaRes}) and~(\ref{eq-defSigma2}), which can be viewed as a particular case of estimating the matrix (actually a constant for the two first expressions) 
$$\Mat{M}({\ParVT})=D^{-d} \sum_{|k|\leq 1}\Esp \left( \Vect{Y}_{\Delta_0(D)}\left(\Phi;\ParVT\right) \tr{\Vect{Y}_{\Delta_k(D)}\left(\Phi;\ParVT\right) }
\right),$$
where, according to the assumptions involved in Propositions \ref{prop-fwk1} and \ref{prop-fwk2}, for any bounded domain $\Lambda$, $\Vect{Y}_{\Lambda}(\Phi;\ParV)$ is a random vector of dimension $q$ ($q=1$ or $s$) depending on $\ParV$, such that for any bounded domains $\Lambda,\Lambda_1,\Lambda_2$  ($\Lambda_1,\Lambda_2$ disjoint), for any $\ParV \in \mathcal{V}(\ParVT)$, for any $j=1,\ldots,q$ and any $\varphi\in \Omega$
\begin{itemize}
\item[(i)] $\Vect{Y}_{\Lambda_1\cup \Lambda_2}(\varphi;\ParV) = \Vect{Y}_{\Lambda_1}(\varphi;\ParV)+\Vect{Y}_{\Lambda_2}(\varphi;\ParV),$
\item[(ii)] $\Esp\left( \left(\Vect{Y}_{\Lambda}(\Phi;\ParV)\right)_j^2\right)<+\infty$,
\item[(iii)] $\Esp \left( \left. \left( \Vect{Y}_{\Lambda}\left( \Phi ; \ParV \right)
\right)_j \right| \Phi_{\Lambda^c}\right) =0,$
\item[(iv)] for any sequence of bounded domains  $\Gamma_n$, $\Esp\left(\left( \Vect{Y}_{\Gamma_n}\left( \Phi ; \ParV \right)
\right)_j^2 \right)\longrightarrow 0\quad\text{as}\quad\Gamma_n\to 0,$
\item[(v)] $\Vect{Y}_{\Lambda}(\varphi;\ParV)$ depends only on $\varphi_{\Lambda\oplus D}$.
\end{itemize}

\subsection{Positive definiteness of $\Mat{M}(\ParVT)$}

Let us consider the following assumption.

\begin{itemize}
\item[\textbf{[PD]}] For some $\Ltt:=\cup_{|i|\leq \left\lceil \frac{D}{\overline{\delta}} \right\rceil } \Delta_i(\overline{\delta})$ with $\overline{\delta}>0$, there exists $B\in\mathcal F$ and $A_0,\ldots,A_{\ell}$, ($\ell\geq 1$) disjoint events of $\cSptt_B:=\left\{\varphi \in \cSp: \varphi_{\Delta_i(\overline{\delta})}\in B, 1\leq |i| \leq 2\left\lceil \frac{D}{\overline{\delta}}\right\rceil \right\}$ such that 
\begin{itemize}
\item for $j=0,\ldots,\ell$, $P_{\ParVT}(A_j)>0$.
\item for all $\left(\varphi_0,\ldots,\varphi_{\ell} \right)\in A_0 \times \cdots \times A_{\ell}$ the $(\ell, q)$ matrix with entries $\left(\ensuremath{\Vect{Y}_{\Ltt} \left( \varphi_i ; \ParVT  \right)} \right)_j - \left(\ensuremath{\Vect{Y}_{\Ltt} \left( \varphi_0 ; \ParVT  \right)} \right)_j$ is injective, which means:
$$
\left(\forall \Vect{y}\in \RR[q], \tr{\Vect{y}} \left( \ensuremath{\Vect{Y}_{\Ltt} \left( \varphi_i ; \ParVT  \right)}-\ensuremath{\Vect{Y}_{\Ltt} \left( \varphi_0 ; \ParVT  \right)}\right)
  =0 \right)\Longrightarrow \Vect{y}=0.
$$

\end{itemize}
\end{itemize}

\begin{proposition} \label{prop-sdp} 
From the  definition of $\Vect{Y}_{\Lambda}(\Phi;\ParV)$  and under \textbf{[PD]}, the matrix $\Mat{M}(\ParVT)$ is positive-definite.
\end{proposition}

\begin{remark}
The assumption \textbf{[PD]} is associated to some characteristics of the point process $\Phi$. The parameter $\overline{\delta}$ is independent of the parameters involved in the different estimators (e.g. $D^\vee$ or $\delta$ arising in the next section). Given a model, the event $B$ and  $\overline{\delta}$ are chosen in order to let the different configurations sets $A_0,A_1,\ldots,A_\ell$ as simple as possible. For most models, a convenient choice is $B=\emptyset$ and  $\overline{\delta}\geq D$ (see  the examples treated in Appendix \ref{exemplesPD} for instance). 
\end{remark}

\subsection{Estimation of $\Mat{M}(\ParVT)$} \label{sec-defEstMat}

The dependence of  $\Mat{M}(\ParVT)$ on $D$ may be lightened thanks to the following lemma, whose proof is relegated to section \ref{proofdelta}.

\begin{lemma}\label{delta}
The matrix $\Mat{M}(\ParVT)$ can be rewritten for any $\delta>0$ and any $D^{\vee}\geq D$ as
$$
\Mat{M}(\ParVT) = \delta^{-d}  \sum_{|k|\leq \left\lceil \frac{D^{\vee}}{\delta}\right\rceil}\Esp \left( \Vect{Y}_{\Delta_0(\delta)}\left(\Phi;\ParVT\right) \tr{\Vect{Y}_{\Delta_k(\delta)}\left(\Phi;\ParVT\right) }\right),
$$
where $\Delta_k(\delta)$ is the cube with side-length $\delta$ centered at $k\delta$.
\end{lemma}

From this result, to achieve an estimation of $\Mat{M}(\ParVT)$, it is required to estimate the involved expectation  and $\ParVT$ (by $ \widehat{\ParV}_n$).  This is enough for the estimation of $\lambda_{Inn}$ for which $\Vect{Y}_\Lambda(\varphi;\ParV)=I_{\Lambda}(\varphi;\ParV)$. But when $\Vect{Y}_\Lambda(\varphi;\ParV)=R_{\infty,\Lambda}(\varphi;h,\ParV)$ or $\Vect{Y}_\Lambda(\varphi;\ParV)=\Vect{R}_{\infty,\Lambda}(\varphi;\Vect{h},\ParV)$, which appears in $\Mat{\Sigma}_1(\ParVT)$ and $\Mat{\Sigma}_2(\ParVT)$,  it can be noticed that $\Vect{Y}_\Lambda$ still depends on an expectation with respect to $P_{\ParVT}$, through the vector $\Vect{\mathcal{E}}(h,\ParVT)$ defined by \eqref{eq-defE}. Moreover, the vector $\Vect{U}_\Lambda$ in \textbf{[E2(bis)]} may also depend on such a term (this is the case for example when considering the maximum pseudolikelihood estimate as shown in Section~\ref{sec-hypEst}). This means that $\Vect{Y}_\Lambda(\varphi;\ParV)$ cannot be estimated only by $\Vect{Y}_\Lambda(\varphi; \widehat{\ParV}_n)$, but by $\widehat{\Vect{Y}}_{n,\Lambda}(\varphi; \widehat{\ParV}_n)$, where $\widehat{\Vect{Y}}_n$ is an estimator of $\Vect Y$. We assume in the sequel that  $\widehat{\Vect{Y}}_n$ satisfies the same properties $(i)-(v)$ as $\Vect{Y}$ and is a good estimator of $\Vect{Y}$ (see Proposition \ref{consistance}).  The explicit form of  $\widehat{\Vect{Y}}_n$ depends strongly on the estimate $\widehat{\ParV}_n$ (e.g. through $\Vect{U}_\Lambda$ in \textbf{[E2(bis)]}). When $\widehat{\ParV}_n$ is the maximum pseudolikelihood estimator, we provide explicit formulas for $\widehat{\Vect{Y}}_n$ in Section \ref{sec-estMatMPLE}.

Let us now specify an estimator of $\Mat{M}(\ParVT)$. Assume that the point process is observed in the domain $\Lambda_{n_0}\oplus D^+$ where $D^+\geq D$ and $\Lambda_{n_0}$ is a cube. For any $\delta$ such that $|\Lambda_{n_0}|\delta^{-d}\in\NN$, we may consider the decomposition $\Lambda_{n_0}=\cup_{k\in \mathcal K_{n_0}} \Delta_k(\delta)$, where  the $\Delta_k(\delta)$'s are disjoint cubes with side-length $\delta$ and centered, without loss of generality, at $k\delta$. For any such $\delta$, according to Lemma~\ref{delta}, a natural estimator of  $\Mat{M}(\ParVT)$ is, for any  $D^{\vee}\geq D$,
\begin{equation} \label{eq-estMn}
\widehat{\Mat{M}}_{n_0}(\varphi; \widehat{\ParV}_{n_0}(\varphi),\delta,D^\vee) = 
|\Lambda_{n_0}|^{-1}\!\!\!
\sum_{k\in \mathcal K_{n_0}} \sum_{j\in \mathbbm{B}_k\left(\left\lceil \frac{D^\vee}{\delta} \right\rceil \right) \cap \mathcal K_{n_0}} \!\! \!\!\!\widehat{\Vect{Y}}_{n_0,\Delta_j(\delta)}\left(\varphi;\widehat{\ParV}_{n_0}(\varphi)\right) \tr{\widehat{\Vect{Y}}_{n_0,\Delta_k(\delta)}\left(\varphi;\widehat{\ParV}_{n_0}(\varphi)\right) }.
\end{equation}
\begin{remark}
As suggested by Lemma \ref{delta}, the parameter $\delta$ in (\ref{eq-estMn}) may be chosen arbitrarily. Yet, while $\Mat{M}(\ParVT)$ is actually independent of $\delta$, its estimate $\widehat{\Mat{M}}_{n_0}$ may depend on it due to edge effects. 
\end{remark}

The following proposition provides a framework to study the asymptotic properties of (\ref{eq-estMn}) and shows the consistency of $\widehat{\Mat{M}}_{n_0}$ when the domain $\Lambda_{n}$ increases up to $\infty$ as $n\to\infty$. Its proof is relegated to section \ref{preuve-consistance}.


\begin{proposition}\label{consistance}
Under \textbf{[Mod]}, \textbf{[E1]}, assume that for any $\ParV$ in a neighborhood $\mathcal{V}(\ParVT)$ of $\ParVT$, for any bounded domain $\Lambda$, for any $\varphi\in \Omega$ and for $j=1,\ldots,p$, $\left(\widehat{\Vect{Y}}_{n,\Lambda}(\varphi;\cdot \right)_j$ is a continuous function. Assume moreover  that 
\begin{equation} \label{eq-hypYhat}
\sup_{k \in \mathcal{K}_n} \left|\widehat{\Vect{Y}}_{n,\Delta_k(\delta_n)}(\Phi;\ParV)- {\Vect{Y}}_{\Delta_k(\delta_n)} (\Phi;\ParV) \right| \stackrel{\mathbb{P}}{\rightarrow} 0,
\end{equation}
where, for any $\delta>0$ as above,  $(\delta_n)_{n\in\NN}$ is a sequence satisfying  $|\Lambda_{n}|\delta_n^{-d}\in\NN$, $\delta_{n_0}=\delta$ and $\delta_n\to\delta$ as $n\to\infty$.
Then, for any $D^\vee\geq D$,
$$\widehat{\Mat{M}}_n \left( \Phi; \widehat{\ParV}_n(\Phi), \delta_n, D^\vee\right)\overset{\mathbb{P}}{\longrightarrow}\Mat{M}(\ParVT).$$
\end{proposition}

\begin{remark}
The choice of the sequence $(\delta_n)_{n\in\NN}$ is always possible (see the proof). Since we allow the domain $\Lambda_n$ to grow continuously up to $\RR[d]$, its decomposition in sub-cubes with side-length $\delta$ is not always possible. The sequence  $(\delta_n)_{n\in\NN}$ is thus mandatory to make a decomposition of the domain available when $n$ increases. We chose it by respecting as most as possible the initial choice of the practicioner.
\end{remark}


\section{Goodness-of-fit tests for stationary marked Gibbs point processes} \label{sec-gof}
We present in this section three goodness-of-fit tests, based on the residuals computed according to the different frameworks considered in Section \ref{sec-asymp}.
We assume that the point process is observed in the domain $\Lambda_{n_0}\oplus D^+$ where $D^+\geq D$ and $\Lambda_{n_0}$ is a cube. 


\subsection{Quadrat-type test with $|J|-1$ degrees of freedom}\label{sec-test-quadrat}
According to the setting detailed in Section \ref{sec-fwk1}, we divide the domain $\Lambda_{n_0}$ into a fixed number of subdomains, namely $\Lambda_{n_0} := \bigcup_{j \in \mathcal J} \Lambda_{j,n_0}$ where $\mathcal J$ is a finite set and all the $\Lambda_{j,n_0}$ are disjoint cubes with the same volume $\left| \Lambda_{0,n_0}\right|$.
Moreover, in each sub-domain, we consider the decomposition $\Lambda_{j,n_0}=\cup_{k\in \mathcal K_{j,n_0}}\Delta_k(\delta)$, for any $\delta$ such that $|\Lambda_{0,n_0}|\delta^{-d}\in\NN$, where  the $\Delta_k(\delta)$'s are disjoint cubes with side-length $\delta$.

Following (\ref{eq-estMn}), we consider, for any $\delta>0$ as above and any $D^{\vee}\geq D$, the estimator
\begin{equation}\label{est-Inn}
\widehat{\lambda}_{n_0,Inn} = |\Lambda_{n_0}|^{-1} \;\sum_{i\in \mathcal K_{n_0}} \;\;\sum_{j \in \mathbbm{B}_i\left( \left\lceil \frac{D^\vee}{\delta}\right\rceil \right)\cap \mathcal K_{n_0}}  \!\!\!\!\!\!
I_{\Delta_i(\delta)}\left(\varphi; \widehat{\ParV}_{n_0}(\varphi) \right) I_{\Delta_j(\delta)}\left(\varphi; \widehat{\ParV}_{n_0}(\varphi) \right),\end{equation}
where $\mathcal K_{n_0}=\cup_{j\in \mathcal{J}} \mathcal K_{j,n_0}$. Note that  $I_{\Delta_i(\delta)}\left(\varphi; \widehat{\ParV}_{n_0}(\varphi) \right) =R_{\Delta_i(\delta)}\left(\varphi; \widehat{\ParV}_{n_0}(\varphi) \right) $ but we preserve this redundant notation in the sequel.

The following corollary is an immediate consequence of Corollary \ref{cor-quadrat} and Proposition \ref{consistance}. 

\begin{corollary}\label{cor-test-quadrat}
Under the assumptions of Proposition \ref{prop-fwk1} and if \textbf{[PD]} holds for $\Vect{Y}_{\Ltt} \left(\Phi;\ParVT\right)={I}_{\Ltt} \left(\Phi;\ParVT\right)$, then, for any $\delta>0$, one can construct a sequence $(\delta_n)_{n\in\NN}$  satisfying  $|\Lambda_{0,n}|\delta_n^{-d}\in\NN$, $\delta_{n_0}=\delta$ and $\delta_n\to\delta$, such that as $n\to+\infty$
\begin{equation}\label{test-quadrat}
T_{1,n}:=|\Lambda_{0,n}|^{-1} \;  \widehat\lambda_{n,Inn}^{-1}  \; \times \; \| 
 \Vect{R}_{\mathcal{J},n}(\Phi;h)- \overline{\Vect{R}}_{\mathcal{J},n}(\Phi;h)
\|^2 \stackrel{d}{\longrightarrow} \chi^2(|\mathcal{J}|-1). 
\end{equation}
\end{corollary}

This result leads to a goodness-of-fit test for $H_0: \Phi\sim P_{\ParVT}$ versus $H_1:\Phi \nsim P_{\ParVT}$.  Let us  briefly summarize the different steps to implement the test for a given asymptotic level $\alpha\in (0,1)$.
\begin{itemize}
\item {\bf Step 1} Consider a parametric model of a stationary marked Gibbs point process with finite range $D$ observed on the domain $\Lambda_{n_0}\oplus D^+$ with $D^+\geq D$.
\item {\bf Step 2} Choose an estimation method satisfying the assumptions \textbf{[E1]}, \textbf{[E2(bis)]} (for example the MPLE) and compute the estimate $\widehat\ParV_{n_0}$ on $\Lambda_{n_0}$.
\item {\bf Step 3} 
	\begin{itemize}
	\item[a)] Consider a test function $h$ (satisfying \textbf{[C1-2]}, \textbf{[N1-3]} and \textbf{[PD]}), divide $\Lambda_{n_0}$ into $|\mathcal{J}|$ cubes and compute the $h-$residuals on each different cube.
	\item[b)] Estimate $\lambda_{Inn}$ by (\ref{est-Inn}).
	\item[c)] Compute the test statistic $T_{1,n_0}$ involved in (\ref{test-quadrat}).
  	\end{itemize}
\item {\bf Step 4} Reject the model if $T_{1,n_0}(\varphi)>\chi^2_{1-\alpha}(|\mathcal{J}|-1)$. 
\end{itemize}

Let us note that in the particular case of a homogeneous Poisson point process with intensity $z$ and when considering the raw residuals ($h=1$), this test is exactly the Poisson dispersion test applied to the $|\mathcal J|$ quadrat counts, also called quadrat counting test, see  \cite{B-Dig03} for instance. Indeed, in this case, $\Vect{R}_{\mathcal{J},n}(\varphi;h)- \overline{\Vect{R}}_{\mathcal{J},n}(\varphi;h)$ is the vector of quadrat counts and $\lambda_{Inn}=z$. Considering $|\Lambda_{0,n}| \widehat{\lambda}_{n_0,Inn}$ as an estimation of the intensity on $\Lambda_{0,n}$, the statistic $T_{1,n}$ reduces to the ratio of the sum of squares of the quadrat counts over their estimated mean.  

\begin{remark} \label{rem-lInn}
The condition \textbf{[PD]} in Corollary~\ref{cor-test-quadrat} has to be verified with $\Vect{Y}_{\Ltt} \left(\Phi;\ParVT\right)={I}_{\Ltt} \left(\Phi;\ParVT\right)$ which is not so difficult (see Proposition~\ref{prop-lInn} for a general result). Indeed, contrarily to Corollary~\ref{cor-T1prime} and~\ref{cor-T2}, this condition does not depend on the form of the estimator $\widehat{\ParV}_n$. 
Moreover, as emphasized in Section~\ref{sec-fwk1}, the assumptions of Proposition~\ref{prop-fwk1} are satisfied for many models (this will be explored in details for exponential models in Section~\ref{sec-MCN}). This means that the proposed goodness-of-fit test based on~\eqref{test-quadrat} may be used for many models and many choices of function $h$.
\end{remark}

\begin{remark} \label{drawback}
The weakeness of this testing procedure (and the next ones) could be the estimation (\ref{est-Inn}) of $\lambda_{Inn}$  (and in general the estimator (\ref{eq-estMn})). The choice of the parameters $\delta$ and $D^{\vee}$ in (\ref{est-Inn}) is crucial. For instance, for fixed $n$,  in the extreme cases $\delta\to 0$ or $D^{\vee}\to\infty$, we get $\widehat{\lambda}_{n,Inn}\approx 0$. A careful simulation study should help for these choices. Another improvement could be to estimate $\lambda_{Inn}$ via Monte-Carlo methods.


\end{remark}

\subsection{Quadrat-type test with $|J|$ degrees of freedom}\label{sec-T1prime}
Under the same setting as above, assume moreover that \textbf{[PD]} holds for $\Vect{Y}_{\Ltt} \left(\varphi;\ParVT\right)=R_{\infty,\Ltt}( \varphi;h,\ParVT)$. Let us define the normalized residuals
$$
\widetilde{\Vect{R}}_{1,n_0} (\varphi;h):= \widehat{\lambda}_{n_0,Inn}^{-1/2} \Vect{R}_{\mathcal{J},n_0}(\varphi;h) +  \left(\widehat{\lambda}_{n_0,Res}^{-1/2} - \widehat{\lambda}_{n_0,Inn}^{-1/2}\right) \overline{\Vect{R}}_{\mathcal{J},n_0}(\varphi;h),
$$
where  $\widehat{\lambda}_{n_0,Inn}$ is defined in (\ref{est-Inn}) and  $\widehat{\lambda}_{n_0,Res}$ is an estimate of $\lambda_{Res}$ following  (\ref{eq-estMn}). When considering the MPLE, explicit formulas for $\widehat{\lambda}_{n_0,Res}$  are given in Section~\ref{sec-estMatMPLE}. It is easy to check that $\widetilde{\Vect{R}}_{1,n_0} (\varphi;h)  = \widehat{\Mat{\Sigma}_1}_{n_0}^{-1/2}  \; \Vect{R}_{\mathcal{J},n_0}(\varphi;h)$. Therefore the following corollary is deduced from Propositions~\ref{prop-fwk1} and \ref{consistance}.

\begin{corollary} \label{cor-T1prime}
Under the assumptions of Propositions \ref{prop-fwk1} and \ref{consistance}, assuming that  \textbf{[PD]} holds for $\Vect{Y}_{\Ltt} \left(\Phi;\ParVT\right)={I}_{\Ltt} \left(\Phi;\ParVT\right)$ and $\Vect{Y}_{\Ltt} \left(\Phi;\ParVT\right)=R_{\infty,\Ltt}( \Phi;h,\ParVT)$, then, for any $\delta>0$, one can construct a sequence $(\delta_n)_{n\in\NN}$ which satisfies $|\Lambda_{0,n}|\delta_n^{-d}\in\NN$, $\delta_{n_0}=\delta$ and $\delta_n\to\delta$ as $n\to\infty$, such that
 as $n\to+\infty$,
\begin{equation} 
\widetilde{T}_{1,n} (\Phi):= |\Lambda_{0,n}|^{1/2} \| \widetilde{\Vect{R}}_{1,n} (\Phi;h)\|^2 \stackrel{d}{\longrightarrow} \chi^2(|\mathcal{J}|) \label{eq-T1n} 
\end{equation}
\end{corollary}

A goodness-of-fit test with asymptotic size $\alpha\in (0,1)$ is deduced similarly as in the previous section. The steps to follow in practice are the same except that in {\bf Step 3} b), one has to estimate both $\lambda_{Inn}$ and $\lambda_{Res}$, and in {\bf Step 4} we reject the model if  $\widetilde{T}_{1,n_0}(\varphi)>\chi^2_{1-\alpha}(|\mathcal{J}|)$.

\begin{remark}
Let us emphasize that, with respect to Corollary~\ref{cor-test-quadrat}, Corollary~\ref{cor-T1prime} involves an additional more complex assumption: \textbf{[PD]} has to be satisfied for $\Vect{Y}_{\Ltt} \left(\Phi;\ParVT\right)=R_{\infty,\Ltt}( \Phi;h,\ParVT)$. This kind of assumption deeply depends  on the nature of the estimate $\widehat{\ParV}$. This problem is investigated in Proposition~\ref{prop-exemples} for particular examples. Furthermore, we show in Proposition \ref{prop-PDfails} that  $\lambda_{Res}=0$  occurs for many models and many choices of $h$ including the Poisson model when $h=1$. These two remarks underline the fact that the test relying on $\widetilde{T}_{1,n}$ is more restrictive than the previous one with $T_{1,n}$.

\end{remark}

\subsection{Empty space function type test}\label{sec-T2}

Let us consider the setting of section  \ref{sec-fwk2}, where $s$ different residuals are computed on the same full domain $\Lambda_{n_0}$. We consider the decomposition $\Lambda_{n_0}=\cup_{k\in \mathcal K_{n_0}}\Delta_k(\delta)$, for any $\delta$ such that $|\Lambda_{n_0}|\delta^{-d}\in\NN$, where  the $\Delta_k(\delta)$'s are disjoint cubes with side-length $\delta$.

Under the notation of Proposition~\ref{prop-fwk2}, assuming \textbf{[PD]} holds for $\Vect{Y}_{\Ltt} \left(\varphi;\ParVT\right)=\Vect{R}_{\infty,\Ltt}(\varphi,\Vect{h},\ParVT)$, let us define
\begin{equation*}
\widetilde{\Vect{R}}_{2,n_0} (\varphi;\Vect{h},\widehat{\ParV})  :=\widehat{\Mat{\Sigma}_2}_{n_0}^{-1/2}  \; \left( R_{\Lambda_{n_0}}(\varphi;{h_j},\widehat{\ParV})\right)_{j=1,\ldots,s}
\end{equation*}
where $\widehat{\Mat{\Sigma}_2}_{n_0}^{-1/2}:=\widehat{\Mat{\Sigma}_2}_{n_0}^{-1/2}(\varphi,\widehat{\ParV};\delta,D^\vee)$ is an estimation of $\Mat{\Sigma}_2(\ParVT)$ as in (\ref{eq-estMn}).  See explicit formulas in Section~\ref{sec-estMatMPLE} when considering the MPLE.


From Propositions~\ref{prop-fwk2} and \ref{consistance}, we get the following corollary.
\begin{corollary} \label{cor-T2}
Assuming \textbf{[PD]} with $\Vect{Y}_{\Ltt} \left(\varphi;\ParVT\right)=\Vect{R}_{\infty,\Ltt}(\varphi,\Vect{h},\ParVT)$, under the assumptions of Propositions \ref{prop-fwk2} and~\ref{consistance}, then, for any $\delta>0$ as above, one can construct a sequence $(\delta_n)_{n\in\NN}$ which satisfies $|\Lambda_{n}|\delta_n^{-d}\in\NN$, $\delta_{n_0}=\delta$ and $\delta_n\to\delta$ as $n\to\infty$, such that, as $n\to+\infty$,
\begin{equation}\label{eq-T2n}
\widetilde T_{2,n} (\Phi) := |\Lambda_{n}|^{1/2} \| \widetilde{\Vect{R}}_{2,n} (\Phi;\Vect{h},\widehat{\ParV})\|^2 \stackrel{d}{\rightarrow} \chi^2(s) .
\end{equation}
\end{corollary}


A goodness-of-fit test for $H_0: \Phi\sim P_{\ParVT}$ versus $H_1:\Phi \nsim P_{\ParVT}$, with asymptotic size $\alpha\in(0,1)$ is deduced as before. From a practical point of view, the steps detailed in \ref{sec-test-quadrat} are modified into:
\begin{itemize}
\item {\bf Step 3(framework 2)}
	\begin{itemize}
	\item[a)] Consider $s$ different test functions (satisfying \textbf{[C1-2]}, \textbf{[N1-3]} and \textbf{[PD]}), and compute the $s$ different $h_j-$residuals on the same initial domain $\Lambda_{n_0}$.
	\item[b)] Estimate the matrix $\Mat{\Sigma}_2(\ParVT)$ by~\eqref{eq-estMn} and compute $\widehat{\Mat{\Sigma}_2}_{n_0}^{-1/2}$ with any numerical routine (e.g. a choleski decomposition or a singular value decomposition).
	\item[c)] Compute the test statistic $\widetilde T_{2,n_0}(\varphi)$  defined by~\eqref{eq-T2n}.
	\end{itemize}
\item {\bf Step 4} Fix $\alpha \in (0,1)$ and reject the model if $\widetilde T_{2,n_0}(\varphi)>\chi^2_{1-\alpha}(s)$.
\end{itemize}


\section{Application to exponential models and the MPLE} \label{sec-disc}
Through Sections~\ref{sec-asymp}, \ref{sec-estMat} and \ref{sec-gof} three sets of assumptions have been considered. The first one deals with integrability and regularity of the model and the test function(s) and gathers \textbf{[Mod]}, \textbf{[C]} and \textbf{[N1-4]}. The second one is about the estimator $\widehat{\ParV}_n$ and involves \textbf{[E1]} and \textbf{[E2(bis)]}. Finally, the third one, assumption \textbf{[PD]} is very specific and deals with the positive definiteness of covariance matrices. We prove in this section that these assumptions are in general fulfilled for exponential family models and the MPLE.

\subsection{Assumptions \textbf{[Mod]}, \textbf{[C]} and \textbf{[N1-4]} for exponential family models} \label{sec-MCN}

The energy function of exponential family models is given for any $\Lambda\in\sB(\RR[d])$  by $V_{\Lambda}(\varphi;\ParV) = \tr{\ParV} \Vect v_{\Lambda} (\varphi)$,
where $\Vect v_{\Lambda}(\varphi)$ is the vector of sufficient statistics given by
$\Vect v_{\Lambda}(\varphi)=\tr{(v_{1,\Lambda}(\varphi),\ldots, v_{p,\Lambda}(\varphi))}$. The local energy is then expressed as $\VIPar{x^{m}}{\varphi}{\ParV}= \tr{\ParV} \SExIV{x^{m}}{\varphi}$, 
where $\SExIV{x^{m}}{\varphi}= (\SExI{1}{x^{m}}{\varphi}, \ldots, \SExI{p}{x^{m}}{\varphi}):=\Vect v_{\Lambda}(\varphi\cup \{x^{m}\}) - \Vect v_{\Lambda}(\varphi)$. Let us consider the following assumption:
\begin{itemize}
\item[\textbf{[Exp]}] For $i=1,\cdots,p$, there exist $\cSEx{\inf}{i},\cSEx{\sup}{i}\geq 0$, $k_i\in \NN$  such that one of both following assumptions is satisfied for all $(m,\varphi) \in \mSp\times \Omega$:
$$
\Par_i \geq 0 \mbox{ and } 
-\cSEx{\inf}{i}\leq \SExI{i}{0^{m}}{\varphi}=\SExI{i}{0^{m}} {\varphi_{\mathcal{B}(0,D)}} \leq \cSEx{\sup}{i} |\varphi_{\mathcal{B}(0,D)} |^{k_i}.
$$
or
$$
-\cSEx{\inf}{i}\leq \SExI{i}{0^{m}}{\varphi}=\SExI{i}{0^{m}} {\varphi_{\mathcal{B}(0,D)}} \leq \cSEx{\sup}{i}.
$$
\end{itemize}
The assumption \textbf{[Exp]} has already been considered by \cite{A-BilCoeDro08}. It is fulfilled for a large class of examples including the overlap area point process, the multi-Strauss marked point process, the $k-$nearest-neighbor multi-Strauss marked point process, the Strauss type disc process, the Geyer's triplet point process, the area interaction point process,\ldots. 

\begin{proposition} \label{prop-C12N13}
Under \textbf{[Exp]}, the assumptions \textbf{[Mod]}, \textbf{[C]} and \textbf{[N1-4]} are satisfied for the raw residuals, inverse residuals, Pearson residuals or residuals based on the empty space function. 
\end{proposition}

\begin{proof}
The assumption \textbf{[Exp]}  implies that the local energy function is local and stable, which, from results of \cite{A-BerBilDro99}, implies that \textbf{[Mod]} is fulfilled. A direct consequence of  \textbf{[Exp]}  is that for every $\alpha>0$, for all $\ParV \in \SpPar$ and for all $i=1,\ldots,p$
\begin{equation} \label{eq-propInteg}
\Esp\left( |\SExI{i}{0^{M}}{\Phi}|^\alpha
e^{ - \tr{\ParV} \SExIV{0^{M}}{\Phi} }
\right) <+ \infty,
\end{equation}
which ensures the integrability assumptions  \textbf{[C]} and \textbf{[N1-2]} for the residuals considered in the proposition.  The locality assumption \textbf{[N4]} is contained in \textbf{[Exp]}. Finally,  an application of the dominated convergence theorem, with the help of (\ref{eq-propInteg}), shows \textbf{[N3]}.
\end{proof}

\begin{remark} \label{rem-LJ}
Our setting is not restricted to locally stable exponential family  models. As an example,  following ideas of \cite{A-CoeDro09}, one may prove that  \textbf{[C]} and \textbf{[N1-4]}  are fulfilled  for Lennard-Jones type models. 
\end{remark}

\subsection{Assumptions \textbf{[E1]} and \textbf{[E2(bis)]} for the MPLE} \label{sec-hypEst}
Among the different parametric estimation methods available for spatial point patterns, the maximum pseudolikelihood is of particular interest. Indeed, unlike the maximum likelihood estimation method, it does not require the computation of the partition function, it is quite easy to implement and asymptotic results are now well-known (see \cite{A-JenMol91}, \cite{A-JenKun94}, \cite{A-BilCoeDro08}, \cite{A-DerLav09} and \cite{A-CoeDro09}). The MPLE is obtained by maximizing the log-pseudolikelihood contrast, given for exponential models by
\begin{equation} \label{eq-defLPL}
LPL_{\Lambda_n}(\varphi;\ParV)  = -\ism[\Lambda_n \times \mSp]{e^{-\tr{\ParV} \Vect{v}(x^m|\varphi)}} \; - \; \tr{\ParV} \; \sum_{x\in \varphi} \Vect{v}(x^m|\varphi\setminus x^m).
\end{equation}

\begin{proposition} \label{prop-MPLE} Under assumption \textbf{[Exp]} (and an additional indentifiability condition), \textbf{[E1]} and \textbf{[E2(bis)]} are fulfilled for the MPLE. The vector $\Vect{U}_\Lambda(\varphi;\ParVT)$ in \textbf{[E2(bis)]} is then expressed as follows
\begin{equation} \label{eq-UMPLE}
\Vect{U}_\Lambda(\varphi;\ParVT) := \Mat{H}(\ParVT)^{-1} \Vect{LPL}^{(1)}_\Lambda(\varphi;\ParVT) ,
\end{equation}
where $\Vect{LPL}^{(1)}_\Lambda(\varphi;\ParVT)$ is the gradient vector of the log-pseudolikelihood given by
\begin{equation} \label{eq-defLPL1}
 \Vect{LPL}_\Lambda^{(1)}(\varphi;\ParVT) : =\ism{\Vect{v}(x^m|\varphi;\ParVT) e^{-\tr{\ParVT}\Vect{v}(x^m|\varphi)} } - \sum_{x^m \in \varphi_\Lambda} \Vect{v}(x^m|\varphi\setminus x^m;\ParVT)
\end{equation}
and where $\Mat{H}(\ParVT)$ is the symmetric matrix  
\begin{equation} \label{eq-defH}
\Mat{H}(\ParVT) := \Esp\left( \Vect{v}(0^M|\Phi;\ParVT) \tr{\Vect{v}(0^M|\Phi;\ParVT)} 
\; e^{-\VIPar{0^M}{\Phi}{\ParVT}}
\right).
\end{equation}
\end{proposition}

\begin{proof} \textbf{[E1]} is proved by \cite{A-BilCoeDro08} (under \textbf{[Exp]} and the identifiability condition \textbf{[Ident]} arising p.244 in \cite{A-BilCoeDro08}). Let $\Vect{Z}_n(\varphi;\ParVT):=-  |\Lambda_n|^{-1} \Vect{LPL}_{\Lambda_n}(\varphi;\ParVT)$. If $\widehat{\ParV}_n(\varphi)=\widehat{\ParV}_n^{MPLE}(\varphi)$ denotes the maximum pseudolikelihood estimate, one derives
$$
\Vect{Z}_n^{(1)}(\varphi;\widehat{\ParV}_n)- \Vect{Z}_n^{(1)}(\varphi;\ParVT) =0 - \Vect{Z}_n^{(1)}(\varphi;\ParVT) =\Mat{H}_n(\varphi;\ParVT,\widehat{\ParV}_n)( \widehat{\ParV}_n(\varphi)-\ParVT)$$ with 
$\Mat{H}_n(\varphi;\ParVT,\widehat{\ParV}_n)= \int_0^t \Mat{Z}_n^{(2)} \left( \varphi;\ParVT + t(\widehat{\ParV}_n(\varphi)-\ParVT)\right)dt.$
Under assumptions \textbf{[Exp]} and \textbf{[Ident]}, then, for $n$ large enough, $\Mat{H}_n$ is invertible and converges almost surely towards the matrix $\Mat{H}(\ParVT)$ given by~\eqref{eq-defH}. Moreover, following the proof of Theorem~2 of \cite{A-BilCoeDro08} (see condition $(iii)$ p.257-258), we derive
$\Var (\Vect{Z}_n^{(1)}(\Phi;\ParVT ))=\mathcal{O}(|\Lambda_n|^{-1})$. So 
\begin{eqnarray*}
|\Lambda_n|^{1/2} \!\! \left( (\widehat{\ParV}_n(\Phi) -\ParVT )  + \Mat{H}^{-1}(\ParVT) \Vect{Z}_n^{(1)}(\Phi;\ParVT)\right) \!\!\!&=& \!\!\!- |\Lambda_n|^{1/2}\left( \Mat{H}_n^{-1}(\Phi;\widehat{\ParV}_n,\ParVT) - \Mat{H}^{-1}(\ParVT) \right) \Vect{Z}_n^{(1)}(\Phi;\ParVT) \\
 &\to& 0,
\end{eqnarray*}
in probability as $n \to +\infty$. This implies (\ref{eq-UMPLE}). Finally, $\Vect{U}_\Lambda(\varphi;\ParVT)$ fulfills properties $(i)-(v)$ in  \textbf{[E2(bis)]} for the same reasons as in the proof of Proposition \ref{prop-C12N13} and, for $(iii)$, from the proof of Theorem 2 (step 1, p. 257) in  \cite{A-BilCoeDro08}.
\end{proof}

\begin{remark}
In the same spirit as Remark~\ref{rem-LJ}, let us underline that the MPLE also satisfies \textbf{[E1]} and \textbf{[E2(bis)]} for some non locally stable and non exponential family models,  including Lennard-Jones type models (provided a locality assumption).
\end{remark}

\subsection{Estimation of asymptotic covariance matrices when considering the MPLE} \label{sec-estMatMPLE}

We still focus on exponential family models. As in Section~\ref{sec-defEstMat}, we assume that the point process is observed in the domain $\Lambda_{n_0}\oplus D^+$ where $D^+\geq D$ and $\Lambda_{n_0}$ is a cube. Moreover,  we consider the decomposition $\Lambda_{n_0}=\cup_{k\in \mathcal K_{n_0}}\Delta_k(\delta)$, for any $\delta$ such that $|\Lambda_{n_0}|\delta^{-d}\in\NN$, where  the $\Delta_k(\delta)$'s are disjoint cubes with side-length $\delta$ and centered, without loss of generality, at $k\delta$.

From (\ref{def-Rinfty}) and \eqref{eq-UMPLE}, we have under the assumptions \textbf{[Exp]} and when considering the MPLE 
\begin{equation} \label{eq-RinftyMPLE}
R_{\infty,\Lambda}(\varphi;h,\ParVT) := I_\Lambda(\varphi;h,\ParVT) - \tr{\Vect{LPL}^{(1)}(\varphi;\ParVT)} \Vect{W}(h,\ParVT)
\end{equation}
where $\Vect{W}(h,\ParV) := \Mat{H}(\ParV)^{-1} \Vect{\mathcal{E}}(h,\ParV)$. A natural estimator of $\Vect{W}(h,\ParVT)$ is given by $\widehat{\Vect{W}}_{n_0}(\varphi;h,\widehat{\ParV}_{n_0}):=\widehat{\Mat{H}}_{n_0}(\varphi;\widehat{\ParV}_{n_0})^{-1} \widehat{\Vect{\mathcal{E}}}_{n_0}(\varphi;h,\widehat{\ParV}_{n_0})$ where 
\begin{eqnarray*}
\widehat{\Mat{H}}_{n_0}(\varphi;\widehat{\ParV}_{n_0}) &=& |\Lambda_{n_0}|^{-1} \ism[\Lambda_{n_0} \times \mSp]{ \Vect{v}(x^m|\varphi) \tr{\Vect{v}(x^m|\varphi)}
e^{-\tr{\widehat{\ParV}_{n_0}}\Vect{v}(x^m|\varphi) }}, \\
\widehat{\Vect{\mathcal{E}}}_{n_0}(\varphi;h,\widehat{\ParV}_{n_0}) &=& |\Lambda_{n_0}|^{-1} \ism[\Lambda_{n_0} \times \mSp]{ 
h(x^m,\varphi; \widehat{\ParV}_{n_0})
\Vect{v}(x^m|\varphi) e^{-\tr{\widehat{\ParV}_{n_0}}\Vect{v}(x^m|\varphi) } }. \\
\end{eqnarray*}
In this spirit, let $\widehat{R}_{n_0,\infty,\Lambda}(\varphi;h,\widehat{\ParV}_{n_0}):= I_{\Lambda}(\varphi;h,\widehat{\ParV}_{n_0})-\tr{\Vect{LPL}^{(1)}_\Lambda(\varphi;\widehat{\ParV}_{n_0}) } \widehat{\Vect{W}}_{n_0}(\varphi;h,\widehat{\ParV}_{n_0})$ and \\ $\widehat{\Vect{R}}_{n_0,\infty,\Lambda}(\varphi;\Vect{h},\widehat{\ParV}_{n_0}):= \left( \widehat{R}_{n_0,\infty,\Lambda}(\varphi;h_j,\widehat{\ParV}_{n_0})\right)_{j=1,\ldots,s}$. Based on these notation, we obtain the following estimations for $\lambda_{Inn}, \lambda_{Res}$ and $\Mat{\Sigma}_2(\ParVT)$
\begin{eqnarray*}
\widehat{\lambda}_{n_0,Inn}(\varphi,\widehat{\ParV}_{n_0}(\varphi),\delta,D^\vee) &=& \!\!\!|\Lambda_{n_0}|^{-1} \;\sum_{i\in \mathcal K_{n_0}} \;\;\sum_{j \in \mathbbm{B}_i\left( \left\lceil \frac{D^\vee}{\delta}\right\rceil \right)\cap \mathcal K_{n_0}}  \!\!\!\!\!\!
I_{\Delta_i(\delta)}\left(\varphi; \widehat{\ParV}_{n_0}(\varphi) \right) I_{\Delta_j(\delta)}\left(\varphi; \widehat{\ParV}_{n_0}(\varphi) \right), \\
\widehat{\lambda}_{n_0,Res}(\varphi,\widehat{\ParV}_{n_0}(\varphi),\delta,D^\vee) &=&\!\!\! |\Lambda_{n_0}|^{-1} \;\sum_{i\in \mathcal K_{n_0}} \;\;\sum_{j \in \mathbbm{B}_i\left( \left\lceil \frac{D^\vee}{\delta}\right\rceil \right)\cap \mathcal K_{n_0}} \!\!\!\!\!\!
\widehat{R}_{\infty,\Delta_i(\delta)}(\varphi;h,\widehat{\ParV}_{n_0}) 
\widehat{R}_{\infty,\Delta_j(\delta)}(\Phi;h,\widehat{\ParV}_{n_0}), \\ 
\widehat{\Mat{\Sigma}_2}_{n_0}(\varphi,\widehat{\ParV}_{n_0}(\varphi),\delta,D^\vee) &=& \!\!\!|\Lambda_{n_0}|^{-1} \;\sum_{i\in \mathcal K_{n_0}} \;\;\sum_{j \in \mathbbm{B}_i\left( \left\lceil \frac{D^\vee}{\delta}\right\rceil \right)\cap \mathcal K_{n_0}} \!\!\!\!\!\!
\widehat{\Vect{R}}_{\infty,\Delta_i(\delta)}(\varphi;\Vect{h},\widehat{\ParV}_{n_0}) 
\tr{ \widehat{\Vect{R}}_{\infty,\Delta_j(\delta)}(\varphi;\Vect{h},\widehat{\ParV}_{n_0})}.
\end{eqnarray*}

\begin{corollary} \label{cor-PD} Under the notation and assumptions of Propositions~\ref{prop-fwk1} and~\ref{prop-fwk2}, and under  \textbf{[Exp]}, then, for any $\delta>0$ as above, one can consider a sequence $\delta_n$  which satisfies $\delta_{n_0}=\delta$ and $\delta_n\to\delta$, such that for any $D^\vee\geq D$,  the estimators
$\widehat{\lambda}_{n,Inn}(\Phi,\widehat{\ParV}_n(\Phi),\delta_n,D^\vee)$, $\widehat{\lambda}_{n,Res}(\Phi,\widehat{\ParV}_n(\Phi),\delta_n,D^\vee)$ and 
$\widehat{\Mat{\Sigma}_2}_n(\Phi,\widehat{\ParV}_n(\Phi),\delta_n,D^\vee)$ converge in probability (as $n\to +\infty$) towards respectively $\lambda_{Inn}$, $\lambda_{Res}$ and $\Mat{\Sigma}_2(\ParVT)$.
\end{corollary}

\begin{proof}
We apply Proposition \ref{consistance}, where for any $\ParV \in \mathcal{V}(\ParVT)$, we set
\begin{itemize}
\item for $\lambda_{Inn}$: $\widehat{\Vect{Y}}_\Lambda(\varphi;\ParV)=\Vect{Y}_\Lambda(\varphi;\ParV)= I_\Lambda(\varphi;h,\ParV)$.
\item for $\lambda_{Res}$: $\Vect{Y}_\Lambda(\varphi;\ParV)=R_{\infty,\Lambda}(\varphi;h,\ParV)$ and $\widehat{\Vect{Y}}_{n,\Lambda}(\varphi;\ParV) = \widehat{R}_{n,\infty,\Lambda}(\varphi;h,\ParV)$.
\item for $\Mat{\Sigma}_2(\ParVT)$: $\Vect{Y}_\Lambda(\varphi;\ParV)=\Vect{R}_{\infty,\Lambda}(\varphi;\Vect{h},\ParV)$ and  $\widehat{\Vect{Y}}_{n,\Lambda}(\varphi;\ParV)=\widehat{\Vect{R}}_{n,\infty,\Lambda}(\varphi;\Vect{h},\ParV)$.
\end{itemize}
The result is obvious for $\lambda_{Inn}$. For $\lambda_{Res}$ (the proof is similar for $\Mat{\Sigma}_2(\ParVT)$), it remains to prove that  for any $\ParV \in \SpPar$, $\sup_{k \in \mathcal{K}_n}\left| \widehat{R}_{n,\infty,\Delta_k(\delta_n)}(\Phi;h,\ParV) -{R}_{n,\infty,\Delta_k(\delta_n)}(\Phi;h,\ParV)\right| \to 0$ in probability as $n\to +\infty$. For any $k \in \mathcal K_n$, we derive
$$
\widehat{R}_{n,\infty,\Delta_k(\delta_n)}(\Phi;h,\ParV) - R_{\infty,\Delta_k(\delta_n)}(\Phi;h,\ParV) = \tr{\Vect{LPL}^{(1)}_{\Delta_k(\delta_n)}(\varphi;\ParV)} \left( \Vect{W}(h,\ParV)-\widehat{\Vect{W}}_n(\varphi;h,\ParV) \right).
$$
Assumption \textbf{[E2(bis)]} (implied by \textbf{[Exp]}, see Proposition \ref{prop-MPLE}) ensures that $|\Vect{LPL}^{(1)}_{\Delta_k(\delta_n)\setminus \Delta_k(\delta)}(\Phi;\ParV)|$ converges to 0 in quadratic mean. In particular, the convergence of  $\Vect{LPL}^{(1)}_{\Delta_k(\delta_n)}(\Phi;\ParV)$ towards $\Vect{LPL}^{(1)}_{\Delta_k(\delta)}(\Phi;\ParV)$ holds in probability. Moreover under the assumptions \textbf{[N1]} and \textbf{[E2(bis)]}, the ergodic theorem of \cite{A-NguZes79} may be applied to prove that $\widehat{\Vect{W}}_n(\Phi;h,\ParV)$ converges almost surely towards $\Vect{W}(h,\ParV)$, as $n\to +\infty$. Slutsky's theorem ends the proof.
\end{proof}

\begin{remark}
If the model is not an exponential model, Corollary~\ref{cor-PD} still holds by replacing the vector of the sufficient statistics, $\Vect{v}(x^m|\varphi)$, by the gradient vector of the local energy function, $\Vect{V}^{(1)}(x^m|\varphi)$ in the different definitions.
\end{remark}

\subsection{Positive definiteness of covariance matrices when considering the MPLE} \label{sec-PDMatMPLE}

Let us now focus on the positive-definitess of the above quantities. According to Proposition \ref{prop-sdp}  the key assumption to check is \textbf{[PD]}.

As adressed in Remark~\ref{rem-lInn}, we begin by giving a general result ensuring that $\lambda_{Inn}>0$.

\begin{proposition} \label{prop-lInn}
Under the assumption \textbf{[Exp]}, then $\lambda_{Inn}>0$ for the raw residuals, the Pearson residuals and the inverse residuals.
\end{proposition}

\begin{proof}
In \textbf{[PD]}, we fix $\overline{\delta}=D$ and $B=\emptyset$. Let us  write $\overline{\Omega}:=\overline{\Omega}_{\emptyset}$. Consider the following events for some $n\geq 1$
$$
A_0= \left\{ \varphi \in \overline{\Omega}: \varphi(\Delta_0(\overline{\delta}) )=0 \right\} \quad \mbox{ and } \quad
A_n = \left\{ \varphi \in \overline{\Omega}: \varphi(\Delta_0(\overline{\delta}))=n \right\},
$$
and let $\varphi_0\in A_0$ and $\varphi_n\in A_n$. Recall that the local stability property (ensured by \textbf{[Exp]})  asserts that there exists $K\geq 0$ such that $\VIPar{x^m}{\varphi}{\ParVT} \geq -K$ for any $x^m\in \mathbb{S}$ and any $\varphi\in \Omega$. Now, let us consider the three type of residuals.\\
{\it Raw residuals $(h=1)$.} From the local stability property
$$
| I_{\Ltt}(\varphi_n;h,\ParVT) - I_{\Ltt}(\varphi_0;h,\ParVT) | \geq n - \left| \ism[\Ltt\times \mSp]{ e^{-\VIPar{x^m}{\varphi_n}{\ParVT}}-e^{-\VIPar{x^m}{\varphi_0}{\ParVT}} }\right| \geq n - 2|\Ltt| e^{K}>0,
$$
for $n$ large enough.  And so assuming that the left-hand-side is zero leads to a contradiction, which proves \textbf{[PD]}.\\
{\it Inverse residuals ($h=e^V$).} Again, from the local stability property
$$
\left| I_{\Ltt}(\varphi_n;h,\ParVT) - I_{\Ltt}(\varphi_0;h,\ParVT) \right| = \left|\sum_{x^m\in {\varphi_n}_{\Ltt}} e^{\VIPar{x^m}{\varphi_n\setminus x^m}{\ParVT}} \right|\geq n e^{-K}>0,
$$
which proves \textbf{[PD]} similarly to the previous case. \\
{\it Pearson residuals ($h=e^{V/2}$).} From the same argument
\begin{eqnarray*}
\left| I_{\Ltt}(\varphi_n;h,\ParVT) - I_{\Ltt}(\varphi_0;h,\ParVT) \right| &\geq& \bigg|\sum_{x^m\in {\varphi_n}_{\Ltt}} e^{\VIPar{x^m}{\varphi_n\setminus x^m}{\ParVT}/2} \bigg| - \\
&&\left| \ism[\Ltt\times \mSp]{ e^{-\VIPar{x^m}{\varphi_n}{\ParVT}/2}-e^{-\VIPar{x^m}{\varphi_0}{\ParVT}/2} }\right| \\
&\geq & n e^{-K/2} - 2 |\Ltt| e^{K/2} >0,
\end{eqnarray*}
for $n$ large enough, which ends the proof.
\end{proof}

Proposition~\ref{prop-lInn} asserts that \textbf{[PD]} is fullfilled for $\Vect{Y}_{\Ltt} \left(\Phi;\ParVT\right)={I}_{\Ltt} \left(\Phi;\ParVT\right)$. Therefore, the combination of Propositions~\ref{prop-C12N13},~\ref{prop-MPLE} and~\ref{prop-lInn} and Corollary~\ref{cor-PD} ensures all the conditions of Corollary~\ref{cor-test-quadrat} hold. So a goodness-of-fit test based on~\eqref{test-quadrat} is valid for exponential family models satisfying \textbf{[Exp]} and for the raw residuals, the Pearson residuals and the inverse ones.

Now, let us focus on tests based on Corollary~\ref{cor-T1prime} and~\ref{cor-T2}. The following result is important from a practical point of view. It  asserts that $\lambda_{Res}$ (and so $\Mat{\Sigma}_1(\ParVT)$), and  $\Mat{\Sigma}_2(\ParVT)$ may fail to be positive-definite for an inappropriate choice of test function.

\begin{proposition} \label{prop-PDfails}
Let us consider an exponential family model, let $\widehat{\ParV}:=\widehat{\ParV}^{MPLE}$ and let us choose a test function of the form $h(x^m,\varphi;\ParV)= \tr{\Vect{\omega}}  \Vect{v}(x^m|\varphi)$ for some $\Vect{\omega}\in \RR[p]\setminus 0$, then $\lambda_{Res}=0$ and the matrices $\Mat{\Sigma}_1(\ParVT)$ and $\Mat{\Sigma}_2(\ParVT)$ in Propositions~\ref{prop-fwk1} and~\ref{prop-fwk2}  are only semidefinite-positive matrices.
\end{proposition}

\begin{proof}
The result is proved by noticing that 
\begin{eqnarray*}
\Mat{H}(\ParVT) \; \Vect{\omega} &=& \Esp \left( \Vect{v} (0^M|\Phi) \tr{\Vect{v} (0^M|\Phi)} e^{-\VIPar{0^M}{\Phi}{\ParVT}} \right) \; \Vect{\omega} \\ 
&=& \Esp \left(  \Vect{v} (0^M|\Phi) \tr{ \left( \tr{\Vect{\omega}}\Vect{v} (0^M|\Phi) \right)} e^{-\VIPar{0^M}{\Phi}{\ParVT}} \right) \\
&=& \Esp \left( h(0^M,\Phi;\ParVT) \Vect{v} (0^M|\Phi)  e^{-\VIPar{0^M}{\Phi}{\ParVT}} \right) =\Vect{\mathcal E}( \tr{\Vect{\omega}}  \Vect{v}, \ParVT).
\end{eqnarray*}
Therefore, $\Vect{W}\left( \tr{\Vect{\omega}}  \Vect{v}, \ParVT\right)= \Mat{H}(\ParVT)^{-1} \Vect{\mathcal{E}}(\tr{\omega} \Vect{v}, \ParVT)=
\Vect{\omega}$, which means that for any $\varphi\in \Omega$ and any bounded domain $\Lambda$
$$
R_{\infty,\Lambda}(\varphi;\tr{\Vect{\omega}}\Vect{v},\ParVT) = I_\Lambda(\varphi;\tr{\Vect{\omega}}\Vect{v},\ParVT) - \tr{\Vect{LPL}^{(1)}_{\Lambda}(\varphi;\ParVT) }\Vect{\omega} =0.
$$
This means that if,  for the framework~1, the test function is of the form $h=\tr{\Vect{\omega}}\Vect{v}$ then $\lambda_{Res}=0$ and
if one of the test functions, for the framework~2,  is of the form $h=\tr{\Vect{\omega}}\Vect{v}$, then $\Mat{\Sigma}_2(\ParVT)$ is necessary singular.
\end{proof}

\begin{remark}
As for Corollary~\ref{cor-PD}, the result of Proposition~\ref{prop-PDfails} still holds in general by replacing the vector  $\Vect{v}(x^m|\varphi)$ by the gradient vector of the local energy function $\Vect{V}^{(1)}(x^m|\varphi)$.
\end{remark}

As a consequence of Proposition~\ref{prop-PDfails},  the two goodness-of-fit tests based on $T'_{1,n}$ and $T_{2,n}$  in Section~\ref{sec-T1prime} and \ref{sec-T2} are not available (for the MPLE) if the test function $h$ is a linear combination of the sufficient statistics $\Vect{v}(x^m|\varphi)$. Since for most classical models, the value 1 can be obtained from a linear combination of $\Vect{v}(x^m|\varphi)$, the raw residuals ($h=1$) are not an appropriate choice for these two tests. This is the case for the  two following examples: the area-interaction point process and the 2-type marked Strauss point process, which are presented in details in Appendix \ref{exemplesPD}. The following result proves that for a different choice of $h-$residuals,  $\Mat{\Sigma}_1(\ParVT)$ and  $\Mat{\Sigma}_2(\ParVT)$ are positive-definite.

\begin{proposition}\label{prop-exemples}
For the 2-type marked Strauss point process and  the area-interaction point process, when considering the MPLE as an estimator of $\ParVT$, then
\begin{itemize}
\item the matrix $\Mat{\Sigma}_1(\ParVT)$ obtained  in Framework 1 from the inverse residuals $h=e^V$, 
\item the matrix $\Mat{\Sigma}_2(\ParVT)$ obtained  in Framework 2 from the empty space residuals, which are constructed for $0<r_1<\ldots<r_s<+\infty$ from the family of test functions
 $$h_j(x^m,\varphi;\ParV)= \mathbf{1}_{[0,r_j]}(d(x^m,\varphi)) e^{V\left(x^m|\varphi;\ParV\right)},\quad j=1,\ldots,s,$$
\end{itemize}
are positive-definite.
\end{proposition}
The proof of this result is postponed in Appendix~\ref{exemplesPD}. The combination of Propositions~\ref{prop-C12N13},~\ref{prop-MPLE},~\ref{prop-exemples} and Corollary~\ref{cor-PD} ensures all the conditions of Corollary~\ref{cor-T1prime} and~\ref{cor-T2} hold. So a goodness-of-fit test based on~\eqref{eq-T1n} (resp.~\eqref{eq-T2n}) is valid for the 2-type marked Strauss point process and  the area-interaction point process and for the inverse residuals (resp. the family of test functions based on the empty space function). 

Following the Proof of Proposition~\ref{prop-exemples}, it is the belief of the authors that such a result holds for other models and other choices of test functions. However, another model and/or test functions will lead to a specific proof. Therefore, this result cannot be as general as the one presented in Proposition~\ref{prop-lInn}.

\section{The non-hereditary case}\label{non-hereditary}
Up to here, we have assumed through \textbf{[Mod-E]} that the  family of energies is hereditary. We consider in this section the non-hereditary case. This particular situation can only occur in presence of a hardcore interaction. From a general point of view, we say that a family of energies involves a hardcore interaction if some point configurations have an infinite energy. Many classical models of Gibbs measures include a hardcore part, as the hard ball model.

A family of energies involving a hardcore part is hereditary if (\ref{heredite}) holds. This is a common assumption done for Gibbs energies and it appears to be fulfilled in most classical models, including the hard ball model. However, one may encounter some non-hereditary models, in the sense that  (\ref{heredite}) does not hold. Intuitively, in this case, when one removes a point from an allowed point configuration, it is possible to obtain a forbidden point configuration. This occurs for instance for Gibbs Delaunay-Voronoï tessellations or  forced-clustering processes (see \cite{A-DerLav09} and \cite{A-DerLav10}).

In the non-hereditary case, the GNZ formula (\ref{GNZnonstat}), which is the basis to define the residuals, becomes false (see Remark 2 in \cite{A-DerLav09}). It is extended to non-hereditary interactions in \cite{A-DerLav09}. This generalization requires to introduce the notion of removable points.

\begin{definition}
Let $\varphi\in\Omega$ and $x\in\varphi$, then $x$ is removable from $\varphi$ if there exists $\Lambda\in\sB(\RR[d])$ such that $x\in\Lambda$ and $V_{\Lambda}(\varphi-x; \ParV)<\infty$. The set of removable points in $\varphi$ is denoted by $\rem(\varphi)$.
\end{definition}

Notice that in the hereditary case, $\rem(\varphi)=\varphi$.

The GNZ formula is then generalized to the non-hereditary case as follows. Assuming for any $\ParV \in \SpPar$  that a Gibbs measure exists for the family of energies  $(V_{\Lambda}(.;\ParV))_{\Lambda\in\sB(\RR[d])}$, then, for any function $h(\cdot,\cdot,\ParV): \sSp\times \Omega\to \RR$  such that the following quantities are finite, 
\begin{equation}\label{GNZnon-hereditary}
\Esp\left( \ism[ \mathbb{R}^d \times \mSp]{h\left(x^m,\Phi;\ParV\right) e^{- \VIPar{x^m}{\Phi}{\ParVT}}} \right) = 
\Esp\left( \sum_{x^m \in \rem(\Phi) } h\left(x^m,\Phi\setminus x^m;\ParV\right)  \right).
\end{equation}

We can therefore define the $h-$residuals for (possibly) non-hereditary interactions.  For any bounded domain $\Lambda$, if  $\widehat{\ParV}$ is an estimate of $\ParVT$, the $h-$residuals are 
\begin{equation}\label{resnon-hereditary}
R_{\Lambda}\left( \varphi;h,\widehat{\ParV} \right) =  \ism[\Lambda\times \mSp]{h \left( x^m , \varphi ; \widehat{\ParV} \right) e^{-\VIPar{x^m}{\varphi}{\widehat{\ParV}}}} - \sum_{x^m\in \rem(\varphi_\Lambda)} h \left( x^m , \varphi\setminus x^m ; \widehat{\ParV} \right).\end{equation}

If the set of  removable points $\rem(\varphi)$ does not depend on $\ParV$, it is straightforward to extend all the asymptotic results obtained for the residuals in the preceding sections  to (\ref{resnon-hereditary}).

If the set of removable points depends on $\ParV$, this is false. Even in the hereditary case, if $\ParV$ is a hardcore parameter (as the hardcore distance in the hard ball model) then $\widehat{\ParV}$ behaves as an estimator of the support of the distribution $P_{\ParV}$. In this case assumption \textbf{[E2]} has typically few chances to hold and the asymptotic law of the residuals is unknown. In \cite{A-DerLav10} Figure 15, some simulations of raw-residuals for Gibbs Voronoï tessellations are presented, involving an estimated hardcore parameter in a non-hereditary setting : they show that the distribution of the residuals does not seem to be gaussian in this case.


\section{Proofs} \label{sec-proofs}

Since any stationary Gibbs measure can be represented as a mixture of ergodic measures (see \cite{B-Pre76}), it is sufficient to prove the different convergences involved in this paper for ergodic measures. We therefore assume from now on that  $P_{\ParVT}$ is ergodic.


\subsection{Proof of Proposition~\ref{prop-cons}}

(a) Under \textbf{[C1]}, the ergodic theorem of \cite{A-NguZes79} holds for both terms appearing in the definition of~$I_{\tilde\Lambda_n}\left( \varphi;h,\ParVT \right)$. Then, as $n\to +\infty$, one has $P_{\ParVT}-$a.s.
$$
|\tilde\Lambda_n|^{-1} I_{\tilde\Lambda_n}\left( \Phi;h,\ParVT \right) \rightarrow
\Esp\left( h \left( 0^M , \Phi ; \ParVT \right) e^{- V^{}\left( 0^M|\Phi; \ParVT \right)} 
\right) -\Esp\left( h \left( 0^M , \Phi\setminus 0^M ; \ParVT \right)  
\right), 
$$
which equals to 0 from the GNZ formula~(\ref{GNZ}).\\
(b) The aim is to prove that the difference $|\tilde\Lambda_n|^{-1}R_{\tilde\Lambda_n}\left( \varphi;h,\widehat{\ParV}_n(\varphi) \right)-|\tilde\Lambda_n|^{-1}I_{\tilde\Lambda_n}\left( \varphi;h,\ParVT \right)$ converges towards 0 for $P_{\ParVT}-$a.e. $\varphi$. Let us write 
$$
R_{\tilde\Lambda_n}\left( \varphi;h,\widehat{\ParV}_n(\varphi) \right)-I_{\tilde\Lambda_n}\left( \Phi;h,\ParVT \right):=T_1(\varphi)-T_2(\varphi)
$$
with 
\begin{eqnarray}
T_1(\varphi)&:=&  \ism[\tilde\Lambda_n\times \mSp]{\left(
f \left( x^m , \varphi ; \widehat{\ParV}_n(\varphi) \right) - f \left( x^m , \varphi ; \ParVT \right)
\right)} \label{eq-defT1}
\\
T_2(\varphi) &:=&   \sum_{x^m \in \varphi_{\tilde\Lambda_n}}
h \left( x^m , \varphi\setminus x^m ; \widehat{\ParV}_n(\varphi) \right) - h \left( x^m , \varphi\setminus x^m ; \ParVT \right). \label{eq-defT2}
\end{eqnarray}
Under the Assumptions~\textbf{[C2]} and \textbf{[E1]}, from the ergodic theorem and the GNZ formula, there exists $n_0\in \NN$ such that for all $n\geq n_0$
\begin{eqnarray}
|\tilde\Lambda_n|^{-1}T_1(\varphi) &\leq& \frac2{|\tilde\Lambda_n|} \ism[\tilde\Lambda_n\times \mSp]{ \tr{\left(\widehat{\ParV}_n(\varphi)-\ParVT \right)}
 \Vect{f}^{(1)}_{} \left( x^m , \varphi ;   \ParVT \right) } \nonumber \\
&\leq & 2 \|\widehat{\ParV}_n(\varphi)-\ParVT \| \times \frac1{|\tilde\Lambda_n|}\ism[\tilde\Lambda_n\times \mSp]{\|  \Vect{f}^{(1)}_{} \left( x^m , \varphi ;   \ParVT \right) \|}  \nonumber\\
&\leq&  4 \|\widehat{\ParV}_n(\varphi)-\ParVT \| \times \Esp\left( \| \Vect{f}^{(1)}_{} \left( 0^M , \Phi ;   \ParVT \right)\|\right) , \label{eq-T1}
\end{eqnarray}
and 
\begin{eqnarray}
|\tilde\Lambda_n|^{-1}T_2(\varphi) &\leq& \frac2{|\tilde\Lambda_n|} \sum_{x^m \in \varphi_{\tilde\Lambda_n}}  \tr{\left(\widehat{\ParV}_n(\varphi)-\ParVT \right)} \Vect{h}^{(1)}_{} \left( x^m , \varphi \setminus x^m ;   \ParVT \right) \nonumber \\
&\leq&  4 \|\widehat{\ParV}_n(\varphi)-\ParVT \| \times \Esp\left( \| \Vect{h}^{(1)}_{} \left( 0^M , \Phi ;   \ParVT \right)\|e^{- V^{}\left( 0^M|\Phi; \ParVT \right)}\right). \label{eq-T2}
\end{eqnarray}
Equations~(\ref{eq-T1}) and~(\ref{eq-T2}) lead to
$$
|\tilde\Lambda_n|^{-1}R_{\tilde\Lambda_n}\left( \varphi;h,\widehat{\ParV}_n(\varphi) \right)-|\tilde\Lambda_n|^{-1} I_{\tilde\Lambda_n}\left( \varphi;h,\ParVT \right) \leq c \|\widehat{\ParV}_n(\varphi)-\ParVT \|,
$$
for $n$ large enough, with $c=4 \times \Esp\left(\| \Vect{f}^{(1)}_{} \left( 0^M , \Phi ;   \ParVT \right)\|+\| \Vect{h}^{(1)}_{} \left( 0^M , \Phi ;   \ParVT \right)\|e^{- V^{}\left( 0^M|\Phi; \ParVT \right)}\right)$.


\subsection{Proof of Proposition~\ref{prop-equiv}}

Recall that 
$$ R_{\tilde\Lambda_n}\left( \varphi;h,\widehat{\ParV}_n(\varphi) \right) - I_{\tilde\Lambda_n}\left( \varphi;h,\ParVT \right)= T_1(\varphi)-T_2(\varphi)$$
 where $T_1(\varphi)$ and $T_2(\varphi)$ are defined by~(\ref{eq-defT1}) and~(\ref{eq-defT2}). Let us write
\begin{eqnarray*}
T_1(\varphi) &=& \ism[\tilde\Lambda_n \times \mSp]{\tr{\left(\widehat{\ParV}_n(\varphi)-\ParVT\right)}  \Vect{f}^{(1)}_{} \left( x^m , \varphi ;   \ParVT \right)} + T^\prime_1(\varphi) \\
T_2(\varphi) &=& \sum_{x^m \in \varphi_{\tilde\Lambda_n}} \tr{\left(\widehat{\ParV}_n(\varphi)-\ParVT\right)}  \Vect{h}^{(1)}_{} \left( x^m , \varphi\setminus x^m ;   \ParVT \right) + T^\prime_2(\varphi),
\end{eqnarray*}
with 
\begin{eqnarray*}
T_1^\prime(\varphi) &:=&  \ism[\tilde\Lambda_n \times \mSp]{A_1 \left( x^m , \varphi ; \widehat{\ParV}_n(\varphi) \right)}\\
T_2^\prime(\varphi) &=&  \sum_{x^m \in \varphi_{\tilde\Lambda_n}}  A_2 \left( x^m , \varphi\setminus x^m ; \widehat{\ParV}_n(\varphi) \right)
\end{eqnarray*}
and
\begin{eqnarray*}
A_1 \left( x^m , \varphi ; \widehat{\ParV}_n(\varphi) \right) &:=&  f \left( x^m , \varphi ; \widehat{\ParV}_n \right)- f \left( x^m , \varphi ; \ParVT \right) - \tr{\left(\widehat{\ParV}_n(\varphi)-\ParVT\right)}  \Vect{f}^{(1)}_{} \left( x^m , \varphi ;   \ParVT \right)\\
A_2 \left( x^m , \varphi ; \widehat{\ParV}_n(\varphi) \right) &:=&h \left( x^m , \varphi ; \widehat{\ParV}_n \right)- h \left( x^m , \varphi ; \ParVT \right) - \tr{\left(\widehat{\ParV}_n(\varphi)-\ParVT\right)}  \Vect{h}^{(1)}_{} \left( x^m , \varphi ;   \ParVT \right).
\end{eqnarray*}
From the mean value theorem, there exist for $j=1,\ldots,p$,\\ $\mathbf{\xi}_{1,j},\mathbf{\xi}_{2,j} \in [\min(\widehat{\theta}_1,\theta_1^\star),\max(\widehat{\theta}_1,\theta_1^\star)]\times \ldots\times
[\min(\widehat{\theta}_p,\theta_p^\star),\max(\widehat{\theta}_p,\theta_p^\star)]$ such that
\begin{eqnarray}
A_1 \left( x^m , \varphi ; \widehat{\ParV}_n(\varphi) \right) \!\!\!&=&\!\!\! \sum_{j=1}^p (\widehat{\theta}_j -\theta_j^\star) \left( 
 {f}^{(1)}_{j} \left( x^m , \varphi ;   \mathbf{\xi}_{1,j} \right)- {f}^{(1)}_{j} \left( x^m , \varphi ;   \ParVT \right)\right) \label{eq-res1}\\
\; A_2 \left( x^m , \varphi\setminus x^m ; \widehat{\ParV}_n(\varphi) \right) \!\!\!&=&\!\!\! \sum_{j=1}^p (\widehat{\theta}_j -\theta_j^\star) \left( 
 {h}^{(1)}_{j} \left( x^m , \varphi\setminus x^m ;   \mathbf{\xi}_{2,j} \right)- {h}^{(1)}_{j} \left( x^m , \varphi\setminus x^m ;   \ParVT \right)\right). \label{eq-res2}
\end{eqnarray}
Let $j \in \{1,\ldots,p\}$, again from the mean value theorem, there exist for $\ell=1,2$ and for $k=1,\ldots,p$, $\mathbf{\eta}_{\ell,j,k} \in [\min(\mathbf{\xi}_{\ell,j,1},\theta_1^\star),\max(\mathbf{\xi}_{\ell,j,1},\theta_1^\star)]\times \ldots\times[\min(\mathbf{\xi}_{\ell,j,p},\theta_p^\star),\max(\mathbf{\xi}_{\ell,j,p},\theta_p^\star)]
$ such that
\begin{eqnarray}
 {f}^{(1)}_{j} \left( x^m , \varphi ;   \mathbf{\xi}_{1,j} \right)- {f}^{(1)}_{j} \left( x^m , \varphi ;   \ParVT \right) &=& \sum_{k=1}^p \left( \mathbf{\xi}_{1,j,k}-\theta_k^\star\right)  f^{(2)}_{jk} \left( x^m , \varphi ;  \mathbf{\eta}_{1,j,k} \right) \label{eq-res3}\\
\quad {h}^{(1)}_{j} \left( x^m , \varphi\setminus x^m ;   \mathbf{\xi}_{2,j} \right)-{h}^{(1)}_{j} \left( x^m , \varphi\setminus x^m ;   \ParVT \right) &=& \sum_{k=1}^p \left( \mathbf{\xi}_{2,j,k}-\theta_k^\star\right)  h^{(2)}_{jk} \left( x^m , \varphi\setminus x^m ;  \mathbf{\eta}_{2,j,k} \right). \label{eq-res4}
\end{eqnarray}
By combining~(\ref{eq-res1}), (\ref{eq-res2}), (\ref{eq-res3}) and~(\ref{eq-res4}) and under \textbf{[N1]}, we can deduce the existence of $n_0 \in \NN$ such that for all $n\geq n_0$, one has for $P_{\ParVT}-$a.e. $\varphi$
\begin{eqnarray*}
|\tilde\Lambda_n|^{-1} |T_1^\prime(\varphi)| &\leq & \frac2{|\tilde\Lambda_n|} \ism[\tilde\Lambda_n \times \mSp]{
\sum_{j,k} \left| (\widehat{\theta}_j-\theta_j^\star)(\widehat{\theta}_k-\theta_k^\star)  f^{(2)}_{jk} \left( x^m , \varphi ;  \ParVT \right) \right|} \\
&\leq& 2 \| \widehat{\ParV}_n(\varphi)- \ParVT \|^2 \times \frac1{|\tilde\Lambda_n|} \ism[\tilde\Lambda_n \times \mSp]{ \| \Vect{\underline{f}}^{(2)}_{} \left( x^m , \varphi ;  \ParVT \right)\|} \\
&\leq& 4 \| \widehat{\ParV}_n(\varphi)- \ParVT \|^2 \times \Esp\left( \| \Vect{\underline{f}}^{(2)}_{} \left( 0^M , \Phi ;  \ParVT \right)\|\right)
\end{eqnarray*}
and
\begin{eqnarray*}
|\tilde\Lambda_n|^{-1} |T_2^\prime(\varphi)| &\leq & \frac2{|\tilde\Lambda_n|} \sum_{x^m \in \varphi_{\tilde\Lambda_n}}  
\sum_{j,k} \left| (\widehat{\theta}_j-\theta_j^\star)(\widehat{\theta}_k-\theta_k^\star)  h^{(2)}_{jk} \left( x^m , \varphi\setminus x^m ;  \ParVT \right) \right| \\
&\leq& 2 \| \widehat{\ParV}_n(\varphi)- \ParVT \|^2 \times \frac1{|\tilde\Lambda_n|}  \sum_{x^m \in \varphi_{\tilde\Lambda_n}} \| \Vect{\underline{h}}^{(2)}_{} \left( x^m , \varphi\setminus x^m ;  \ParVT \right) \| \\
&\leq& 4 \| \widehat{\ParV}_n(\varphi)- \ParVT \|^2 \times \Esp\left( \| \Vect{\underline{h}}^{(2)}_{} \left( 0^M , \Phi ;  \ParVT \right)\|e^{- V^{}\left( 0^M|\Phi; \ParVT \right)}
\right)
\end{eqnarray*}
Since
\begin{eqnarray*}
{}|\tilde\Lambda_n|^{1/2} \| \widehat{\ParV}_n(\varphi)- \ParVT \|^2 = \left(\frac{|\tilde\Lambda_n|}{|\Lambda_n|}\right)^{1/2} \| \; |\Lambda_n|^{1/2} (\widehat{\ParV}_n(\varphi)- \ParVT)\| \times \| \widehat{\ParV}_n(\varphi)- \ParVT\|{}
\end{eqnarray*}
then, under the assumptions of Proposition~\ref{prop-equiv}, one has, from Slustsky's theorem, the following convergence in probability as $n \to +\infty$
$$
{}|\tilde\Lambda_n|^{1/2} \| \widehat{\ParV}_n(\Phi)- \ParVT \|^2 \stackrel{P}{\longrightarrow} 0.
$$
By combining all these results, one obtains the following convergence in probability, as $n\to+\infty$
$$
{}|\tilde\Lambda_n|^{-1/2}\left(T_1(\Phi) - T_2(\Phi)- |\tilde\Lambda_n| \tr{\left( \widehat{\ParV}_n(\Phi)-\ParVT\right)} \Vect{X}_{\tilde\Lambda_n}(\Phi) \right)
= |\tilde\Lambda_n|^{-1/2} \left( T^\prime_1(\Phi) - T^\prime_2(\Phi) \right)
 \stackrel{P}{\longrightarrow} 0.
$$
where $\Vect{X}_{\tilde\Lambda_n}(\Phi)$ is the random vector defined for all $j=1,\ldots,p$ by
$$
\left(\Vect{X}_{\tilde\Lambda_n}(\Phi)\right)_j := \frac{1}{|\tilde\Lambda_n| }\ism[\tilde\Lambda_n\times \mSp]{ {f}^{(1)}_{j} \left( x^m , \Phi ;   \ParVT \right)} -\frac{1}{|\tilde\Lambda_n| } \sum_{x^m \in \Phi_{\tilde\Lambda_n}} {h}^{(1)}_{j} \left( x^m , \Phi\setminus x^m ;   \ParVT \right).
$$
By using the ergodic theorem and the GNZ formula, one has $P_{\ParVT}-$a.s. as $n\to +\infty$
$$
\left(\Vect{X}_{\tilde\Lambda_n}(\Phi)\right)_j \rightarrow \Esp\left(  {f}^{(1)}_{j} \left( 0^M , \Phi ;   \ParVT \right) -  {h}^{(1)}_{j} \left( 0^M , \Phi ;   \ParVT \right) e^{- V^{}\left( 0^M|\Phi; \ParVT \right)}\right).
$$
Finally, let us notice that for all $(m,\varphi) \in \mSp\times \Omega$ and for all $j=1,\ldots,p$ 
\begin{eqnarray*}
{f}^{(1)}_{j} \left( 0^m , \varphi ;   \ParVT \right) &=& \frac{\partial}{\partial \theta_j} \left. \left( h \left( 0^m , \varphi ; \ParV \right) e^{- V^{}\left( 0^m|\varphi; \ParV \right)} \right) \right|_{\ParV=\ParVT}\\
&=&  {h}^{(1)}_{j} \left( 0^m , \varphi ;   \ParVT \right) e^{- V^{}\left( 0^m|\varphi; \ParVT \right)} - h \left( 0^m , \varphi ; \ParV \right)  {V}^{(1)}_{j}\left( 0^m|\varphi; \ParVT \right)e^{- V^{}\left( 0^m|\varphi; \ParVT \right)}.
\end{eqnarray*}
Therefore $\left(\Vect{X}_{\tilde\Lambda_n}(\Phi)\right)_j  \rightarrow -\Vect{\mathcal{E}}_j(h,\ParVT)$ $P_{\ParVT}-$a.s. as $n\to +\infty$. This finally leads to the following convergence in probability, as $n\to+\infty$
$$
|\tilde\Lambda_n|^{-1/2}\left(T_1(\Phi) -T_2(\Phi) +|\tilde\Lambda_n| \tr{\left( \widehat{\ParV}_n(\Phi)-\ParVT\right)} \Vect{\mathcal{E}}(h;\ParVT)) \right) \stackrel{P}{\longrightarrow} 0.
$$


\subsection{Proof of Proposition~\ref{prop-fwk1}}
Let us first state a result widely used in the following.
\begin{lemma} \label{lem-centCond}
For any bounded domain $\Lambda$ and for any test function $h$ 
\begin{equation} \label{eq-centCond}
\Esp \left( I_\Lambda ( \Phi ;h,\ParVT ) | \Phi_{\Lambda^c} \right)  = 0.
\end{equation}
\end{lemma}
The proof of Lemma~\ref{lem-centCond} is omitted since it corresponds to the proof of Theorem~2 (Step~1, p.~257) of \cite{A-BilCoeDro08} by subsituting $v_j(x^m|\varphi)$ by the test function $h(x^m,\varphi;\ParVT)$. \\

For all $n\in\NN$, the domain $\Lambda_n$ is assumed to be a cube divided as $\Lambda_n = \bigcup_{j \in \mathcal J} \Lambda_{j,n}$ where for all $j\in\mathcal J$,  the $\Lambda_{j,n}$'s are disjoint cubes. So $|\Lambda_n|=|\mathcal J||\Lambda_{j,n}|=|\mathcal J||\Lambda_{0,n}|$. Moreover, for all $j\in\mathcal J$, we can decompose each $\Lambda_{j,n}$ in the following way :
 \begin{equation}\label{decomposition}
 \Lambda_{j,n} := \bigcup_{k \in \mathcal K_{j,n}} \Delta_k(D_n)\end{equation}
where the $\Delta_k(D_n)$'s are disjoint cubes with side-length $D_n$ and $\mathcal K_{j,n}\subset\ZZ[d]$.  The side-length $D_n$ is chosen greater than $D$ and as close as possible to $D$, leading to 
$$D_n= \frac{|\Lambda_n|^{1/d}}{|\mathcal J|^{1/d}\left\lfloor \frac{|\Lambda_n|^{1/d}}{|\mathcal J|^{1/d}D}\right\rfloor}.$$
This choice implies $D_n\to D$ when $n\to\infty$ and guarantees $D\leq D_n \leq  2D$ as soon as $|\Lambda_n|\geq|\mathcal J|D^d$.  
The cubes $\Lambda_{j,n}$'s are therefore divided into $|\mathcal K_{j,n}|=|\Lambda_{0,n}|D_n^{-d}$ cubes whose volumes are closed to $D^d$. Denoting $\mathcal K_{n}=\bigcup_{j \in \mathcal J} \mathcal K_{j,n}$, we have $|\mathcal K_n|=|\Lambda_n|D_n^{-d}=|\mathcal J||\mathcal K_{j,n}|$ and finally
\begin{equation}\label{eq-decompLnF1}
\Lambda_n  = \bigcup_{j \in \mathcal J} \bigcup_{k \in \mathcal K_{j,n}} \Delta_k(D_n) = \bigcup_{k \in \mathcal K_n} \Delta_k(D_n).
\end{equation}

From Proposition~\ref{prop-equiv} and under Assumption~\textbf{[E2(bis)]}, one has for any $j\in\mathcal J$
$$
{}|\Lambda_{j,n}|^{-1/2} R_{\Lambda_{j,n}}\left( \Phi;h,\widehat{\ParV}_n(\Phi) \right) = {}|\Lambda_{j,n}|^{-1/2} R_{\infty, \Lambda_{j,n}}\left( \Phi;h,{\ParVT} \right) + o_P(1),
$$
where $R_{\infty, \Lambda_{j,n}}\left( \Phi;h,{\ParVT} \right)$ is defined in (\ref{def-Rinfty}).

Therefore the proof of Proposition~\ref{prop-fwk1} reduces to the proof of the asymptotic normality of the vector $\left(|\Lambda_{j,n}|^{-1/2} R_{\infty,\Lambda_{j,n}}\left( \Phi;h,\ParVT \right)\right)_{j\in\mathcal J}$. Now
\begin{align}\label{somme}
{}|\Lambda_{j,n}|^{-1/2} R_{\infty,\Lambda_{j,n}}\left( \Phi;h,\ParVT \right) &= 
{}|\Lambda_{0,n}|^{-1/2}\left( {I}_{\Lambda_{j,n}}\left( \Phi ; h, \ParVT \right) -
\frac{|\Lambda_{0,n}|}{|\Lambda_n|} \tr{{\Vect{U}}_{\Lambda_n}\left( \Phi ; \ParVT \right) }\Vect{\mathcal{E}}(h;\ParVT)
\right) \nonumber\\
& =  \frac{|\Lambda_{0,n}|^{1/2}}{|\Lambda_n|}  \left(|\mathcal J|\times I_{\Lambda_{j,n}}\left( \Phi ; h, \ParVT \right) - \tr{ {\Vect{U}}_{\Lambda_n}\left( \Phi ; \ParVT \right) }\Vect{\mathcal{E}}(h;\ParVT)\right) \nonumber\\
& =  \frac{1}{D_n^{d/2} |\mathcal J|^{1/2}} \frac{1}{|\mathcal K_n|^{1/2}} \sum_{k \in \mathcal K_n} W_{j,n,\Delta_k(D_n)}\left( \Phi ; \ParVT \right),
\end{align}
where for any $\varphi \in \Omega$
\begin{equation}\label{defW}
W_{j,n,\Delta_k(D_n)}\left( \varphi ; \ParVT \right) = \left\{ \begin{array}{ll}
{}W^{(1)}_{\Delta_k(D_n)}\left( \varphi ; \ParVT \right): =|\mathcal J|\times I_{\Delta_k(D_n)}\left( \varphi ; h, \ParVT \right) &\\
\qquad\qquad \qquad\qquad-\tr{\Vect{U}_{\Delta_k(D_n)}\left( \varphi ; \ParVT \right) }\Vect{\mathcal{E}}(h;\ParVT) & \mbox{ if } k\in \mathcal K_{j,n},\\
W^{(2)}_{\Delta_k(D_n)}\left( \varphi ; \ParVT \right):= - \tr{{\Vect{U}}_{\Delta_k(D_n)}\left( \varphi ; \ParVT \right) }\Vect{\mathcal{E}}(h;\ParVT) & \mbox{ if } k\in \mathcal K_n\setminus \mathcal K_{j,n}.
\end{array} \right.
\end{equation}

Therefore, to prove a central limit theorem for the vector  $\left(|\Lambda_{j,n}|^{-1/2} R_{\infty,\Lambda_{j,n}}\left( \Phi;h,\ParVT \right)\right)_{j\in\mathcal J}$, it suffices to apply Theorem \ref{tcl} (see Appendix \ref{annexe-tcl}), where in its statement we choose $\Vect Z_{n,k}=(W_{j,n,\Delta_k(D_n)}(\Phi;\ParVT))_{j\in\mathcal J}$, $X_{n,i}=\Phi_{\Delta_i(D_n)}$ and $p=|\mathcal J|$. For this, we first have to specify the asymptotic variance matrix $\Mat \Sigma$, then to check the assumptions of Theorem \ref{tcl}.\\

{\it \underline{First step}: computation of the asymptotic variance.}

Let us fix a cartesian coordinate system such that $0$ is the center of $\Lambda_n$. 
We assume, without loss of generality, that $|\mathcal J|$ is odd. Moreover, we can always choose an odd number $|\mathcal K_{j,n}|$ of cubes $\Delta_k(D_n)$ in (\ref{decomposition}). Consequently, $\Lambda_{0,n}$ may be centered at $0$ and each $\Delta_k(D_n)$ is centered at $kD_n$, $k\in\ZZ[d]$. Note that if $|\mathcal J|$ was even, each $\Delta_k(D_n)$ would be centered at $kD_n/2$.
So, in this system, $\mathcal K_n$ is a subset of $\ZZ[d]$, independent of $D_n$, with $|\mathcal K_n|=|\mathcal J|\left\lfloor \frac{|\Lambda_n|^{1/d}}{|\mathcal J|^{1/d}D}\right\rfloor^d$ elements.

Set, for all $k,k' \in \ZZ[d]$,
\begin{eqnarray*}
\left\{ \begin{array}{ll}
E_{k,k^\prime}^{(1)} (D_n):= \Esp\left(
 W^{(1)}_{\Delta_k(D_n)}\left( \Phi ; \ParVT \right)W^{(1)}_{\Delta_{k^\prime}(D_n)}\left( \Phi ; \ParVT \right)
\right) \\
 E_{k,k^\prime}^{(12)}(D_n) :=\Esp\left(  
W^{(1)}_{\Delta_k(D_n)}\left( \Phi ; \ParVT \right)W^{(2)}_{\Delta_{k^\prime}(D_n)}\left( \Phi ; \ParVT \right)
\right) \\
E_{k,k^\prime}^{(2)} (D_n):= \Esp\left( 
W^{(2)}_{\Delta_k(D_n)}\left( \Phi ; \ParVT \right)W^{(2)}_{\Delta_{k^\prime}(D_n)}\left( \Phi ; \ParVT \right)
 \right)  
\end{array}\right. 
\end{eqnarray*}
Note that from the stationarity of the point process, we have $E_{k,k^\prime}^{(l)}(D_n)=E_{0,k-k^\prime}^{(l)}(D_n)$, for $l=1,12,2$. 
Moreover, under Assumptions \textbf{[N4]} and  \textbf{[E2(bis)]}, for any $k\in \mathcal K_n$ and for any configuration $\varphi$, since $D_n\geq D$, $W^{(i)}_{\Delta_k(D_n)}\left( \varphi ; \ParVT \right)$, $i=1,2$,  depends only on $\varphi_{\Delta_l(D_n)}$ for $|l-k|\leq 1$ that is $l\in\mathbbm{B}_k(1)$. As a consequence, if  $k^\prime\in\mathbbm{B}^c_k(1)$, $W^{(i)}_{\Delta_{k^\prime}(D_n)}\left( \Phi ; \ParVT \right)$ is a measurable function of $\Phi_{\Delta_k^c(D_n)}$. This  leads, for $i,j=1,2$, to
\begin{align}\label{markov}
\Esp \left( W^{(i)}_{\Delta_k(D_n)}\left( \Phi ; \ParVT \right)
W^{(j)}_{\Delta_{k^\prime}(D_n)}\left( \Phi ; \ParVT \right)
\right) 
&= \Esp\left( \Esp\left( 
W^{(i)}_{\Delta_k(D_n)}\left( \Phi ; \ParVT \right)
W^{(j)}_{\Delta_{k^\prime}(D_n)}\left( \Phi ; \ParVT \right) | \Phi_{\Delta_k^c(D_n)}
\right) \right)\nonumber\\
&= \Esp\left( W^{(i)}_{\Delta_{k^\prime}(D_n)}\left( \Phi ; \ParVT \right)
\Esp\left( 
W^{(j)}_{\Delta_k(D_n)}\left( \Phi ; \ParVT \right)
 {}| \Phi_{\Delta_k^c(D_n)}
\right) \right).
\end{align}
From Lemma~\ref{lem-centCond} and under \textbf{[E2(bis)]} then for any $k\in \ZZ[d]$ and for $i=1,2$, 
\begin{equation}\label{eq-propMarkovW}
\Esp\left( 
W^{(i)}_{\Delta_k(D_n)}\left( \Phi ; \ParVT \right) | \Phi_{\Delta_k^c(D_n)}
\right) =0.
\end{equation}
From (\ref{markov}) and (\ref{eq-propMarkovW}), we deduce that, for $l=1,12,2$, 
\begin{equation}\label{portee}
k^\prime\in\mathbbm{B}^c_k(1) \Longrightarrow E^{(l)}_{k,k'}(D_n)=0.
\end{equation}

We are now in position to compute the covariance. For any $i$ and $j$ in $\mathcal J$, from (\ref{somme}),
\begin{multline}\label{cov}
cov\left(|\Lambda_{i,n}|^{-1/2} R_{\infty,\Lambda_{i,n}}\left( \Phi;h,\ParVT \right),|\Lambda_{j,n}|^{-1/2} R_{\infty,\Lambda_{j,n}}\left( \Phi;h,\ParVT \right)\right)\\= \frac{1}{D_n^d |\mathcal J|} \Esp\left(\frac{1}{|\mathcal K_n|} \sum_{k \in \mathcal K_n} \sum_{k' \in \mathcal K_n} W_{i,n,\Delta_k(D_n)}\left( \Phi ; \ParVT \right)W_{j,n,\Delta_{k'}(D_n)}\left( \Phi ; \ParVT \right)\right).
\end{multline}

Let us first consider the case $i=j$. We may write 
\begin{multline*}
\Esp\left(\frac{1}{|\mathcal K_n|} \sum_{k \in \mathcal K_n} \sum_{k' \in \mathcal K_n} W_{i,n,\Delta_k(D_n)}\left( \Phi ; \ParVT \right)^2\right)\\ =
\frac1{|\mathcal K_n|}\bigg( \underbrace{\sum_{k,k^\prime\in \mathcal K_{i,n}} E_{k,k^\prime}^{(1)}(D_n)}_{:=S_1}+2 \underbrace{\sum_{k\in \mathcal K_{i,n},k^\prime\in \mathcal K_n\setminus \mathcal K_{i,n}} E_{k,k^\prime}^{(12)}(D_n)}_{:=S_2} + \underbrace{\sum_{k,k^\prime\in \mathcal K_n\setminus \mathcal K_{i,n}} E_{k,k^\prime}^{(2)}(D_n)}_{:=S_3}\bigg).
\end{multline*}

The following lemma will be useful to drop the dependence on $D_n$ in each term $S_1$, $S_2$, $S_3$ above.
\begin{lemma}\label{dropDn}
For any $i,j=1,2$, denoting $\overline\Delta_{0}(\tau)=\cup_{k \in \mathbbm B_0(1)}\Delta_{k}(\tau)$ (for some $\tau>0$), we have
$$ W^{(i)}_{\Delta_0(D_n)}\left( \Phi ; \ParVT \right)W^{(j)}_{\overline\Delta_{0}(D_n)}\left( \Phi ; \ParVT \right)\overset{L_1}{\longrightarrow} W^{(i)}_{\Delta_0(D)}\left( \Phi ; \ParVT \right)W^{(j)}_{\overline\Delta_{0}(D)}\left( \Phi ; \ParVT \right).$$
\end{lemma}

\begin{proof}
For any $i=1,2$,  $W^{(i)}_{\Delta_0(D_n)}$ is a linear combination of $I_{\Delta_0(D_n)}$ and  $\Vect U_{\Delta_0(D_n)}$, which converge respectively in $L^2$ to $I_{\Delta_0(D)}$ and  $\Vect U_{\Delta_0(D)}$ by \textbf{[N3]}  and \textbf{[E2(bis)]}, since $D_n\to D$. Thus  $W^{(i)}_{\Delta_0(D_n)}$ converges in $L^2$ to  $W^{(i)}_{\Delta_0(D)}$ as  $n\to\infty$. Similarly, for any $j=1,2$, $W^{(j)}_{\overline\Delta_{0}(D_n)}$ tends in $L^2$ to $W^{(j)}_{\overline\Delta_{0}(D)}$. The convergence stated in Lemma \ref{dropDn} then follows.
\end{proof}

Let us focus on the asymptotic of each term $S_1$, $S_2$, $S_3$.

\underline{Term $S_1$}: from (\ref{portee}),
$$
S_1 = \sum_{k\in \mathcal K_{i,n}} \bigg( 
\sum_{k^\prime \in \mathbbm{B}_k(1)\cap \mathcal K_{i,n}} E_{k,k^\prime}^{(1)}(D_n) + \underbrace{\sum_{k^\prime \in \mathbbm B_k^c(1)\cap \mathcal K_{i,n}} E_{k,k^\prime}^{(1)}(D_n)}_{=0} 
\bigg) = \sum_{k\in \mathcal K_{i,n}} \sum_{k^\prime \in \mathbbm B_k(1)\cap \mathcal K_{i,n}} E_{k,k^\prime}^{(1)}(D_n).
$$
Let $\widetilde{\mathcal K}_{i,n}:= \mathcal K_{i,n} \cap \left( \cup_{j\in \partial \mathcal K_{i,n}} \mathbbm B_j(1)\right)$ and note that $\frac{|\widetilde{\mathcal K}_{i,n}|}{|\mathcal K_{i,n}|} \to 0$ as $n\to +\infty$. Then,
$$
S_1 = \sum_{k\in \mathcal K_{i,n} \setminus \widetilde{\mathcal K}_{i,n}} \sum_{k^\prime \in \mathbbm B_k(1)\cap \mathcal K_{i,n}} E_{k,k^\prime}^{(1)}(D_n) +\underbrace{\sum_{k\in \widetilde{\mathcal K}_{i,n}} \sum_{k^\prime \in \mathbbm B_k(1)\cap \mathcal K_{i,n}} E_{k,k^\prime}^{(1)}(D_n) }_{:=A_1}.
$$
Since, 
$$
{}\frac{1}{|\mathcal K_n|}\times |A_1|\leq \frac{|\widetilde{\mathcal K}_{i,n}|}{|\mathcal K_n|} \sum_{k \in \mathbbm B_0(1)} |E_{0,k}^{(1)}(D_n)|{} \stackrel{n\to +\infty}{\longrightarrow} 0,
$$
(because $D\leq D_n\leq 2D$ and $\frac{|\widetilde{\mathcal K}_{i,n}|}{|\mathcal K_{i,n}|} \to 0$),  we obtain,  as $n\to +\infty$,
\begin{eqnarray*}
\frac1{|\mathcal K_n|} \; S_1 &\sim& \frac{|\mathcal K_{i,n} \setminus \widetilde{\mathcal K}_{i,n}|}{|\mathcal K_n|} \; \sum_{k \in \mathbbm B_0(1)} E_{0,k}^{(1)}(D_n){} 
\sim \frac{|\mathcal K_{i,n} |}{|\mathcal K_n|} \; \sum_{k \in \mathbbm B_0(1)} E_{0,k}^{(1)}(D_n). 
\end{eqnarray*}

From Lemma \ref{dropDn}, $$\sum_{k \in \mathbbm B_0(1)} E_{0,k}^{(1)}(D_n)=\Esp\left(W^{(1)}_{\Delta_0(D_n)}\left( \Phi ; \ParVT \right)W^{(1)}_{\overline\Delta_{0}(D_n)}\left( \Phi ; \ParVT \right)\right)\longrightarrow \sum_{k \in \mathbbm B_0(1)} E_{0,k}^{(1)}(D).$$
Therefore,
$$\frac1{|\mathcal K_n|} \; S_1 \sim\frac1{|\mathcal J|} \; \sum_{k \in \mathbbm B_0(1)} E_{0,k}^{(1)}(D){} .$$ 

\underline{Term $S_2$}: with similar arguments as above, we obtain
\begin{eqnarray*}
S_2 &=& \underbrace{\sum_{k \in \mathcal K_{i,n}\setminus \widetilde{\mathcal K}_{i,n}} \sum_{k^\prime \in \mathcal K_n \setminus \mathcal K_{i,n}} E_{k,k^\prime}^{(12)}(D_n)}_{=0} \\ &&\hspace*{.5cm}  +
\sum_{k \in \widetilde{\mathcal K}_{i,n}} \bigg(
\sum_{k^\prime \in \mathbbm B_k(1)\cap (\mathcal K_n \setminus \mathcal K_{i,n})} E_{k,k^\prime}^{(12)}(D_n) +
\underbrace{\sum_{k^\prime \in \mathbbm B_k^c(1)\cap (\mathcal K_n \setminus \mathcal K_{i,n})} E_{k,k^\prime}^{(12)}(D_n)}_{=0} 
\bigg) 
\end{eqnarray*}
Therefore, since $\frac{|\widetilde{\mathcal K}_{i,n}|}{|\mathcal K_{n}|} \to 0$ and $D\leq D_n\leq 2D$,
$$\frac1{|\mathcal K_n|} \; S_2 \leq \frac{|\widetilde{\mathcal K}_{i,n} |}{|\mathcal K_n|}\sum_{k \in \mathbbm B_0(1)} |E_{0,k}^{(12)}(D_n)|\stackrel{n\to +\infty}{\longrightarrow} 0.
$$


\underline{Term $S_3$}: 
$$
S_3 = 
\sum_{k \in \mathcal K_n \setminus \mathcal K_{i,n}} \sum_{k^\prime \in \mathbbm B_k(1)\cap (\mathcal K_n \setminus \mathcal K_{i,n})} E_{k,k^\prime}^{(2)}(D_n)
+\underbrace{\sum_{k \in \mathcal K_n \setminus \mathcal K_{i,n}} \sum_{k^\prime \in \mathbbm B_k^c(1)\cap (\mathcal K_n \setminus \mathcal K_{i,n})} E_{k,k^\prime}^{(2)}(D_n)}_{=0}.
$$
Let $\widetilde{\mathcal K}_n = (\mathcal K_n\setminus \mathcal K_{i,n}) \cap \left( \cup_{j \in \partial(\mathcal K_n\setminus \mathcal K_{i,n})} \mathbbm B_j(1) \right)$ and note that $\frac{|\widetilde{\mathcal K}_n|}{|\mathcal K_n|} \to 0$, as $n\to +\infty$. Then,
$$
S_3 = \sum_{k \in \mathcal K_n\setminus \widetilde{\mathcal K}_n} \sum_{k^\prime \in \mathbbm B_k(1)\cap(\mathcal K_n \setminus \mathcal K_{i,n})} E_{k,k^\prime}^{(2)}(D_n)+
\underbrace{\sum_{k \in \widetilde{\mathcal K}_n} \sum_{k^\prime \in \mathbbm B_k(1)\cap(\mathcal K_n \setminus \mathcal K_{i,n})}E_{k,k^\prime}^{(2)}(D_n)}_{:=A_3}.
$$
Since,
$$
\frac1{|\mathcal K_n|} |A_3| \leq  \frac{| \widetilde{\mathcal K}_n|}{|\mathcal K_n|} \sum_{k \in \mathbbm B_0(1)} |E_{0,k}^{(2)}(D_n)| \stackrel{n\to +\infty}{\longrightarrow} 0,
$$
we obtain, from Lemma \ref{dropDn},
$$
\frac{1}{|\mathcal K_n|} S_3 \sim \frac{|\mathcal K_n\setminus \widetilde{\mathcal K}_n|}{|\mathcal K_n|} \sum_{k \in \mathbbm B_0(1)} E_{0,k}^{(2)}(D_n) \sim
\frac{|\mathcal J|-1}{|\mathcal J|}\sum_{k \in \mathbbm B_0(1)} E_{0,k}^{(2)}(D).
$$
Combining the three terms $S_1$, $S_2$ and $S_3$, we have, as $n\to +\infty$
\begin{equation}\label{i=j}
\Esp\left(\frac{1}{|\mathcal K_n|} \sum_{k \in \mathcal K_n} \sum_{k' \in \mathcal K_n} W_{i,n,\Delta_k(D_n)}\left( \Phi ; \ParVT \right)^2\right) \sim \sum_{k \in \mathbbm B_0(1)} \left( \frac1{|\mathcal J|} E_{0,k}^{(1)}(D) + \frac{|\mathcal J|-1}{|\mathcal J|} E_{0,k}^{(2)}(D)\right).
\end{equation}

When $i\not=j$, there are three main cases in (\ref{cov}), according to $k,k' \in \mathcal K_{i,n}$, $k,k'\in \mathcal K_{j,n}$, or $k,k' \in \mathcal K_n\setminus(\mathcal K_{i,n}\cup \mathcal K_{j,n})$.  As for the case $i=j$ treated before, the other situations involve non-zero correlations on edges sets like $\widetilde{\mathcal K}_{i,n}$, which are negligible with respect to $|\mathcal K_n|$.
The covariance is therefore equivalent, up to $D_n^d|\mathcal J|$, to

$$\frac{1}{|\mathcal K_n|} \sum_{k,k' \in \mathcal K_{i,n}}E_{k,k^\prime}^{(12)}(D_n) + \frac{1}{|\mathcal K_n|} \sum_{k,k' \in \mathcal K_{j,n}} E_{k,k^\prime}^{(12)}(D_n)
+ \frac{1}{|\mathcal K_n|} \sum_{k,k' \in \mathcal K_n\setminus(\mathcal K_{i,n}\cup \mathcal K_{j,n})} E_{k,k^\prime}^{(2)}(D_n).$$
The simplification occurs as for the case $i=j$ and, since $|\mathcal K_{i,n}|=|\mathcal K_{j,n}|$, we obtain the asymptotic equivalent for the covariance  (\ref{cov})
\begin{equation}\label{inot=j}
\frac{1}{D^d |\mathcal J|} \sum_{k \in \mathbbm B_0(1)} \left( \frac 2{|\mathcal J|} E_{0,k}^{(12)}(D) + \frac{|\mathcal J|-2}{|\mathcal J|} E_{0,k}^{(2)} (D)\right).
\end{equation}

Finally, from (\ref{i=j}) and (\ref{inot=j}), we deduce that $\Mat \Sigma_1(\ParVT)$, defined in Proposition \ref{prop-fwk1}, corresponds to the asymptotic variance of $\left(|\Lambda_{j,n}|^{-1/2} R_{\infty,\Lambda_{j,n}}\left( \Phi;h,\ParVT \right)\right)_{j\in\mathcal J}$.\\

{\it \underline{Second step}: application of Theorem \ref{tcl}.}

We apply Theorem \ref{tcl} with $\Vect Z_{n,k}=(W_{j,n,\Delta_k(D_n)})_{j\in\mathcal J}$, $X_{n,i}=\Phi_{\Delta_i(D_n)}$, $p=|\mathcal J|$ and $\Mat \Sigma=\Mat \Sigma_1(\ParVT)$, which is a symmetric positive-semidefinite matrix as the limit of a covariance matrix (from the first step of the proof).

The assumption  (\ref{subordinate}) holds from \textbf{[N4]}, \textbf{[E2(bis)]} and because $D_n\geq D$. Assumptions $(i)$, $(ii)$ and $(iii)$ are  direct consequences of  \textbf{[E2(bis)]}, \textbf{[N2]} and Lemma~\ref{lem-centCond}. It remains to prove $(iv)$. Assuming $\Mat\Sigma=(\Sigma_{ij})$ for $1\leq i,j\leq p$, from the definition of the Frobenius norm, we have

\begin{multline}\label{CVvariance}
\Esp \left\| |\mathcal K_n|^{-1} \sum_{k\in \mathcal K_n}\sum_{k'\in \mathbbm B_k(1)\cap \mathcal K_n} \Vect Z_{n,k}\tr{\Vect Z_{n,k'}} - \Mat \Sigma \right\|\\ \leq \sum_{i=1}^p\sum_{j=1}^p \Esp\left| |\mathcal K_n|^{-1} \sum_{k\in \mathcal K_n}\sum_{k'\in \mathbbm B_k(1)\cap \mathcal K_n} W_{i,n,\Delta_k(D_n)}W_{j,n,\Delta_{k'}(D_n)}-\Sigma_{ij}\right|.\end{multline}

Let us first assume that $i\not=j$ are fixed and denote  $Y_{n,k}(D_n)= W_{i,n,\Delta_k(D_n)}$, $S_n^k(D_n)=\sum_{k'\in \mathbbm B_k(1)\cap \mathcal K_n} W_{j,n,\Delta_{k'}(D_n)}$. We have
\begin{align*}
\Esp \left| |\mathcal K_n|^{-1} \sum_{k\in \mathcal K_n}\sum_{k'\in \mathbbm B_k(1)\cap \mathcal K_n} W_{i,n,\Delta_k(D_n)}W_{j,n,\Delta_{k'}(D_n)}-{\Sigma}_{ij}\right| &=
|\mathcal K_n|^{-1} \Esp \left| \sum_{k\in \mathcal K_n} Y_{n,k}(D_n) S_n^k(D_n)  -{\Sigma}_{ij}  \right|
\\
&\leq E_1+E_2+E_3+E_4,
\end{align*}
where
$$E_1=\frac{|\mathcal K_{i,n}|}{|\mathcal K_n|}\Esp\left| |\mathcal K_{i,n}|^{-1}\sum_{k\in \mathcal K_{i,n}} (Y_{n,k}(D_n)S_n^k(D_n)-\Esp(Y_{n,k}(D_n)S_n^k(D_n)))\right|,$$
$$E_2=\frac{|\mathcal K_{j,n}|}{|\mathcal K_n|}\Esp\left| |\mathcal K_{j,n}|^{-1}\sum_{k\in \mathcal K_{j,n}} (Y_{n,k}(D_n)S_n^k(D_n)-\Esp(Y_{n,k}(D_n)S_n^k(D_n)))\right|,$$
\begin{eqnarray*}
E_3&=&\frac{|\mathcal K_n\setminus(\mathcal K_{i,n}\cup \mathcal K_{j,n})|}{|\mathcal K_n|}\times \\
&&\qquad \Esp\left| |\mathcal K_n\setminus(\mathcal K_{i,n}\cup \mathcal K_{j,n})|^{-1} 
 \!\! \sum_{k\in \mathcal K_n\setminus(\mathcal K_{i,n}\cup \mathcal K_{j,n})}\!\!\!\!\!\!\! (Y_{n,k}(D_n)S_n^k(D_n)-\Esp(Y_{n,k}(D_n)S_n^k(D_n)))\right|,
\end{eqnarray*}
\begin{multline*}
E_4=\Bigg| \frac{|\mathcal K_{i,n}|}{|\mathcal K_n|}\sum_{k\in \mathcal K_{i,n}}\Esp(Y_{n,k}(D_n)S_n^k(D_n))+\frac{|\mathcal K_{i,n}|}{|\mathcal K_n|}\sum_{k\in \mathcal K_{j,n}}\Esp(Y_{n,k}(D_n)S_n^k(D_n))\\+\frac{|\mathcal K_n\setminus(\mathcal K_{i,n}\cup \mathcal K_{j,n})|}{|\mathcal K_n|}\sum_{k\in \mathcal K_n\setminus(\mathcal K_{i,n}\cup \mathcal K_{j,n})}\Esp(Y_{n,k}(D_n)S_n^k(D_n))-{\Sigma}_{ij} \Bigg|.\end{multline*}

The first three terms $E_1$, $E_2$ and $E_3$ can be handled similarly. Let us focus on $E_1$:
\begin{align}\label{E1}
\frac{|\mathcal K_n|}{|\mathcal K_{i,n}|}E_1&\leq  |\mathcal K_{i,n}|^{-1} \sum_{k\in \mathcal K_{i,n}} \Esp\left|Y_{n,k}(D_n)S_n^k(D_n)-Y_{n,k}(D)S_n^k(D)\right|\nonumber\\
&\quad+ |\mathcal K_{i,n}|^{-1} \Esp\left|\sum_{k\in \mathcal K_{i,n}} (Y_{n,k}(D)S_n^k(D)-\Esp(Y_{n,k}(D)S_n^k(D)))\right|\nonumber\\
&\quad+ |\mathcal K_{i,n}|^{-1} \sum_{k\in \mathcal K_{i,n}}\left|\Esp\left(Y_{n,k}(D)S_n^k(D)\right)-\Esp\left(Y_{n,k}(D_n)S_n^k(D_n)\right)\right|.
\end{align} 
Up to the edge effects which are negligible with respect to $|\mathcal K_{i,n}|$, $\left(Y_{n,k}(D)S_n^k(D)\right)_k$ is stationary when $k\in\mathcal K_{i,n}$, since in this case, from (\ref{defW}), $W_{i,n,\Delta_k(D)}=W^{(1)}_{\Delta_k(D)}$ does not depend on $n$. Therefore the second term in (\ref{E1}) tends to 0 by the mean ergodic theorem. For a fixed $n$, we have also by stationarity (up to the edge effects)
\begin{align*}
\Esp\left|Y_{n,k}(D_n)S_n^k(D_n)-Y_{n,k}(D)S_n^k(D)\right|&=\Esp\left|Y_{n,0}(D_n)S_n^0(D_n)-Y_{n,0}(D)S_n^0(D)\right|\\
&=\Esp\left|W^{(1)}_{\Delta_0(D_n)}W^{(2)}_{\overline\Delta_0(D_n)}-W^{(1)}_{\Delta_0(D)}W^{(2)}_{\overline\Delta_0(D)}\right|,
\end{align*} 
where $\overline\Delta_0(D_n)=\cup_{k'\in \mathbbm B_0(1)}\Delta_{k'}(D_n)$. From Lemma \ref{dropDn}, this term tends to 0, therefore the first term in (\ref{E1}) asymptotically vanishes. The same argument shows that the third term in  (\ref{E1}) also tends to 0 as $n\to\infty$. 
As a consequence, $E_1$ tends to 0.

The same decomposition as in (\ref{E1}) may be done for $E_2$ and $E_3$, which leads by similar arguments to $E_2\to0$ and $E_3\to0$.
The last term $E_4$ involves the difference between $\Sigma_{ij}$ and its empirical counterpart. The same calculations as in the first step of the proof shows that $E_4\to 0$.

Therefore, we have proved that the terms in the double-sum (\ref{CVvariance}) corresponding to $i\not=j$ asymptotically vanish. The same result can be proved similarly when $i=j$. Thus assumption $(iv)$ in Theorem \ref{tcl} holds and the convergence in law is deduced.


\subsection{Proof of Proposition~\ref{prop-fwk2}}

We can decompose $\Lambda_{n}$ in the following way :
 $$\Lambda_{n} := \bigcup_{k \in \mathcal K_{n}} \Delta_k(D_n)$$
where the $\Delta_k$'s are disjoint cubes with side-length $D_n$ and $\mathcal K_{n}\subset\ZZ[d]$  satisfies $|\mathcal K_{n}|=|\Lambda_{n}|D_n^{-d}$. Similarly as in the proof of Proposition \ref{prop-fwk1}, we choose 
$$D_n= \frac{|\Lambda_n|^{1/d}}{\left\lfloor \frac{|\Lambda_n|^{1/d}}{D}\right\rfloor},$$
which implies $D_n\to D$ when $n\to\infty$ and guarantees $D\leq D_n \leq  2D$ as soon as $|\Lambda_n|\geq D^d$.  

From Proposition~\ref{prop-equiv} and under Assumption~\textbf{[E2(bis)]}, for all $i=1,\dots,s$,
\begin{eqnarray*}
{}|\Lambda_{n}|^{-1/2} R_{\Lambda_n}\left( \Phi;h_i,\widehat{\ParV}_n(\Phi) \right) & = & {}|\Lambda_{n}|^{-1/2} R_{\infty,\Lambda_n}\left( \Phi;h_i,\ParVT \right) +o_P(1) \\
&=& 
{}|\Lambda_n|^{-1/2}\left( {I}_{\Lambda_n}\left( \Phi ; h_i, \ParVT \right) -
 \tr{{\Vect{U}}_{\Lambda_n}\left( \Phi ; \ParVT \right) }\Vect{\mathcal{E}}(h_i;\ParVT)
\right) +  o_{P}(1)\\
& = & \frac{1}{D_n^{d/2}}\frac{1}{|\mathcal K_n|^{1/2}} \sum_{k \in \mathcal K_n} W_{\Delta_k(D_n)}\left( \Phi ; h_i, \ParVT \right)+  o_{P}(1),
\end{eqnarray*}
where for any $\varphi \in \Omega$
$$
W_{\Delta_k(D_n)}\left( \varphi ; h_i, \ParVT \right) := I_{\Delta_k(D_n)}\left( \varphi ; h_i, \ParVT \right) + \tr{\Vect{U}_{\Delta_k(D_n)}\left( \varphi ; \ParVT \right) }\Vect{\mathcal{E}}(h_i;\ParVT).$$

We apply Theorem \ref{tcl} in the simpler case when $f_{n,k}=f$ for all $n\in\NN$ and all $k\in\mathcal K_n$. If $D_n=D$ for all $n$, this framework would reduce to a stationary setting similar to Theorem 2.1 in \cite{A-JenKun94}. But as $\Lambda_n$ is allowed to increase continuously up to $\RR[d]$, $D_n\equiv D$ is impossible. We will therefore apply Theorem \ref{tcl} in  Appendix \ref{annexe-tcl} with $\Vect Z_{n,k}=\left(W_{\Delta_k(D_n)}\left( \Phi ; h_j,  \ParVT \right)\right)_{j=1\dots s}$, $X_{n,i}=\Phi_{\Delta_i(D_n)}$ and  $p=s$.

Let us first compute the covariance matrix of $\left(|\Lambda_{n}|^{-1/2} R_{\infty,\Lambda_n}\left( \Phi;h_i,\ParVT \right)\right)_{i=1,\dots,s}$. By the same calculations as for the term $S_1$ in the proof of Proposition~\ref{prop-fwk1}, we obtain
\begin{align}\label{sigma2}
cov \bigg(|\Lambda_n|^{-1/2} R_{\infty,\Lambda_n}\left( \Phi; \right. & \left(h_i, \ParVT \right)  ,|\Lambda_n|^{-1/2} R_{\infty,\Lambda_n}\left( \Phi;h_j, \ParVT \right)\bigg) \nonumber\\ 
&= \frac{1}{D_n^d} \Esp\left(\frac{1}{|\mathcal K_n|} \sum_{k \in \mathcal K_n} \sum_{k' \in \mathcal K_n} W_{\Delta_k(D_n)}\left( \Phi ; h_i,  \ParVT \right)W_{\Delta_{k'}(D_n)}\left( \Phi ; h_j, \ParVT \right)\right)  \nonumber\\
&\sim \frac{1}{D^d}  \sum_{k \in \mathbbm B_0(1)} \Esp\left(W_{\Delta_0(D)}\left( \Phi ; h_i,  \ParVT \right)W_{\Delta_{k}(D)}\left( \Phi ; h_j, \ParVT \right)\right).
\end{align}

The asymptotic covariance matrix is thus $\Mat \Sigma_2(\ParVT)$ defined in Proposition~\ref{prop-fwk2}. We can now apply  Theorem~\ref{tcl} in the appendix with $\Mat \Sigma=\Mat \Sigma_2(\ParVT)$.  The assumption (\ref{subordinate}) hods because $D_n\geq D$ and from  \textbf{[N4]} and \textbf{[E2(bis)]}. The assumptions $(i)$, $(ii)$ and $(iii)$ follow from \textbf{[E2(bis)]}, \textbf{[N2]} and Lemma~\ref{lem-centCond}. Assumption $(iv)$ may be checked easily as in the second step of the proof of Proposition~\ref{prop-fwk1}, by using (\ref{sigma2}).

\subsection{Proof of Lemma~\ref{delta}}\label{proofdelta}

For simplicity, let $\Vect{Y}_\Lambda:=\Vect{Y}_\Lambda(\Phi;\ParV)$. Let us denote $\overline{\Delta}(\delta,D^\vee):= \cup_{|j|\leq \lceil D^\vee/\delta\rceil} \Delta_j(\delta)$. From the additivity property of $\Vect Y$, proving Lemma \ref{delta} reduces to prove that for any $\delta>0$ and any $D^\vee\geq D$, $D^d A(\delta,D^\vee)=\delta^d A(D,D)$ where
$$
A(\delta,D^\vee):=\Esp\left(\Vect Y_{\Delta_0(\delta)}\tr{\Vect Y_{\overline{\Delta}(\delta,D^\vee)}}\right).
$$

Since $D^\vee\geq D$, we can write $\overline{\Delta}(\delta,D^\vee)=\overline{\Delta}(\delta,D)\cup \Delta'$, where $\Delta'\subset \left(\overline{\Delta}(\delta,D)\right)^c$.  From the locality assumption, $\Vect Y_{\Delta'}$ is only a function of  $\Phi_{\Delta_0^c(\delta)}$. So 
\begin{equation}\label{decorrelation}\Esp\left(\Vect Y_{\Delta_0(\delta)}\tr{\Vect Y_{\Delta'}}\right)=\Esp\left(\Esp\left(\Vect Y_{\Delta_0(\delta)}\tr{\Vect Y_{\Delta'}}|\Phi_{\Delta_0^c(\delta)}\right)\right)=\Esp\left(\Esp\left(\Vect Y_{\Delta_0(\delta)}|\Phi_{\Delta_0^c(\delta)}\right)\tr{\Vect Y_{\Delta'}}\right)=0,\end{equation}
which yields $A(\delta,D^\vee)=A(\delta,D)$. By denoting $A(\delta):=A(\delta,D^\vee)=A(\delta,D)$ and $\overline{\Delta}(\delta):=\overline{\Delta}(\delta,D)$, we must prove $D^d A(\delta)=\delta^d A(D)$.

Let us first assume $\delta=k D$ with $k\in\NN$. We may write $\overline{\Delta}(kD)=\left(\Delta_0(kD)\oplus D\right) \cup \Delta'$ and may assert that $\Vect Y_{\Delta'}$ depends only on a function of  $\Phi_{\Delta_0^c(kD)}$. By a similar argument as in (\ref{decorrelation}), we obtain $A(\delta)= \Esp ( \Vect Y_{\Delta_0(kD)} \tr{\Vect Y_{\Delta_0(kD)\oplus D}})$. From the disjoint decomposition $\Delta_0(kD)=\cup_{j\in\mathcal K} \Delta_j(D)$ where $|\mathcal K|=k^d$, we have, by the same decorrelation argument as above and by stationarity,
\begin{eqnarray*}
A(\delta) &=& \sum_{j\in\mathcal K} \Esp(\Vect Y_{\Delta_j(D)}\tr{\Vect Y_{\Delta_0(kD)\oplus D}})=\sum_{j\in\mathcal K} \Esp(\Vect Y_{\Delta_j(D)}\tr{\Vect Y_{\Delta_j(D)\oplus D}})\\
&=&k^d \Esp(\Vect Y_{\Delta_0(D)}\tr{\Vect Y_{\Delta_0(D)\oplus D}})=\frac{\delta^d}{D^d}A(D).
\end{eqnarray*}

Let us now assume $D=k\delta$ with $k\in\NN$. First notice that in this case $\overline{\Delta}(D)= \overline{\Delta}(\delta)\oplus \frac{D}{2}(1-1/k)$.  The following decomposition holds: $\Delta_0(D)=\cup_{j\in\mathcal K} \Delta_j(\delta)$ where $|\mathcal K|=k^d$. For any $j\in\mathcal K$, $|j|\leq \frac{D}{2}(1-1/k)$, so $\overline{\Delta}(D)$ contains any translation of the set $\overline{\Delta}(\delta)$ with respect to $j$. Let us denote this translated set by $\tau_j \overline{\Delta}(\delta)$. From the same decorrelation argument as above and by stationarity, we have
\begin{equation}\label{cas2}
A(D)=\sum_{j\in\mathcal K} \Esp(\Vect Y_{\Delta_j(\delta)}\tr{\Vect Y_{\overline{\Delta}(D)}})=\sum_{j\in\mathcal K} \Esp(\Vect Y_{\Delta_j(\delta)}\tr{\Vect Y_{\tau_j \overline{\Delta}(\delta)}})=\sum_{j\in\mathcal K} \Esp(\Vect Y_{\Delta_0(\delta)}\tr{\Vect Y_{\overline{\Delta}(\delta)}})=\frac{D^d}{\delta^d}A(\delta).\end{equation}

Let us now consider the case $D/\delta=k^\prime/k$, where $(k,k^\prime)\in\NN[2]$. Let $\delta^\prime=\delta/k$,  then $D=k^\prime\delta^\prime$ and according to (\ref{cas2}), ${\delta^\prime}^d A(D) = D^d A(\delta^\prime)$. In the same way as we have proved  $D^d A(\delta)=\delta^dA(D)$ when $\delta=kD$, this is not difficult to show that for any $\delta=k\delta^\prime$ with $\delta^\prime\leq D$, ${\delta^\prime}^d A(\delta)=\delta^dA(\delta^\prime)$. As a consequence when $D/\delta=k^\prime/k$, we obtain 
\begin{equation}\label{casrationnel}A(D)=\frac{D^d }{{\delta^\prime}^d}A(\delta^\prime)=\frac{D^d }{\delta^d}A(\delta).\end{equation}

In the general case, one may find a sequence of rational numbers $(q_n)_{n\in\NN}$ which converges to $D/\delta$. Let $\delta_n=q_nD$, we have from (\ref{casrationnel}), $A(D)=\frac{D^d}{\delta_n^d}A(\delta_n)$.  Since we have assumed $\Esp(\Vect Y^2_{\Gamma_n})\to 0$ when $\Gamma_n\to0$, the additivity of $\Vect Y$ and $\delta_n\to\delta$ yield $$A(\delta_n)=\Esp\left(\Vect Y_{\Delta_0(\delta_n)}\tr{\Vect Y_{\overline{\Delta}(\delta_n)}}\right)\to\Esp\left(\Vect Y_{\Delta_0(\delta)}\tr{\Vect Y_{\overline{\Delta}(\delta)}}\right)=A(\delta)$$ as $n$ goes to infinity. Therefore, the identity (\ref{casrationnel}) holds for any $\delta>0$, which concludes the proof.


\subsection{Proof of Proposition~\ref{prop-sdp}}

The proof follows arguments presented by \citeauthor{A-JenKun94} in \cite{A-JenKun94}. Let $C_n(\overline\delta)=[-n\overline{\delta}-\overline{\delta}/2,n\overline{\delta}+\overline{\delta}/2]^d$, so $C_n(\overline\delta)=\cup_{k\in\mathcal K_n}\Delta_k(\overline\delta)$, where $\mathcal K_n=[-n,n]^d\cap\ZZ[d]$ and $\Delta_k(\overline\delta)$ is the cube centered at $k\overline\delta$ with side-length $\overline\delta$. We have 
\begin{align*}
\Var\left( |C_n(\overline\delta)|^{-1/2} \Vect{Y}_{C_n(\overline\delta)}(\Phi;\ParVT)\right)  &= |C_n(\overline\delta)|^{-1} \sum_{i,j\in \mathcal K_n} \Esp\left( 
\Vect{Y}_{\Delta_i(\overline\delta)}\left(\Phi;\ParVT\right) \tr{\Vect{Y}_{\Delta_j(\overline\delta)}\left(\Phi;\ParVT\right) }
\right) \nonumber \\
&= |C_n(\overline\delta)|^{-1} \sum_{i\in \mathcal K_n} \sum_{j \in \mathbbm{B}_i\left( \left\lceil \frac{D}{\overline\delta}\right\rceil \right)\cap \mathcal K_n }\Esp\left( 
\Vect{Y}_{\Delta_i(\overline\delta)}\left(\Phi;\ParVT\right) \tr{\Vect{Y}_{\Delta_j(\overline\delta)}\left(\Phi;\ParVT\right) }
\right) \nonumber. \end{align*}
Since $|C_n(\overline\delta)|=\overline{\delta}^d|\mathcal K_n|$, from the ergodic theorem, 
$$\Var\left( |C_n(\overline\delta)|^{-1/2} \Vect{Y}_{C_n(\overline\delta)}(\Phi;\ParVT)\right) \longrightarrow \overline\delta^{-d}  \sum_{|k|\leq \left\lceil \frac{D}{\overline\delta}\right\rceil}\Esp \left( \Vect{Y}_{\Delta_0(\overline\delta)}\left(\Phi;\ParVT\right) \tr{\Vect{Y}_{\Delta_k(\overline\delta)}\left(\Phi;\ParVT\right) }
\right)$$ which is $\Mat{M}(\ParVT)$ by Lemma \ref{delta}.

Therefore, to prove that $\Mat{M}(\ParVT)$ is positive-definite, it is sufficient to prove that the covariance matrix $\Var\left( |C_n(\overline\delta)|^{-1/2} \Vect{Y}_{C_n(\overline\delta)}(\Phi;\ParVT)\right)$ is  positive-definite for $n$ large enough. Let $\Vect{x}\in \RR[q]\setminus\{0\}$, we must show that  $$V:=\tr{\Vect{x}} \Var\left( |C_n(\overline\delta)|^{-1/2} \Vect{Y}_{C_n(\overline\delta)}(\Phi;\ParVT)\right) \Vect{x} >0.$$
Since, for two random variables $X,X^\prime$ with finite variance 
$$
\Var(X) = \Esp (\Var(X|X^\prime)) \;\;+ \;\; \Var( \Esp(X|X^\prime)) \geq \Esp( \Var(X|X^\prime)),
$$
we have, by denoting $L:= \left(2\cD[\frac{D}{\overline{\delta}}]+1\right)\ZZ[d]$, 

\begin{eqnarray*}
V &\geq& |C_n(\overline{\delta})|^{-1} \; \Esp \left( \Var \left( \tr{\Vect{x}}
\Vect{Y}_{C_n(\overline{\delta})} (\Phi;\ParVT) | \; \Phi_{\Delta_k(\overline{\delta})}, k\notin L
\right) \right) \\
&=& |C_n(\overline{\delta})|^{-1} \tr{\Vect{x}} \Esp\bigg( \Var\bigg(
\sum_{\ell \in L \cap \mathcal K_n} \;\; \underbrace{\sum_{i \in \mathbbm{B}_\ell\left(\left\lceil\frac{D}{\overline{\delta}} \right\rceil \right)\cap \mathcal K_n } \Vect{Y}_{\Delta_i(\overline{\delta})}  (\Phi;\ParVT) }_{:=\Vect{S}_{\ell,n}(\Phi)}
 | \; \Phi_{\Delta_k(\overline{\delta})}, k\notin L \bigg) \bigg) \Vect{x}
\end{eqnarray*}
Note that from the locality property, $\Vect{S}_{\ell,n}(\Phi)$ depends only on $\Phi_{\Delta_{j}(\overline{\delta})}$ for $j \in \mathbbm{B}_\ell\left(2 \left\lceil \frac{D}{\overline{\delta}} \right\rceil\right)$. Therefore, conditionally on $\Phi_{\Delta_k(\overline{\delta})}, k \notin L$, the variables $\Vect{S}_{\ell,n}(\Phi)$ and $\Vect{S}_{\ell^\prime,n}(\Phi)$ (for $\ell\neq \ell^\prime$) are independent. Now, let $\overline\Delta(\overline\delta):=\cup_{|i|\leq \left\lceil \frac{D}{\overline{\delta}} \right\rceil}\Dom[i](\overline{\delta})$, from the stationarity we have for $n$ large enough
\begin{eqnarray*}
V &\geq & |C_n(\overline{\delta})|^{-1} \tr{\Vect{x}} \sum_{\ell \in L\cap \mathcal K_n} \Esp\bigg(\Var\bigg(
\Vect{S}_{\ell,n}(\Phi) 
| \; \Phi_{\Delta_k(\overline{\delta})}, k\notin L \bigg) \bigg) \Vect{x}\\
&\geq &\frac{{\overline{\delta}}^{-d}}2 \frac{
|L\cap \mathcal K_n| }{|\mathcal K_n|}\times
\Esp  \left( \Var \left( \tr{\Vect{x}} \ensuremath{\Vect{Y}_{\overline\Delta\left(\overline\delta\right)} \left( \Phi ; \ParVT  \right)} \big| \Phi_{\Dom[k](\overline{\delta})}, 1\leq|k|\leq2\cD[\frac{D}{\overline{\delta}}] \right)\right) \\
&\geq& \kappa(\overline{\delta},D,d)  \times 
\Esp  \left( \Var \left( \tr{\Vect{x}} \ensuremath{\Vect{Y}_{\overline\Delta\left(\overline\delta\right)}\left( \Phi ; \ParVT  \right)} \big| \Phi_{\Dom[k](\overline{\delta})}, 1\leq|k|\leq2 \cD[\frac{D}{\overline{\delta}}]\right)\right) ,
\end{eqnarray*}
where $\kappa(\overline{\delta},D,d)$ is a positive constant. Assume there exists some positive constant $c$ such that $P_{\ParVT}-$a.s. $\tr{\Vect{x}} \ensuremath{\Vect{Y}_{\overline\Delta\left(\overline\delta\right)} \left( \Phi ; \ParVT  \right)}=c$ when the variables $\Phi_{\Dom[k](\overline{\delta})}, 1\leq |k|\leq 2\cD[\frac{D}{\overline{\delta}}]$ are  fixed to belong to $B$, where $B\in\mathcal F$ is involved in \textbf{[PD]}. It follows that for any $\varphi_i \in A_i$ for $i=0,\ldots,\ell$ (with $\ell\geq 1$), where the $A_i$'s come from \textbf{[PD]}, $\tr{\Vect{x}} \left( \ensuremath{\Vect{Y}_{\overline\Delta\left(\overline\delta\right)} \left( \varphi_i ; \ParVT  \right)}-\ensuremath{\Vect{Y}_{\overline\Delta\left(\overline\delta\right)} \left( \varphi_0 ; \ParVT  \right)}\right)=0$. Since for all $\left(\varphi_0,\ldots,\varphi_{\ell}\right) \in A_0\times\ldots\times A_{\ell}$, the matrix with entries $\left(\ensuremath{\Vect{Y}_{\overline\Delta\left(\overline\delta\right)} \left( \varphi_i ; \ParVT  \right)}\right)_j-\left(\ensuremath{\Vect{Y}_{\overline\Delta\left(\overline\delta\right)} \left( \varphi_0 ; \ParVT  \right)}\right)_j$  is assumed to be injective, this leads to $\Vect{x}=0$ and hence to some contradiction. Therefore, when the variables $\Phi_{\Dom[k](\overline{\delta})}$, $1\leq |k|\leq 2\cD[\frac{D}{\overline{\delta}}]$ are for example assumed to belong to $B$, the variable $ \tr{\Vect{x}}\ensuremath{\Vect{Y}_{\overline\Delta\left(\overline\delta\right)} \left( \Phi ; \ParVT  \right)}$ is almost surely not a constant and so $V>0$, which proves that $\Mat{M}(\ParVT)$ is a symmetric positive-definite matrix.


\subsection{Proof of Proposition~\ref{consistance}}\label{preuve-consistance}

Since for any $\varphi\in \Omega$, $\widehat{\Mat{M}}_n(\varphi; \cdot,\delta,D^\vee)$ is continuous in a neighborhood $\mathcal{V}(\ParVT)$ of $\ParVT$ and according to \textbf{[E1]}, it is sufficient to prove that for any $\ParV \in \mathcal{V}(\ParVT)$, $\widehat{\Mat{M}}_n(\Phi; \ParV,\delta_n,D^\vee)$ converges in probability towards $\Mat{M}(\ParV)$.

We choose the sequence $\delta_n$ as follows :
$$\delta_n= \frac{|\Lambda_n|^{1/d}}{\left\lfloor \frac{|\Lambda_n|^{1/d}}{\delta}\right\rfloor},$$
which guarantees $\delta_{n_0}=\delta$, since $|\Lambda_{n_0}|\delta^{-d}\in\NN$, $\delta\leq\delta_n \leq 2\delta$ for $n$ sufficiently large, and $\delta_n\to\delta$ as $n\to\infty$. This choice allows us to consider, for all $n\in\NN$, the decomposition $\Lambda_{n}=\cup_{k\in \mathcal K_{n}} \Delta_k(\delta_n)$, where  the $\Delta_k(\delta_n)$'s are disjoint cubes with side-length $\delta_n$ and centered at $k\delta_n$. Moreover, since $\delta_n\geq\delta$ and $\delta_n\to\delta$, we have $\left\lceil \frac{D^\vee}{\delta_n} \right\rceil=\left\lceil \frac{D^\vee}{\delta} \right\rceil$ when $n$ is large enough, which is assumed in the sequel of the proof.

Let $\widetilde{\mathcal K}_{n}:= \mathcal K_{n} \cap \left( \cup_{j\in \partial \mathcal K_{n}} \mathbbm B_j\left(\left\lceil \frac{D^\vee}{\delta} \right\rceil \right)\right)$.  Since $|\Lambda_n|=\delta_n^d |\mathcal K_n|$, we have
$$\left| \delta_n^d \widehat{\Mat{M}}_{n}(\Phi; \ParV,\delta_n,D^\vee)- \delta^d \Mat{M}(\ParV)\right|
\leq X_1+X_2+X_3+X_4,$$
where by setting $\overline{\Delta}_k(\tau)=\cup_{j\in \mathbbm{B}_k\left(\left\lceil \frac{D^\vee}{\delta} \right\rceil \right)}\Delta_j(\tau)$ (for some $\tau>0$),
\begin{eqnarray*}
X_1&=&\left| |\mathcal K_n|^{-1}\sum_{k\in \widetilde{\mathcal K}_n}
\sum_{j\in \mathbbm{B}_k\left(\left\lceil \frac{D^\vee}{\delta} \right\rceil \right) \cap \mathcal K_n}
 \widehat{\Vect{Y}}_{n,\Delta_k(\delta_n)}\left(\Phi;\ParV\right) \tr{\widehat{\Vect{Y}}_{n,\Delta_j(\delta_n)}\left(\Phi;\ParV\right)}
\right| \\
X_2&=&\left||\mathcal K_n|^{-1}\sum_{k\in \mathcal K_n\setminus \widetilde{\mathcal K}_{n}}
\left( \widehat{\Vect{Y}}_{n,\Delta_k(\delta_n)}\left(\Phi;\ParV\right) \tr{\widehat{\Vect{Y}}_{n,\overline{\Delta}_k(\delta_n)}\left(\Phi;\ParV\right)}-
{\Vect{Y}}_{\Delta_k(\delta_n)}\left(\Phi;\ParV\right) \tr{{\Vect{Y}}_{\overline{\Delta}_k(\delta_n)}\left(\Phi;\ParV\right)}\right)\right| \\
X_3 &=&  \left||\mathcal K_n|^{-1}\sum_{k\in \mathcal K_n\setminus \widetilde{\mathcal K}_{n}}
\left( {\Vect{Y}}_{n,\Delta_k(\delta_n)}\left(\Phi;\ParV\right) \tr{{\Vect{Y}}_{n,\overline{\Delta}_k(\delta_n)}\left(\Phi;\ParV\right)}-
{\Vect{Y}}_{\Delta_k(\delta)}\left(\Phi;\ParV\right) \tr{{\Vect{Y}}_{\overline{\Delta}_k(\delta)}\left(\Phi;\ParV\right)}\right)\right| \\
X_4 &= &\left||\mathcal K_n|^{-1}\sum_{k\in \mathcal K_n\setminus \widetilde{\mathcal K}_{n}} 
 {\Vect{Y}}_{\Delta_k(\delta)}\left(\Phi;\ParV\right) \tr{{\Vect{Y}}_{\overline{\Delta}_k(\delta)}\left(\Phi;\ParV\right)}
-
\Esp \left( \Vect{Y}_{\Delta_0(\delta)}\left(\Phi;\ParV\right) \tr{\Vect{Y}_{\overline{\Delta}_0(\delta)}\left(\Phi;\ParV\right) }\right)
\right|.
\end{eqnarray*}
We have  from the additivity and the stationartiy of $\widehat{\Vect{Y}}_n$,
\begin{align*}
\Esp |X_1|&\leq |\mathcal K_n|^{-1}\sum_{k\in \widetilde{\mathcal K}_{n}} \sum_{j\in \mathbbm{B}_k\left(\left\lceil \frac{D^\vee}{\delta} \right\rceil \right) \cap \mathcal K_{n}} \Esp\left|\widehat{\Vect{Y}}_{n,\Delta_k(\delta_n)}\left(\varphi;\ParV\right) \tr{\widehat{\Vect{Y}}_{n,\Delta_j(\delta_n)}\left(\varphi;\ParV\right)}\right|\\
&\leq \frac{|\widetilde{\mathcal K}_{n}|}{|\mathcal K_n|} \Esp \left|
\widehat{\Vect{Y}}_{n,\Delta_0(\delta_n)}(\Phi;\ParV)\tr{\widehat{\Vect{Y}}_{n,\overline{\Delta}_0(\delta_n)}(\Phi;\ParV)} \right|,
\end{align*}
which tends to 0 as $n\to +\infty$, because $\frac{|\widetilde{\mathcal K}_{n}|}{|\mathcal K_{n}|} \to 0$ and $\delta\leq\delta_n\leq 2\delta$. Therefore, $X_1$ converges in probability to 0. The second term converges also to 0 in probability from the additivity of $\Vect{Y}$ and $\widehat{\Vect{Y}}$ and from~\eqref{eq-hypYhat}. The expectation of the third term converges to 0 by following the proof of Lemma~\ref{dropDn}. Finally, from the stationarity of $\Vect Y$ and since $|\mathcal K_n|\sim|\mathcal K_n\setminus \widetilde{\mathcal K}_{n}|$, the mean ergodic theorem applies to $\Esp|X_4|$, which, in particular, shows that $X_4\to 0$ in probability. This proves that $$\delta_n^d \widehat{\Mat{M}}_{n}(\Phi; \ParV,\delta_n,D^\vee) \longrightarrow  \delta^d \Mat{M}(\ParV),$$ 
in probability, as $n\to\infty$. Since $\delta_n$ is a deterministic sequence converging to $\delta$, the conclusion of Proposition \ref{consistance} follows.


\appendix

\section{Central Limit Theorem}\label{annexe-tcl}
The following result is a central limit theorem for conditionnally centered random fields. It generalizes Theorem 2.1 in \cite{A-JenKun94} to a non-stationary and non-ergodic setting. A general result has been proved by \cite{A-ComJan98} for self normalized sums, provided a fourth moment condition. Our result is in the same spirit but it is proved for triangular array and without self-normalization, which is well-adapted to the residuals process framework. This allows in particular to avoid the fourth moment assumption.

\begin{theorem}\label{tcl}
Let $X_{n,i}$, $n\in\NN$, $i\in\ZZ[d]$, be a triangular array field in a measurable space $S$. For $n\in\NN$, let $\mathcal K_n\subset\ZZ[d]$ and for $k\in \mathcal K_n$, assume 
\begin{equation}\label{subordinate}\Vect Z_{n,k}=f_{n,k}\left(X_{n,k+i},\ i\in \mathcal I_0\right),\end{equation} where $\mathcal I_0=\{i\in\ZZ[d],\ |i|\leq 1\}$ and $f_{n,k} : S^{\mathcal I_0} \to \RR[p]$. Let $\Vect S_n=\sum_{k\in \mathcal K_n} \Vect Z_{n,k}$. If
\begin{itemize}
\item[(i)] $c_3:=\sup_{n\in\NN} \sup_{k\in \mathcal K_n} \Esp|\Vect Z_{n,k}|^3 <\infty$,
\item[(ii)] $\forall n\in\NN$, $\forall k\in \mathcal K_n$, $\Esp(\Vect Z_{n,k} | X_{n,j},\ j\not=k)=0,$
\item[(iii)] $|\mathcal K_n|\to +\infty$ as $n\to\infty$,
\item[(iv)] There exists a symmetric matrix $\Mat\Sigma\geq 0$ such that $$\Esp \left\| |\mathcal K_n|^{-1} \sum_{k\in \mathcal K_n}\sum_{j\in \mathbbm B_k(1)\cap \mathcal K_n} \Vect Z_{n,k}\tr{\Vect Z_{n,j}} - \Mat\Sigma \right\|\to 0,$$
\end{itemize}
then $|\mathcal K_n|^{-1/2}\Vect S_n\stackrel{d}{\longrightarrow} \mathcal{N}\left( 0, \Mat\Sigma\right)$ as $n\to\infty$.
\end{theorem}

\begin{proof}

Let us first assume that $\Mat\Sigma$ is a positive-definite matrix (i.e. $\Mat\Sigma> 0$). According to the Stein's method (see also \cite{A-Bol82}), it suffices to prove that, for all $\Vect e\in\RR[p]$ such that $\|\Vect e\|=1$ and for all $\lambda\in\RR[]$, 
$$\Esp\left(\left(i\lambda-\tr{\Vect e} |\mathcal K_n|^{-1/2}\Mat \Sigma^{-1/2}\Vect S_n\right)\exp\left(i\lambda \tr{\Vect e} |\mathcal K_n|^{-1/2}\Mat\Sigma^{-1/2}\Vect S_n\right)\right)\to 0.$$
Denoting $\Vect u=\lambda\Vect e$, this is equivalent to prove that for all $\Vect u\in\RR[p]$, 
$$\Esp\left(\underbrace{\left( i\Vect u-|\mathcal K_n|^{-1/2}\Mat\Sigma^{-1/2}\Vect S_n\right)\exp(i\tr{\Vect u}|\mathcal K_n|^{-1/2}\Mat\Sigma^{-1/2}\Vect S_n)}_{:=\Vect{A}}\right)\to  0.$$

We decompose the term $\Vect{A}$ in the same spirit as in \cite{A-Bol82}, \cite{A-JenKun94} and \cite{A-ComJan98}. Let us denote by $\Mat I_p$ the identity matrix of size $p$ and $\Vect S_n^k=\sum_{j\in \mathbbm B_k(1)\cap \mathcal K_n}\Vect Z_{n,j}$. Noting that $\tr{\Vect u}\Mat \Sigma^{-1/2}\Vect S_n^k=\tr{\Vect S_n^{k}}\tr{\Mat \Sigma^{-1/2}}\Vect u$, the decomposition is $\Esp(\Vect{A})=\Esp(\Vect A_1-\Vect A_2-\Vect A_3)$ where
\begin{align*}
\Vect A_1&= i\exp(i\tr{\Vect u}|\mathcal K_n|^{-1/2}\Mat \Sigma^{-1/2}\Vect S_n)\left[\Mat I_p-|\mathcal K_n|^{-1}\Mat \Sigma^{-1/2}\sum_{k\in \mathcal K_n} \Vect Z_{n,k}\tr{\Vect{S}_n^k}
\right]\Vect u , \\ 
\Vect A_2&=\exp(i\tr{\Vect u}|\mathcal K_n|^{-1/2}\Mat \Sigma^{-1/2}\Vect S_n)|\mathcal K_n|^{-1/2}\Mat \Sigma^{-1/2}\\
 &\hspace{2.5cm}\times\sum_{k\in \mathcal K_n} \Vect Z_{n,k}\left(1- \exp(-i\tr{\Vect u}|\mathcal K_n|^{-1/2}\Mat \Sigma^{-1/2}\Vect S_n^k)-i\tr{\Vect u}|\mathcal K_n|^{-1/2}\Mat \Sigma^{-1/2}\Vect S_n^k\right), \\
\Vect A_3&= |\mathcal K_n|^{-1/2}\Mat \Sigma^{-1/2}\sum_{k\in \mathcal K_n} \Vect Z_{n,k} \exp\left[i\tr{\Vect u}|\mathcal K_n|^{-1/2}\Mat \Sigma^{-1/2}(\Vect S_n-\Vect S_n^k)\right].
\end{align*}
The two last terms $\Vect A_2$ and $\Vect A_3$ can be handled as in \cite{A-JenKun94}: $\Esp(\Vect A_3)=0$ by $(ii)$ and $|\Esp(\Vect A_2)|\to 0$ from the same inequalities therein and the sub-multiplicative property of the Frobenius norm. These inequalities rely on two facts: $\forall x\in\RR[], |1-e^{-ix}-ix|\leq x^2/2$ and for all $n$ and for all $(k_1,k_2,k_3)\in \mathcal K_n$, $\Esp(|\Vect Z_{n, k_1}||\Vect Z_{n, k_2}||\Vect Z_{n, k_3}|)\leq \left(\Esp(|\Vect Z_{n, k_1}|^3)\Esp(|\Vect Z_{n, k_2}|^3)\Esp(|\Vect Z_{n, k_3}|^3)\right)^{1/3}$ which is less than $c_3$ by $(i)$.

For $\Vect A_1$, we cannot use a mean ergodic theorem as in \cite{A-JenKun94}, but Assumption $(iv)$ is sufficient. Indeed,
\begin{align*} \|\Esp(\Vect A_1)\| &\leq \|\Vect u\| \Esp\left\| \Mat I_p-|\mathcal K_n|^{-1}\Mat \Sigma^{-1/2}\sum_{k\in \mathcal K_n}\sum_{j\in \mathbbm B_k(1)\cap \mathcal K_n} \Vect Z_{n,k}\tr{\Vect Z_{n,j}}\tr{\Mat \Sigma^{-1/2}}\right\|\\ 
&\leq \|\Vect u\|\left\|\Mat \Sigma^{-1/2}\right\|^2\Esp\left\| |\mathcal K_n|^{-1}\sum_{k\in \mathcal K_n}\sum_{j\in \mathbbm B_k(1)\cap \mathcal K_n} \Vect Z_{n,k}\tr{\Vect Z_{n,j}}-\Mat \Sigma\right\|
\end{align*}
which tends to 0 by $(iv)$.

Now, if  $\Mat\Sigma$ is not a positive-definite matrix, one can find an orthonormal basis $(\Vect f_1,\dots,\Vect f_p)$ of $\RR[p]$, where the $\Vect f_i$'s are eigenvectors of $\Mat\Sigma$. We agree that, if $r<p$ denotes the rank of $\Mat\Sigma$, then $(\Vect f_1,\dots,\Vect f_r)$ is a basis of the image of $\Mat\Sigma$, while $(\Vect f_{r+1},\dots,\Vect f_p)$ is a basis of its kernel. 

Let us denote by $\Mat V_{Im}$ the matrix whose columns are $(\Vect f_1,\dots,\Vect f_r)$ and $\Mat V_{Ker}$ the matrix whose columns are $(\Vect f_{r+1},\dots,\Vect f_p)$. Similarly, for any $\Vect u\in\RR[p]$, let us denote by $u_i$ its $i$-th coordinate in the basis $(\Vect f_1,\dots,\Vect f_p)$  and $\Vect u_{Im}=(u_1,\dots,u_r)$, $\Vect u_{Ker}=(u_{r+1},\dots,u_p)$. Hence $\Vect u=\Mat V_{Im}\Vect u_{Im}+\Mat V_{Ker}\Vect u_{Ker}$.

 The convergence in law of $|\mathcal K_n|^{-1/2}\Vect S_n$ to a Gaussian vector reduces to the convergence of $\tr{\Vect u}|\mathcal K_n|^{-1/2}\Vect S_n$ for all $\Vect u\in\RR[p]$. We have
\begin{equation}\label{decomposition}
\tr{\Vect u}|\mathcal K_n|^{-1/2}\Vect S_n=\tr{\Vect u_{Im}}\tr{\Mat V_{Im}}|\mathcal K_n|^{-1/2}\Vect S_n + \tr{\Vect u_{Ker}}\tr{\Mat V_{Ker}}|\mathcal K_n|^{-1/2}\Vect S_n.
\end{equation}
From $(iv)$ and since $\tr{\Mat V_{Ker}}\Mat\Sigma\ \Mat V_{Ker}=0$, we deduce that  $$\Esp \left\| |\mathcal K_n|^{-1} \sum_{k\in \mathcal K_n}\sum_{j\in \mathbbm B_k(1)\cap \mathcal K_n}\tr{\Mat V_{Ker}}\Vect Z_{n,k}\tr{\Vect Z_{n,j}}\Mat V_{Ker}\right\|\longrightarrow0,$$ which means that $\tr{\Mat V_{Ker}}|\mathcal K_n|^{-1/2}\Vect S_n$ tends to $0$ in quadratic mean.

On the other hand, the assumptions of Theorem \ref{tcl} imply that $(i)-(iv)$ remain true when one replaces $\Vect Z_{n,k}$ by $\tr{\Mat V_{Im}}\Vect Z_{n,k}$ and $\Mat\Sigma$ by $\tr{\Mat V_{Im}}\Mat\Sigma\ \Mat V_{Im}$. Since  $\tr{\Mat V_{Im}}\Mat\Sigma\ \Mat V_{Im}$ is positive-definite, the convergence in law 
of  $\tr{\Mat V_{Im}}|\mathcal K_n|^{-1/2}\Vect S_n$  holds for the same reasons as in the first part of the proof. 

Therefore, we have proved that for all $\Vect u\in\RR[p]$, $\tr{\Vect u}|\mathcal K_n|^{-1/2}\Vect S_n\stackrel{d}{\longrightarrow} \mathcal{N}\left( 0, \tr{\Vect u_{Im}}\tr{\Mat V_{Im}}\Mat\Sigma\ \Mat V_{Im}\Vect u_{Im}\right)$. It is easy to check that $\tr{\Vect u_{Im}}\tr{\Mat V_{Im}}\Mat\Sigma\ \Mat V_{Im}\Vect u_{Im} = \tr{\Vect u}\Mat\Sigma\Vect u$, which concludes the proof.
\end{proof}


\section{Assumption \textbf{[PD]} on two examples}\label{exemplesPD}

In this section, we focus on the two following models, belonging to the exponential family :
\begin{itemize}
\item[1.]  {\bf Two-type marked Strauss point process}~: $\mathbb{M}=\{1,2\}$ and for $\ParV=\tr{\left(\theta_1^{1},\theta_1^{2},\theta_2^{1,1},\theta_2^{1,2},\theta_2^{2,2} \right)}$, for any $\Lambda\in\sB(\RR[d])$,
$$
V_{\Lambda}\left( \varphi ; \ParV \right) = \theta_1^{1} \underbrace{|\varphi_{\Lambda}^{1}|}_{:=v_{\Lambda,1}^{1}(\varphi)} + \theta_1^{2} \underbrace{|\varphi_{\Lambda}^{2}|}_{:=v_{\Lambda,1}^{2}(\varphi)} + \sum_{\begin{subarray}{c}m_1,m_2=1\\ m_1\leq m_2\end{subarray}}^2 \theta_2^{m_1,m_2} \underbrace{\sum_{\begin{subarray}{l}\{x_1^{m_1},x_2^{m_2}\} \in \mathcal{P}_2(\varphi)\\ \{x_1^{m_1},x_2^{m_2}\}\cap\Lambda\not=\emptyset \end{subarray}} \mathbf{1}_{[0,D^{m_1,m_2}]} \left( \|x_2-x_1\|\right)}_{:=v_{\Lambda,2}^{m_1,m_2}(\varphi)}.
$$
Alternatively,
$$ \VIPar{x^m}{\varphi}{\ParV}= \theta_1^m + \sum_{m^\prime=1}^2 \theta_2^{m,m^\prime}\sum_{y^{m^\prime} \in \varphi} \mathbf{1}_{[D_0^{m_1,m_2}, D^{m_1,m_2}]}(\|y-x\|).$$
\end{itemize}

This process is well-defined when $\theta_2^{m_1,m_2}\geq 0$ and $D_0^{m_1,m_2}=0$ (inhibition assumption), or when $\theta_2^{m_1,m_2} \in \RR$ and $D_{0}^{m_1,m_2}=\delta>0$ (hard-core assumption), see  Proposition \ref{prop-C12N13} for instance. The range of the local energy function equals $D=\max\left( D^{1,1},D^{1,2},D^{2,2}\right)$.  

\begin{itemize}
\item[2.] {\bf Area interaction point process}~: $\mSp=\{0\}$ and for $R>0$, $\ParV=(\theta_1 , \theta_2)$ and any $\Lambda\in\sB(\RR[d])$,
$$
V_{\Lambda}\left( \varphi ; \ParV \right) = \theta_1 |\varphi_{\Lambda}| + \theta_2 \SEx{2}{\varphi_{\Lambda}}, \;
\mbox{ with } \;
\SEx{2}{\varphi_{\Lambda}} := \left| \cup_{x\in \varphi_{\Lambda}} B(x,R)\right|.
$$

Alternatively,
$$
\SExI{1}{0}{\varphi} :=1,\qquad\SExI{2}{0}{\varphi} := \left| \cup_{x\in (\varphi_{\mathcal{B}(0,2R)}\cup\{0\})} \mathcal{B}(x,R) \setminus \cup_{x\in \varphi_{\mathcal{B}(0,2R)}} \mathcal{B}(x,R)\right|.
$$
\end{itemize}
This model is well-defined for $\ParV \in \RR[2]$ (see  Proposition \ref{prop-C12N13} for instance) and the range of the local energy equals $D=2R$.

Both these models satisfy \textbf{[C]} and \textbf{[N1-4]}. The aim of the sequel is to prove Proposition \ref{prop-exemples}, which claims that  $\Mat{\Sigma}_1(\ParVT)$ and $\Mat{\Sigma}_2(\ParVT)$, involved respectively in Proposition~\ref{prop-fwk1} and \ref{prop-fwk2}, are positive-definite for these models, when considering the maximum pseudolikelihood estimate for $\widehat{\ParV}_n$ and the two following frameworks
\begin{itemize}
\item Framework~1 (for $\Mat{\Sigma}_1(\ParVT)$): we consider the inverse residuals ($h=e^V$). Let us recall that Proposition~\ref{prop-PDfails} asserts that \textbf{[PD]} fails for the raw residuals ($h=1$) for both the area-interaction and 2-type marked Strauss models.
\item Framework~2 (for $\Mat{\Sigma}_2(\ParVT)$): we consider the family of test functions given for $j=1,\ldots,s$ and $0<r_1<\ldots<r_s<+\infty$ by $$h_j(x^m,\varphi;\ParV)= \mathbf{1}_{[0,r_j]}(d(x^m,\varphi)) e^{V\left(x^m|\varphi;\ParV\right)},$$
related to parametric and nonparametric estimations of the empty space function at distance $r_j$.
\end{itemize}

When considering the MPLE, $R_{\infty,\Lambda}(\varphi;h,\ParVT)$ is given by~\eqref{eq-RinftyMPLE} with $\Vect{LPL}^{(1)}$, $\Mat{H}$ and $\Vect{\mathcal{E}}$ respectively given by~\eqref{eq-defLPL1}, \eqref{eq-defH} and~\eqref{eq-defE}.

\subsection{2-type marked Strauss point process}
We only deal with the inhibition case, that is $\SpPar=\RR[2] \times \mathbb{R}^{3}_+$ and $D_0^{m_1,m_2}=0$ .
The following proofs could  easily be extended to the hard-core case and to the multi-Strauss marked point process (see {\it e.g. } \cite{A-BilCoeDro08}). For any vector $\Vect{z}$ of length 5, we sometimes reparameterize it similarly as the parameter vector, that is $\Vect{z} =\tr{(z_1^{1},z_1^{2},z_2^{1,1},z_2^{1,2},z_2^{2,2})}$.

\subsubsection{Proof that $\Mat{\Sigma}_1(\ParVT)$ is positive-definite for the two-type Strauss model}
From Proposition~\ref{prop-sdp}, proving that $\Mat{\Sigma}_1(\ParVT)$ is positive-definite in Framework 1 reduces to check  Assumption \textbf{[PD]} with  $h=e^V$ and 
\begin{itemize}
\item[$(i)$] $\Vect{Y}_\Lambda(\varphi;\ParVT)=I_\Lambda(\varphi;e^V,\ParVT)$,
\item[$(ii)$] $\Vect{Y}_\Lambda(\varphi;\ParVT)=R_{\infty,\Lambda} (\varphi;e^V,\ParVT)$.
\end{itemize}

$(i)$ is ensured by Proposition~\ref{prop-lInn}. \\


$(ii)$ We  fix $\overline\delta=D$ and $B=\emptyset$ in  \textbf{[PD]}. Let $\overline{\Omega}:=\overline{\Omega}_{\emptyset}$. Without loss of generality, one may assume that ${\theta_2^\star}^{1,1}>0$. Let us define for $n\geq 1$ 
\begin{multline*}
A_{n,-}(\eta)=\left\{ 
\varphi\in \overline{\Omega}: \varphi(\Delta_0(\overline{\delta})\times \{1\})=2n, \varphi(\Delta_0(\overline{\delta})\times\{2\})=0, \right. \\
\left.\varphi\left(\mathcal{B}\left((0,0),\frac{\eta}4 \right)\right)=n, \varphi\left(\mathcal{B}\left((D^{1,1}-\frac{\eta}2,0),\frac{\eta}4 \right)\right)=n
\right\},
\end{multline*}
\begin{multline*}
A_{n,+}(\eta)= \left\{ 
\varphi\in \overline{\Omega}: \varphi(\Delta_0(\overline{\delta})\times\{1\})=2n, \varphi(\Delta_0(\overline{\delta})\times\{2\})=0, \right. \\
\left.\varphi\left(\mathcal{B}\left((0,0),\frac{\eta}4 \right)\right)=n, \varphi\left(\mathcal{B}\left((D^{1,1}+\frac{\eta}2,0),\frac{\eta}4 \right)\right)=n
\right\} .
\end{multline*}


Let $\varphi_{n,-}\in A_{n,-}$ and $\varphi_{n,+}\in A_{n,+}$. Then for $\eta$ small enough
\begin{eqnarray*}
I_{\Ltt}(\varphi_{n,\bullet};e^V,\ParVT) &=& |\Ltt| - \left\{
\begin{array}{ll}
2n e^{{\theta_1^\star}^{1} + (2n-1){\theta_2^\star}^{1,1} } & \mbox{ if } \bullet=-,\\
2n e^{{\theta_1^\star}^{1} + (n-1){\theta_2^\star}^{1,1} } & \mbox{ if } \bullet=+.\\
\end{array}\right.
\end{eqnarray*}

\begin{eqnarray*}
\left( \Vect{LPL}^{(1)}_{\Ltt}(\varphi_{n,\bullet};\ParVT) \right)_{1}^{m^\prime} &=& \ism[\Ltt\times \mSp]{v_1^{m^\prime}(x^m|\varphi_{n,\bullet}) e^{-\VIPar{x^m}{\varphi_{n,\bullet}}{\ParVT}  } } 
- \left\{ \begin{array}{ll}
2n & \mbox{ if } m^\prime=1,\\
0 & \mbox{ if } m^\prime=2.
\end{array} \right.\\
\left( \Vect{LPL}^{(1)}_{\Ltt}(\varphi_{n,-};\ParVT) \right)_{2}^{m^\prime_1,m^\prime_2} &=& \ism[\Ltt\times \mSp]{v_2^{m_1,m_2}(x^m|\varphi_{n,-}) e^{-\VIPar{x^m}{\varphi_{n,-}}{\ParVT}  } } \\
&& \hspace*{5cm}-\left\{ \begin{array}{ll}
2n(2n-1) & \mbox{ if } m_1=m_2=1,\\
0 & \mbox{ otherwise}.
\end{array} \right.\\
\left( \Vect{LPL}^{(1)}_{\Ltt}(\varphi_{n,+};\ParVT) \right)_{2}^{m_1,m_2} &=& \ism[\Ltt\times \mSp]{v_2^{m_1,m_2}(x^m|\varphi_{n,+}) e^{-\VIPar{x^m}{\varphi_{n,+}}{\ParVT}  } } \\
&& \hspace*{5cm}- \left\{ \begin{array}{ll}
n(n-1) & \mbox{ if } m_1=m_2=1,\\
0 & \mbox{ otherwise}.
\end{array} \right.
\end{eqnarray*}


Now,
\begin{align*}
&\Delta R_{\infty,\Ltt}(\varphi_{n,-},\varphi_{n,+}) := R_{\infty,\Ltt}(\varphi_{n,-};e^V,\ParVT) - R_{\infty,\Ltt}(\varphi_{n,+};e^V,\ParVT) \\ 
&=2n\left( e^{  {\theta_1^\star}^{1}+ (n-1){\theta_2^\star}^{1,1} } -e^{{\theta_1^\star}^{1}+ (2n-1){\theta_2^\star}^{1,1} }\right) +\left(\Vect{W}(e^V,\ParVT)\right)_{2}^{1,1} \left( 2n(2n-1)-n(n-1) \right) \\ 
&\qquad+ f(\varphi_{n,-},\varphi_{n,+},\Vect{W},\eta) \\
&= \underbrace{2ne^{ {\theta_1^\star}^{1}+ (n-1){\theta_2^\star}^{1,1}}(1-e^{{n\theta_2^\star}^{1,1}}) + n(3n-1)\left(\Vect{W}(e^V,\ParVT)\right)_2^{1,1} }_{:=x_n} + f(\varphi_{n,-},\varphi_{n,+},\Vect{W},\eta).
\end{align*}
Fix $\varepsilon>0$, there exists $n_0\geq 1$ such that for all $n\geq n_0$, $x_n < -\varepsilon$. Now by a continuity argument, there exists $\eta_0=\eta_0(n_0)$ such that for all $\eta \leq \eta_0(n_0)$, $|f(\varphi_{n_0,-},\varphi_{n_0,+},\Vect{W},\eta)|\leq \varepsilon/2$. Therefore by assuming that $\Delta R_{\infty,\Ltt}(\varphi_{n_0,-},\varphi_{n_0,+})=0$, we obtain for $\eta\leq \eta_0$
$$
0= |\Delta R_{\infty,\Ltt}(\varphi_{n_0,-},\varphi_{n_0,+})| \geq |x_{n_0}| -  |f(\varphi_{n_0,-},\varphi_{n_0,+},\Vect{W},\eta)| \geq \varepsilon/2>0
$$
which leads to a contradiction and proves \textbf{[PD]}.

\subsubsection{Proof that $\Mat{\Sigma}_2(\ParVT)$ is positive-definite for the two-type Strauss model}\label{sec-DP1fwk2Str}

From Proposition~\ref{prop-sdp}, proving that $\Mat{\Sigma}_2(\ParVT)$ is positive-definite in Framework 2 reduces to check  Assumption \textbf{[PD]} with 
$$\Vect{Y}_\Lambda(\varphi;\ParVT)=\Vect{R}_{\infty,\Lambda}(\varphi;\Vect{h};\ParVT),$$
where, for all $j=1,\dots,s$, $\left(\Vect{R}_{\infty,\Lambda}(\varphi;\Vect{h};\ParVT)\right)_j=R_{\infty,\Lambda}(\varphi;h_j,\ParVT)$, $h_j$ is the test function given by  $h_j(x^m,\varphi;\ParV)= \mathbf{1}_{[0,r_j]}(d(x^m,\varphi)) e^{V\left(x^m|\varphi;\ParV\right)}$. We fix as before $\overline\delta=D$
 and  $B=\emptyset$ in \textbf{[PD]}.

Let $0<r_1<\ldots<r_s<+\infty$. Let us also assume that $r_i\neq D$ for $i=1,\ldots,s$ and define
\begin{eqnarray*}
\!\!A_{i,-}^{1,1}(\eta) \!\!\!\!\!&=&\!\!\!\! \!\!\Big\{ \varphi \in \cSptt : \varphi({\Dom[0](\Dtt)})=2, \varphi\left( \mathcal{B}\left( (0,0),\frac{\eta}4 \right)\times\{1\} \right)=1, \varphi\left( \mathcal{B}\left( (r_i-\frac\eta2,0),\frac{\eta}4 \right)\times\{1\} \right)=1\Big\},  \\
A_{i,+}^{1,1}(\eta) \!\!\!\!&=&\!\!\!\!\!\! \Big\{ \varphi \in \cSptt : \varphi({\Dom[0](\Dtt)})=2, \varphi\left( \mathcal{B}\left( (0,0),\frac{\eta}4 \right)\times\{1\} \right)=1, \varphi\left( \mathcal{B}\left( (r_i+\frac\eta2,0),\frac{\eta}4 \right)\times\{1\} \right)=1\Big\}.
\end{eqnarray*}
Let $\varphi_{i,\bullet} \in A_{i,\bullet}^{1,1}(\eta)$ for $\bullet=-,+$ and i$=1,\ldots,s$. Let $\kappa_i$ the constant given by
$$
\kappa_i = \left\{ \begin{array}{ll}
2 e^{{\theta_1^\star}^{1}+{\theta_2^\star}^{1,1} }   & \mbox{ if } r_i<D,\\
2 e^{{\theta_1^\star}^{1} }   & \mbox{ otherwise.}
\end{array} \right.
$$
Then for $i,j=1,\ldots,s$ and for $\eta$ small enough
\begin{eqnarray*}
I_{\Ltt} ( \varphi_{i,-};h_j,\ParVT) = \ism[\Ltt\times \mSp]{h_j(x^m,\varphi_{i,-})e^{-\VIPar{x^m}{\varphi_{i,-} }{\ParVT} }} -
\left\{ \begin{array}{ll}
\kappa_i & \mbox{ if } i\leq j,\\
0& \mbox{ otherwise.}
\end{array} \right. \\
I_{\Ltt} ( \varphi_{i,+};h_j,\ParVT) = \ism[\Ltt\times \mSp]{h_j(x^m,\varphi_{i,+})e^{-\VIPar{x^m}{\varphi_{i,+} }{\ParVT} }} -
\left\{ \begin{array}{ll}
\kappa_i & \mbox{ if } i< j,\\
0& \mbox{ otherwise.}
\end{array} \right. 
\end{eqnarray*}
On the other hand
\begin{align*}
&\left(\Vect{LPL}^{(1)}_{\Ltt} ( \varphi_{i,\bullet};\ParVT)\right)_1^{m^\prime} = \ism[\Ltt\times \mSp]{ v_1^{m^\prime}(x^m|\varphi_{i,\bullet} e^{-\VIPar{x^m}{\varphi_{i,\bullet} }{\ParVT} }} -
\left\{ \begin{array}{ll}
2 & \mbox{ if } m^\prime=1,\\
0& \mbox{ if } m^\prime=2. 
\end{array} \right. \\
&\left(\Vect{LPL}^{(1)}_{\Ltt} ( \varphi_{i,\bullet};\ParVT)\right)_2^{m_1,m_2} = \ism[\Ltt\times \mSp]{ v_2^{m_1,m_2}(x^m|\varphi_{i,\bullet} e^{-\VIPar{x^m}{\varphi_{i,\bullet} }{\ParVT} }} -
\left\{ \begin{array}{ll}
2 & \mbox{ if } m_1=m_2=1 \\
&\mbox{ and }r_i<D\\
0& \mbox{ otherwise.}
\end{array} \right. \\
\end{align*}
Let $\Vect{x} \in \RR[s]\setminus \{0\}$ , then from previous computations
\begin{equation} \label{eq-hyp0}
\tr{\Vect{x}} \left( \Vect{R}_{\infty,\Ltt}(\varphi_{i,+}; \Vect{h},\ParVT) 
-\Vect{R}_{\infty,\Ltt}(\varphi_{i,-}; \Vect{h},\ParVT) 
\right)= 2\kappa_i x_i+ f(\Vect{x},\varphi_{i,+},\varphi_{i,-},\Vect{h}).
\end{equation}
By using a continuity argument, one may prove that for every $\varepsilon>0$ there exists $\eta>0$ such that $|f((\Vect{x},\varphi_{i,+},\varphi_{i,-},\Vect{h})|\leq \varepsilon $. 
Therefore, assuming that the l.h.s. of~(\ref{eq-hyp0}) equals 0 leads to $x_i=0$ for $i=1,\ldots,s$.

\subsection{Area-interaction point process}

We fix for simplicity $d=2$, though the proofs may be extended easily to higher dimensions.

\subsubsection{Proof that $\Mat{\Sigma}_1(\ParVT)$ is positive-definite for the area-interaction model}
From Proposition~\ref{prop-sdp}, the proof reduces to check  Assumption \textbf{[PD]} with $h=e^V$,  $\overline\delta=D$, $B=\emptyset$ and
\begin{itemize}
\item[$(i)$]$\Vect{Y}_\Lambda(\varphi;\ParVT)=I_\Lambda(\varphi;e^V,\ParVT)$,
\item[$(ii)$]$\Vect{Y}_\Lambda(\varphi;\ParVT)=R_{\infty,\Lambda}(\varphi;e^V,\ParVT)$.
\end{itemize}

Again $(i)$ is ensured by Proposition~\ref{prop-lInn} since this model satisfies \textbf{[Exp]}.

$(ii)$ Let us consider for some $\eta,\omega>0$ the two following events:

\begin{eqnarray*}
A_1(\eta,\omega) &:=& \left\{ \varphi\in \overline{\Omega}: \varphi(\Delta_0(\overline{\delta}))=2, \varphi(\mathcal{B}((0,0),\eta))=1,\varphi(\mathcal{B}((0,\omega),\eta))=1\right\} \\
A_2(\eta,\omega) &:=& \left\{ \varphi\in \overline{\Omega}: \varphi(\Delta_0(\overline{\delta}))=3, \varphi(\mathcal{B}((0,0),\eta))=1,\varphi(\mathcal{B}((0,\omega),\eta))=2\right\} 
\end{eqnarray*}
Fix $\eta,\omega$, let $\varphi_j \in A_j(\eta,\omega)$ and denote by $\widetilde{e^V}(\varphi):=\sum_{x\in \varphi_{\Ltt}}e^{V(x|\varphi\setminus x;\ParVT)}$
\begin{eqnarray*}
I_{\Ltt}(\varphi_j;e^V,\ParVT) &=& |\Ltt| - \widetilde{e^V}(\varphi_j).
\end{eqnarray*} 
When $\eta\to 0$, 
$$
\widetilde{e^V}(\varphi_1) \to 2 e^{\theta_1^\star+ \theta_2^\star g(\omega)} 
\quad \mbox{ and } \quad 
\widetilde{e^V}(\varphi_2) \to 2e^{\theta_1^\star} +  e^{\theta_1^\star+ \theta_2^\star g(\omega)} 
$$
where $g(\omega):=| \mathcal{B}((0,0),R)\cup\mathcal{B}((0,\omega),R)|-| \mathcal{B}((0,0),R)|$.
Moreover, by denoting $\widetilde{v}_2(\varphi) = \sum_{x\in \varphi} v_2(x|\varphi\setminus x)$
\begin{eqnarray*}
\left(\Vect{LPL}^{(1)}_{\Ltt} (\varphi_j ; \ParVT)\right)_1 &=& \int_{\Ltt} e^{-\VIPar{x}{\varphi_j}{\ParVT}}dx - \left\{
\begin{array}{ll}
2 & \mbox{ if } j=1 \\
3 & \mbox{ if } j=2 \\
\end{array} \right. \\
\left(\Vect{LPL}^{(1)}_{\Ltt} (\varphi_j ; \ParVT)\right)_2 &=& \int_{\Ltt} v_2(x|\varphi_j)e^{-\VIPar{x}{\varphi_j}{\ParVT}}dx - \widetilde{v}_2(\varphi_j).
\end{eqnarray*}
Again, when $\eta\to 0$, one may note that for $k=1,2$
$$
\int_{\Ltt} v_k(x|\varphi_1) e^{-\VIPar{x}{\varphi_1}{\ParVT}}dx - \int_{\Ltt} v_k(x|\varphi_2) e^{-\VIPar{x}{\varphi_2}{\ParVT}}dx \to 0
$$
and 
$$
\widetilde{v}_2(\varphi_1) \to 2 g(\omega) 
\quad \mbox{ and } \quad 
\widetilde{v}_2(\varphi_2) \to g(\omega)
$$
These computations lead to 
\begin{eqnarray*}
R_{\infty,\Ltt}(\varphi_1;e^V,\ParVT) -R_{\infty,\Ltt}(\varphi_2;e^V,\ParVT) &=& 2e^{\theta_1^\star} - e^{\theta_1^\star+\theta_2^\star g(\omega)} -\left(\Vect{W}(e^V,\ParVT)\right)_1 +g(\omega)\left(\Vect{W}(e^V,\ParVT)\right)_2 \\&&+ f(\varphi_1,\varphi_2,\Vect{W}),
\end{eqnarray*}
where the function $f$ is such that for all $\varepsilon>0$, there exists $\eta$ small enough such that {$|f(\varphi_1,\varphi_2,\Vect{W})|\leq \varepsilon$}.
Let $\varphi_j \in A_j(\eta,0)$, then, since $g(0)=0$, assuming that the l.h.s. of the previous equation equals 0 leads to $\left(\Vect{W}(e^V,\ParVT)\right)_1=e^{\theta_1^\star}$. Now, let $\omega>0$ and again assume that $R_{\infty,\Ltt}(\varphi_1;e^V,\ParVT) =R_{\infty,\Ltt}(\varphi_2;e^V,\ParVT)$, we therefore obtain (by the continuity argument)
$$
\left(\Vect{W}(e^V,\ParVT)\right)_2 = \frac{e^{\theta_1^\star+\theta_2^\star g(\omega)} - e^{\theta_1^\star}}{g(\omega)}.
$$
But $\left(\Vect{W}(e^V,\ParVT)\right)_2$ is a constant and so cannot depend on $\omega$. Therefore one of the assumptions made before is untrue, which proves \textbf{[PD]}.

\subsubsection{Proof that $\Mat{\Sigma}_2(\ParVT)$ is positive-definite for the area-interaction model}\label{sec-DP1fwk2Str}
From Proposition~\ref{prop-sdp}, it suffices to check Assumption \textbf{[PD]}  with
$$\Vect{Y}_\Lambda(\varphi;\ParVT)=\Vect{R}_{\infty,\Lambda}(\varphi;\Vect{h};\ParVT ),$$
where, for all $j=1,\dots,s$, $\left(\Vect{R}_{\infty,\Lambda}(\varphi;\Vect{h};\ParVT)\right)_j=R_{\infty,\Lambda}(\varphi;h_j,\ParVT)$, $h_j$ is the test function given by $h_j(x^m,\varphi;\ParV)= \mathbf{1}_{[0,r_j]}(d(x^m,\varphi)) e^{V\left(x^m|\varphi;\ParV\right)}$, and where, again, we choose $\overline\delta=D$ and $B=\emptyset$.

The proof is quite similar to the one proposed for the 2-type marked Strauss point process (see \ref{sec-DP1fwk2Str}). Let $0<r_1<\ldots<r_s<+\infty$. Let us also assume that $r_i\neq D$ for $i=1,\ldots,s$
\begin{eqnarray*}
A_{i,-}(\eta) &=& \Big\{ \varphi \in \cSptt : \varphi({\Dom[0](\Dtt)})=2, \varphi\left( \mathcal{B}\left( (0,0),\frac{\eta}4 \right) \right)=1, \varphi\left( \mathcal{B}\left( (r_i-\frac\eta2,0),\frac{\eta}4 \right) \right)=1\Big\},  \\
A_{i,+}(\eta) &=& \Big\{ \varphi \in \cSptt : \varphi({\Dom[0](\Dtt)})=2, \varphi\left( \mathcal{B}\left( (0,0),\frac{\eta}4 \right) \right)=1, \varphi\left( \mathcal{B}\left( (r_i+\frac\eta2,0),\frac{\eta}4 \right) \right)=1\Big\}.
\end{eqnarray*}
Let $i,j\in \{1,\ldots,s\}$ and $k\in \{1,2\}$, let $\varphi_{i,-} \in A_{i,-}$ and $\varphi_{i,+} \in A_{i,+}$, then
\begin{eqnarray*}
I_{\Ltt}(\varphi_{i,-};h_j,\ParVT) &=& \int_{\Ltt} h_j(x,\varphi_{i,-};\ParVT) e^{-\VIPar{x}{\varphi_{i,-}}{\ParVT}}dx - 
\left\{ \begin{array}{ll}
\widetilde{e^V}(\varphi_{i,-}) & \mbox{ if } i\leq j \\
0 & \mbox{ otherwise.}
\end{array} \right. \\
I_{\Ltt}(\varphi_{i,+};h_j,\ParVT) &=& \int_{\Ltt} h_j(x,\varphi_{i,+};\ParVT) e^{-\VIPar{x}{\varphi_{i,+}}{\ParVT}}dx - 
\left\{ \begin{array}{ll}
\widetilde{e^V}(\varphi_{i,+}) & \mbox{ if } i< j \\
0 & \mbox{ otherwise.}
\end{array} \right. \\
\left( \Vect{LPL}^{(1)}_{\Ltt}(\varphi_{i;\bullet};\ParVT)\right)_k &=& \int_{\Ltt} v_k(x|\varphi_{i,\bullet})e^{-\VIPar{x}{\varphi_{i,\bullet}}{\ParVT}} dx - \sum_{x\in \varphi_{i,\bullet}} v_k(x|\varphi_{i,\bullet}\setminus x),
\end{eqnarray*}
for $\bullet=-,+$. It is expected that for small $\eta$, $\left( \Vect{LPL}^{(1)}_{\Ltt}(\varphi_{i,-};\ParVT)\right)_k\simeq \left( \Vect{LPL}^{(1)}_{\Ltt}(\varphi_{i,+};\ParVT)\right)_k$ and $\widetilde{e^V}(\varphi_{i,-})\simeq \widetilde{e^V}(\varphi_{i,+})\simeq \kappa_i:=2 e^{\theta_1^\star+\theta_2^\star|\mathcal{B}(0,R)\cup\mathcal{B}(r_i,R)|}$. Let $\Vect{x} \in \RR[s]\setminus \{0\}$ , then from previous computations
\begin{equation} \label{eq-hyp0}
\tr{\Vect{x}} \left( \Vect{R}_{\infty,\Ltt}(\varphi_{i,+}; \Vect{h},\ParVT) 
-\Vect{R}_{\infty,\Ltt}(\varphi_{i,-}; \Vect{h},\ParVT) 
\right)= 2\kappa_i x_i+ f(\Vect{x},\varphi_{i,+},\varphi_{i,-},\Vect{h})
\end{equation}
where for every $\varepsilon>0$ there exists $\eta>0$ such that $|f((\Vect{x},\varphi_{i,+},\varphi_{i,-},\Vect{h})|\leq \varepsilon $. 
Therefore, assuming that the l.h.s. of~(\ref{eq-hyp0}) equals 0 leads to $x_i=0$ for $i=1,\ldots,s$.


\bibliographystyle{plainnat.bst}

\bibliography{residus}

\end{document}